\documentclass[preprint, 10pt]{imsart}

\RequirePackage[OT1]{fontenc}
\RequirePackage{amsthm,amsmath}
\RequirePackage[authoryear,round]{natbib}
\RequirePackage[colorlinks,citecolor=blue,urlcolor=blue]{hyperref}

\usepackage{bibentry}
\usepackage{float}
\usepackage{amsmath,amssymb,epsfig}
\usepackage{graphics}
\usepackage{tabularx}
\usepackage{framed}
\usepackage[mathscr]{euscript}
\usepackage{imsart}

\usepackage{epstopdf}

\usepackage[makeroom]{cancel}

\usepackage{soul,xcolor}

\startlocaldefs
\numberwithin{equation}{section}
\theoremstyle{plain}

\long\def\comment#1{}

\newtheorem{theorem}{Theorem}

\newtheorem{lemma}{Lemma}

\theoremstyle{definition}
\newtheorem{definition}{Definition}

\newcommand{\citen}{\cite}

\renewcommand{\H}{{\mathsf H}}

\newcommand{\eps}{\varepsilon}

\newcommand{\be}{\begin{eqnarray}}
\newcommand{\ee}{\end{eqnarray}}

\newcommand{\E}{\mathscr E}

\newcommand{\thetahat}{\hat \theta}

\newcommand{\ba}{\begin{array}}
\newcommand{\ea}{\end{array}}
\newcommand{\bs}{\begin{align}\begin{split}\nonumber}
\newcommand{\bsnumber}{\begin{align}\begin{split}}
\newcommand{\es}{\end{split}\end{align}}

\newcommand{\n}{n}

\renewcommand{\(}{\left(}
\renewcommand{\)}{\right)}
\renewcommand{\[}{\left[}
\renewcommand{\]}{\right]}
\renewcommand{\hat}{\widehat}

\newcommand{\Ep}{{\mathrm{E}}}
\newcommand{\ctest}{c_{\text{test}}}
\newcommand{\N}{\mathrm N}

\newcommand{\En}{{\mathbb E_n}}

\renewcommand{\Pr}{{\mathrm{P}}}

\def\x{{x}}
\def\supp{{\rm support}}

\renewcommand{\hat}{\widehat}
\renewcommand{\leq}{\leqslant}
\renewcommand{\geq}{\geqslant}
\renewcommand{\supp}{\text{supp}}

\endlocaldefs

\begin{document}

\setstcolor{red}
\setul{}{1pt}

\begin{frontmatter}

\title{\vspace*{-8mm} \large \textcolor{black}{Analysis of} Testing-Based Forward Model Selection \\ }
\runtitle{ATBFMS}

\begin{aug}
\author{\fnms{Damian} \snm{Kozbur}\ead[label=e1]{damian.kozbur@econ.uzh.ch.}}

\thankstext{t1}{ 
 First version: November 1, 2015.  This version \today.  I gratefully acknowledge helpful discussions with Christian Hansen, Tim Conley, Kelly Reeve, Dan Zou, Nicolai Meinshausen, Marloes Maathius, Rahul Mazumder, Ryan Tibshirani, Trevor Hastie, Martin Schonger, Pietro Biroli, Michael Wolf, Lorenzo Casaburi, Hannes Schwandt, Ralph Ossa, Rainer Winkelmann, attendants at the ETH Z\"urich Seminar f\"ur Statistik Research Seminar, attendants at the Center for Law and Economics Internal Seminar, attendants at the Toulouse School of Economics Econometrics Seminar, research assistants Vincent Lohmann and Alexandre Jenni for assisting in several stages of verification of formal arguments, Yuming Pan for assistance with computation, as well as financial support from UZH and from the ETH Fellowship program. The content of this paper draws from two separate working papers posted on ArXiv (https://arxiv.org/abs/1702.01000 (\cite{DK:FSEL}) and https://arxiv.org/abs/1512.02666 (the current paper)).  The projects have been merged in preparation for the publication process. A mechanical statement of a special case of TBFMS appears in \cite{DK:TBFMS:PandP}, which neither claims nor derives any theoretical results.  This paper previously had the title ``Testing-Based Forward Model Selection,'' which is also the title of  \cite{DK:TBFMS:PandP}.}

\runauthor{Damian Kozbur}

\affiliation{University of Z\"urich }

\address{ University of Z\"urich \\ Department of Economics \\ Sch\"onberggasse 1, 8001 Z\"urich  \\
\printead{e1}\\
}

\end{aug}

\

\begin{abstract} 
This paper \textcolor{black}{ analyzes} a procedure called {Testing-Based Forward Model Selection} (TBFMS) in linear regression problems. This procedure inductively selects covariates that add predictive power into a working statistical model before estimating a final regression.  The criterion for deciding which covariate to include next and when to stop including covariates is derived from a profile of traditional statistical hypothesis tests.  This paper proves probabilistic bounds, which depend on the quality of the tests, for prediction error and the number of selected covariates.  As an example, the bounds are then specialized to a case with heteroskedastic data, with tests constructed with the help of Huber-Eicker-White standard errors.  Under the assumed regularity conditions, these tests lead to estimation convergence rates matching other common high-dimensional estimators including Lasso. %TBFMS performance is compared to Lasso and Post-Lasso in simulation studies.   TBFMS is then analyzed as a component into larger post-model selection estimation problems for structural economic parameters.  
\end{abstract}

\begin{keyword}[class=MSC]
\kwd[]{62J05, 62J07, 62L12}
\end{keyword}

\begin{keyword}[class=JEL]
{\hspace{-1mm}\textbf{JEL subject classifications:} C55}
\end{keyword}

\begin{keyword}
\kwd{model selection, forward regression, sparsity, hypothesis testing} 
\end{keyword}

\end{frontmatter}
\thispagestyle{empty}
\pagebreak
\setcounter{page}{1}
\section{Introduction}
%\doublespacing

This paper analyzes a procedure called \textit{Testing-Based Forward Model Selection} (TBFMS) for high-dimensional econometric problems, which are characterized by settings in which the number of observed characteristics per observation in the data is large.\footnote{High-dimensional data may arise in several ways---data may be intrinsically high-dimensional with many characteristics per observation, or alternatively, researchers may obtain a large final set of covariates through forming interactions and transformations of underlying covariates.}    High-dimensional econometrics is a leading area of current research because of recent rapid growth in data availability and computing capacity, coupled with the important need to extract as much useful information from data in a way that allows precise and rigorous testing of scientific hypotheses.  %Working within the flexibility of a high-dimensional framework allows researchers to fully exploit richer data sets both in prediction problems and in structural inference problems.

The primary settings of this paper are high-dimensional sparse linear regression models, in which the number of covariates is allowed to be commensurate with or exceed the sample size.  A key challenge with a high-dimensional data set is that estimation requires dimension reduction or regularization to avoid statistical overfitting.   A \textit{sparsity} assumption imposes that the regression function relating the outcome and the covariates can be approximated by a regression of the outcome on a small, ex ante unknown subset of covariates.  Under sparsity, there are several consistent estimation procedures (further reviewed below) that work by enforcing that the estimated regression function be sparse or small under an appropriate norm.

An appealing class of techniques for high-dimensional regression problems are \textit{greedy} algorithms.  These are procedures that inductively select individual covariates into a working model (i.e., a collection of covariates) until a stopping criterion is met.  A linear regression restricted to the final selected model is then estimated.  A leading example is \textit{Simple Forward Selection},\footnote{Simple Forward Selection is not standard nomenclature, but is used here in order to have a parallel language with Testing-Based Forward Model Selection.  The literature is varied and uses several names including Forward Regression and Forward Stepwise Regression.  ``Model'' is used in the name of TBFMS to avoid confusion with sample selection problems common in econometrics.} which chooses the covariate that gives the highest increase of in-sample R-squared above the previous working model.  This class of techniques is widely used because they are intuitive and simple to implement.  Such methods in the statistics literature date back to at least \cite{efroymson1966stepwise}.

In practice, deciding which covariate gives the best additional predictive power relative to a working model is complicated by the fact that outcomes are observed with noise or are partly idiosyncratic.  For example, in linear regression, a variable associated with a positive increment of in-sample R-squared upon inclusion may not add any predictive power out-of-sample.  Statistical hypothesis tests offer one way to determine whether a variable of interest is likely to improve out-of-sample predictions.  Furthermore, in many econometric and statistical applications, the classical assumption of independent and identically distributed data is not always appropriate.  %One example of this is the presence of heteroskedastic disturbances. %In such settings, higher R-squared resulting from inclusion of one variable relative to another need not be a signal that the first variable is a better choice.  More generally, model selection procedures tailored to the classical assumptions may have inferior performance when applied to more realistic data-generating processes.  
The availability of hypothesis tests for diverse classes of problems and settings motivates the introduction of a testing-based strategy.  Mechanically, TBFMS begins with an empty model.  The procedure then tests whether any covariates provide additional predictive capability in the population.   The selection stops when no tests return a significant covariate.  Selection into the model is then based on the largest value of an associated test statistic.   Note that in this context, the hypothesis tests are solely serving a role as assisting model selection, not ex post inference.

There are several earlier analyses of Simple Forward Selection.\footnote{TBFMS using different tests than proposed here is natively programmed in some statistical software, including SPSS, but is not previously formally  justified \textcolor{black}{in high-dimensional settings}.}     \cite{Wang:UltraHighForwardReg}  gives bounds on the performance and number of selected covariates under a $\beta$-min condition which requires the minimum magnitude of non-zero coefficients to be suitably bounded away from zero.  \cite{Zhang:GreedyLeastSquares} and \cite{Tropp:Greed} prove performance bounds for greedy algorithms under a strong irrepresentability condition, which restricts the empirical covariance matrix of the predictors.   \cite{Submodular:Spectral} prove bounds on the relative performance in population $R$-squared of a forward selection based model (relative to infeasible $R$-squared) when the number of variables allowed for selection is fixed.     In addition to Simple Forward Selection, there are several related procedures in which estimation is done in stages.  These include a method that is not strictly greedy called \textit{Forward-Backward Selection}, which proceeds similarly to Simple Forward Selection but allows previously selected covariates to be discarded from the working model at certain steps (see \cite{Zhang:GreedyLeastSquares}). 

% Another related class of methods are called boosting methods.  Boosting methods inductively select covariates predictive of a linear combination of estimated residuals and the outcome at each step (among many other references, see \cite{BuhlmannBoosting2006}, see also \cite{Martin:Boosting} for additional results and applications in econometrics.)   A recent related paper, \cite{pesaran:testing}, considers a different iterative model selection procedure which also involves using hypothesis tests in high-dimensional sparse linear models.  In \cite{pesaran:testing}, in the first iteration, a marginal regression of the outcome on each potential covariate is run.   Once all marginal regressions are run, all significant covariates are included into a working model.   Each subsequent iteration in \cite{pesaran:testing} works similarly.

As a preliminary, this paper proves new bounds on the predictive performance and number of selected covariates for Simple Forward Selection.  The conditions required here are weaker that those used in \cite{Zhang:GreedyLeastSquares} and \cite{Tropp:Greed} and impose no $\beta$-min restrictions or irrepresentability.  
The convergence rates here are most similar to the analysis of a Forward-Backward Selection in \cite{TongZhang:ForwardBackward}, but require markedly different analysis since there is no chance to correct ``over-selection mistakes.''  %As a part of the analysis in this paper, it is shown that that mistakes, suitably defined, cannot happen too often.  
 
%The analysis of Simple Forward Selection lays the groundwork for deriving statistical performance bounds for TBFMS.   
This paper then gives performance bounds for TBFMS which depend directly on the quality of the profile of tests considered, as measured by five constants which characterize size and power.  
The abstract results for TBFMS are used to derive asymptotic bounds for various sequences of data-generating processes. As an example, concrete tests for heteroskedastic data constructed from Huber-Eicker-White standard errors are used to construct $t$-tests and explicit rates of convergence are calculated.

There are many other sensible approaches to high-dimensional estimation.  An important and common approach to generic high-dimensional estimation problems is the Lasso.  
The Lasso minimizes a least squares criterion augmented with a penalty proportional to the $\ell_1$ norm of the coefficient vector.  %This approach favors a model with good in-sample prediction while still placing high value on parsimony (the structure of the objective sets many coefficients identically to zero).  The Post-Lasso refits based on a least squares objective function on the selected model. 
For theoretical and simulation results for Lasso, see \cite{FF:1993} \cite{T1996}, \cite{elements:book} \cite{CandesTao2007} \cite{BaiNg2008}, \cite{BickelRitovTsybakov2009}, \cite{horowitz:lasso},  \cite{BuhlmannGeer2011}, among many more.  Other related methods include boosting  (see \cite{Freund+Schapire:1996},  \cite{BuhlmannBoosting2006},  \cite{Martin:Boosting}), Least Angle Regression (see \cite{EfronHastieJohnstoneTibshirani2004} ), Post-Lasso (see \cite{BellChenChernHans:nonGauss}), and many others.   A recent related paper, \cite{pesaran:testing}, considers a different iterative model selection procedure which also involves using hypothesis tests.  In \cite{pesaran:testing}, in the first iteration, a marginal regression of the outcome on each potential covariate is run.   Once all marginal regressions are run, all significant covariates are included into a working model.   Each subsequent iteration works similarly.  

The asymptotic estimation rates calculated here for TBFMS, applied to a constructed profile of tests for heteroskedastic data, match those standard for Lasso and Post-Lasso.  
%Though all alternative procedures have benefits in different settings, TBFMS can be applied to a large set of problems for which there is a reliable testing procedure for determining whether any particular variable (or set of variables or new parameters) adds predictive power.  
Relative to the analysis of asymptotic properties of Lasso and related high-dimensional estimation techniques, analysis of TBFMS is complicated by the fact that the procedure is not the optimizer of a simple objective function.  As a result, the theory also departs from the literature on m-estimation in a fundamental way.  

A recent paper, \cite{Hastie:Tibshirani:Tibshirani:ExtendedComparisons}, performs a systematic simulation analysis of statistical and computational performance of simple forward selection as well as a few additional estimators, including Lasso and best subset selection; see \cite{best:subset:MIP}.  The paper reports that in regression models with higher signal-to-noise ratios, forward selection performs favorably relative to Lasso; a finding consistent across many simulation settings.  

This paper complements recent literature on sequential testing (see \cite{GSell:SeqTesting}, \cite{Li:Barber:AccumTestsFDR}, \cite{SelSeqModSel2}, \cite{SelSeqModSel1}).  Sequential testing considers hypothesis testing in stages, in which tests in later stages can depend on testing outcomes in earlier stages.    In various settings, properties like family-wise error rates of proposed testing procedures can be controlled over such sequences of hypothesis tests.  % In all cases, the authors note that the testing procedures are complementary to forward model selection problems as they guide which variables should be selected and offer principled stopping rules.   
While the current paper focuses on statistical properties of estimates after TBFMS given properties of the implemented tests, future work may potentially combine the two types of problems.

%After developing several theoretical bounds, a simulation study illustrates the relative performance of TBFMS to Lasso and Post-Lasso regression.  The simulation study shows that there are data-generating processes under which forward selection outperforms Lasso regression in terms of prediction and estimation error.  

In economic applications, models learned using formal model selection are often used in subsequent estimation steps, with the final goal of learning a structural parameter of interest.  One example is the selection of instrumental variables for later use in a first-stage regression (see \cite{BellChenChernHans:nonGauss}).  Another example is the selection of a conditioning set to properly control for omitted variables bias when there are many control variables (see  \citen{ZhangZhang:CI}, \citen{BCH-PLM}, \citen{vdGBRD:AsymptoticConfidenceSets}, and \citen{JM:ConfidenceIntervals}).  Bounds about the quality of the selected model are used to derive results about the quality of post-model selection estimation and to guide subsequent inference.  Such applications require a model selection procedure with hybrid objectives:  (1) produce a good fit, and (2) return a sparse set of variables.  %\footnote{In simulation studies reported in \cite{TU} and below, cross-validation consistently returns models with a relatively higher number of parameters (relative to e.g. the Lasso penalization procedure suggested in \cite{BellChenChernHans:nonGauss}) .  Additional simulation results documenting size distortion using cross-validated Lasso for post-model selection procedures (post-double selection and first stage instrumental variable selection) can be found below in this paper, and in the reference \cite{TU}.  }  
This paper addresses both objectives by providing sparsity and fit bounds for TBFMS.

% On the other hand, Lasso has good estimation properties relative to forward selection in low-signal-to-noise settings.  Informal intuition for this behavior is that it results from Lasso selecting a higher number of covariates while simultaneously shrinking the corresponding estimated regression coefficients.  For input into post-double selection, such use of Lasso for model selection can cause inferential size-distortion:  simulation evidence on the effects of over-selecting on size distortion in post-double selection appear in \cite{BCHK:Panel}.  For input into a first stage model selection in an instrumental variable regression, size distortion is also a problem associated to over selection. A theoretical explanation of this phenomenon in the instrumental variables setting is given in \citen{RJIVE}.  

In terms of computing, one fast implementation of forward selection depends on what is sometimes referred to as a ``guided QR decomposition.''  Formally, simple forward selection can be computed in $O(npk)$ flops, with $n$ being sample size, $p$ being number of covariates, and $k$ being number of steps (see for example  \cite{Hastie:Tibshirani:Tibshirani:ExtendedComparisons}), and requires the storage of the QR decomposition of at most $k$ variables.  The version of TBFMS presented in the paper for data with heteroskedastic disturbances can be computed with the same order of time and storage requirements.\footnote{A modification of the least angle regression algorithm can be made to implement a similarly efficient computation of Lasso (see \cite{EfronHastieJohnstoneTibshirani2004}, though this may require potentially more iterations, as covariates can multiple times both enter and exit a suitably defined active set which terminates as the selected set).}

\section{Precursor: Sharp Convergence Rates for Simple Forward Selection without $\beta$-min or Irrepresentability Conditions}

This section proves a precursory result about Simple Forward Selection which is new in the high-dimensional econometrics and statistics literature.  The procedure is defined formally below and is similar to TBFMS, but uses a single threshold rather than a profile of hypothesis tests in determining the selection of covariates.   The framework set out in this section is also helpful in terms of outlining minimal structure needed to facilitate the method of analysis in the formal arguments that follow.%\footnote{This section draws material from the draft \cite{DK:FSEL}, which was originally a separate project and posted on ArXiv, but is now merged into the current paper in preparation for the publication process.}  %Despite the fact that the high-dimensional econometrics and statistics literatures are very robust by now, this section gives new analysis for a fundamental problem. 

\subsection{Framework}

A realization of data of sample size $n$ is given by $\mathscr D_n = \{(x_i,y_i)\}_{i=1}^n$ and is generated by a joint distribution $\Pr$.  The data consist of a set of covariates $x_i\in \mathbb R^p,$ as well as outcome variables $y_i \in \mathbb R$ for each observation $i=1,...,n$.  The data satisfy  
$$y_i = x_i'\theta_0 + \eps_i$$
for some unknown parameter of interest $\theta_0 \in \mathbb R^p$ and unobserved disturbance terms $\eps_i \in \mathbb R$.  
%The covariates $x_i$ are normalized so that $\En[x_{ij}]=0$ and $ \En[x_{ij}^2]=1$ for every $j=1,...,p$, where $\En[ \hspace{.6mm} \cdot \hspace{.6mm} ] = \frac{1}{n}\sum_{i=1}^n ( \hspace{.2mm} \cdot \hspace{.2mm}  )$ denotes empirical expectation.  
The parameter $\theta_0$ is sparse in the sense that the set of non-zero components of $\theta_0$, denoted $S_0=\text{supp}(\theta_0)$, has cardinality $s_0<n$.  \textcolor{black}{(Below, exact orthogonality between $\eps_i$ and $x_i$ will not be required; rather a notion of approximate orthogonality will be used.  As a result, the framework also handles notions of approximate sparsity). }  %The interest in this paper is in studying estimates for $x_i'\theta_0$ for $i=1,...,n$.  

Define an empirical loss function $\ell(\theta)$ $$\ell(\theta) = \En[ (y_i - x_i'\theta)^2],$$
where $\En[ \hspace{.6mm} \cdot \hspace{.6mm} ] = \frac{1}{n}\sum_{i=1}^n ( \hspace{.2mm} \cdot \hspace{.2mm}  )$ denotes empirical expectation.   Note that $\ell(\theta)$ depends on $\mathscr D_n$, but this dependence is suppressed from the notation.  
Define also for subsets $S \subseteq\{1,...,p\}$ $$\ell(S) = \min_{\theta: \text{supp}(\theta) \subseteq S} \ell(\theta).$$The estimation strategy proceeds by first searching for a sparse subset $\hat S \subseteq \{1,...,p\}$, with cardinality $\hat s$,  that assumes a small value of
$\ell(S)$,
followed by estimating $\theta_0$ with least squares via
$$\hat \theta \in \text{arg} \min_{\theta:\text{supp}(\theta) \subseteq \hat S} \ell(\theta).$$
This gives the construction of the estimates $ x_i'\hat \theta$ for $i=1,...,n$.
%The paper provides bounds for the prediction error norm defined by $$\En [(x_i' \theta_0 - x_i'\hat \theta )^2]^{1/2}. $$

The set $\hat S$ is selected as follows.    For any $S$ define the incremental loss from the $j$th covariate by
$$\Delta_j \ell(S) = \ell(S \cup \{j\} ) -   \ell(S).$$ 

\noindent Consider the greedy algorithm, which inductively selects the $j$th covariate to enter a working model if $-\Delta_j\ell(S)$ exceeds a threshold $t$:  $$-\Delta_j\ell(S) > t$$ and $-\Delta_j\ell(S) \geq -\Delta_k\ell(S)$ for each $k \neq j$.    The threshold $t$ is chosen by the user; it is the only tuning parameter required.  This procedure is summarized formally here.

\

\noindent \textbf{Algorithm 1.} \textit{Simple Forward Regression}

\noindent \textit{Initialize}  \noindent Set $\hat S =  \varnothing$ 

\noindent \textit{For $1 \leq k \leq p$} 

\textit{If} $-\Delta_j\ell(S) > t$ for some $j\in \{1,...,p\} \setminus \hat S$ 

\ \ \ \ \textit{Set} $\hat j \in \text{arg} \max  \left \{-\Delta_j\ell(S) : -\Delta_j\ell(S) > t \right \}$ 

\ \ \ \ \textit{Update }  $\hat S = \hat S \cup \{\hat j\}$ 

\textit{Else}

\ \ \ \  \textit{Break}

\noindent \textit{Set} $\hat \theta \in \arg \min_{\theta:\supp(\theta) \subset \hat S} \ell({\theta}).$

\subsection{Formal Analysis}

In order to analyze Algorithm 1 and state the first theorem, a few more definitions are convenient.  Define the empirical Gram matrix $G$ by $G = \En[x_ix_i']$. Let $\varphi_{\min}(s)(G)$ denote the minimal $s$-sparse eigenvalues given by $$\varphi_{\min}(s)(G) = \min_{S \subseteq \{1,...,p\}: |S| \leq s} \lambda_{\min}(G_{S})$$ where $G_{S}$ is the principal submatrix of $G$ corresponding to the component set $S$.   The maximal sparse eigenvalues $\varphi_{\max}(s)(G)$ are defined analogously.
Let
%$$  \mathscr C_1 = 2\| \En[\eps_i x_i'] \|_\infty \sqrt{\hat s + s_0} \varphi_{\min}(\hat s + s_0)(G)^{-1} + \sqrt{s_0t\varphi_{\min}(\hat s + s_0)(G)^{-1}}.$$
$$  c_{\text{F}}(\hat s) = (\hat s + s_0)^{1/2}  \varphi_{\min}(\hat s + s_0)(G)^{\textcolor{black}{-1/2}} \[2\| \En[x_i \eps_i] \|_\infty +t^{1/2} \].$$

Finally, for each positive integer $m$, let 
$$
\textcolor{black}{ c_{\text{F}}'(m)  = 80  \times  \varphi_{\min}(m+s_0)(G)^{-4}.}
$$ %
%The  above quantities are useful for displaying results in Theorem 1. Slightly tighter but messier usable quantities than $ c_{\text{fsel}}(\hat s) $ and $c_{\text{fsel}}'(m) $ are derived in the proof.  %Note also that $\mathscr C_1$ depends on $\hat s$.

\begin{theorem}

Consider a data set $\mathscr D_n$ of fixed sample size $n$ with parameter $\theta_0$.  Suppose the normalizations $\En[x_{ij}^2] =1$ hold for each $j \leq p$.  Then under Algorithm 1 with threshold $t$,

$$
\En [(x_i'\theta_0 -\x_{ i}'\hat \theta)^2]^{1/2} \leq   c_{\text{\em F}}(\hat s) . $$
In addition, for every integer $m \geq 0$ with  \textcolor{black}{$m \leq |\hat S \setminus S_0|$} such that  ${t^{1/2}} \geq  2\varphi_{min}(m+s_0)(G)^{-1} \| \En[x_i \eps_i]\|_{\infty}$, it holds that 
$$ \ m \leq   c_{\text{\em F}}'(m)  s_0.$$

\end{theorem}

The above theorem calculates explicit finite sample constants bounding the prediction error norm.  The second statement is a tool for bounding the number of selected covariates.   In particular, setting $m^* = \min \{ m: m> c_{\text{F}}'(m)s_0 \}$ implies that $$\hat s < m^*+s_0$$ provided that the condition on $m^*$ given by $t^{1/2} \geq 2 \varphi_{\min}(m^* + s_0)(G)^{-1} \| \En[x_i \eps_i]\|_\infty$ is met.

The statement in Theorem 1 gives finite sample bounds which are completely deterministic in the sense that they hold for every possible realization of the data.  Furthermore, the proof does not use any random nature of $\mathscr D_n$ at any step.  As a result, the bounds are very general, but it is helpful for interpretation to consider the convergence rates implied by Theorem 1 under asymptotic conditions on $\mathscr D_n$.   Consider a sequence of random data sets $(\mathscr D_n)_{\n \in \mathbb N}$ generated by joint distributions $(\Pr = \Pr_{n})_{n\in \mathbb N}$.   For each $n$, the data again satisfy $y_i = x_i'\theta_0 + \eps_i$. In what follows, the parameters $\theta_0$, the thresholds $t$, distribution $\Pr$, the dimension $p$ of $x_i$, etc.  can all change with $n$.

\

\noindent \textbf{Condition 1} (\textit{Asymptotic Regularity}). The sparsity satisfies $s_0 = o(n)$.  There is a sequence $K_n$ for which $s_0= o(K_n) $ and there is a bound $\varphi_{\min}(K_n)(G)^{-1} = O(1)$ which holds with probability $1 - o(1)$.  The normalizations $\En[x_{ij}^2]=1$ hold a.s. for every $j \leq p$.  The threshold satisfies a bound $t = O({ \log p /n})$.  In addition, ${t^{1/2}} \geq 2 \varphi_{\min}(K_n)(G)^{-1} \| \En[x_i \eps_i]\|_{\infty}$ with probability $1 - o(1)$.\footnote{\textcolor{black}{Formally, for a sequence of random variables $X_n$, the statement ``$X_n = O(1)$ with probability $1 - o(1)$'' is defined as:  ``there is a constant $C$ independent of $n$ such that $\Pr(|X_n| > C) \rightarrow 0$.''}}%\footnote{In Condition 1 and in the rest of the paper, the expression $O$ signifies bounds which do not depend on $n$.  For instance, $t = O(\log p / n)$ means that there exists a constant $C<\infty$ which does not depend on $n$ such that $| t | \leq C$.  The asymptotic $o$ is taken as $n \rightarrow \infty$.}

\

The rates assumed in Condition \textcolor{black}1 reflect typical rates achieved under various possible sets of low-level conditions standard in the literature (ie. \cite{BellChenChernHans:nonGauss}).  Condition 1 asserts three important statements.  The first statement bounds the size of $S_0$ and requires that the sparsity level is small relative to the sample size.
%In practice, it is difficult to measure $s$ and $ c_{\text{sprs}}$ from a dataset and so for application, values for these parameters are often assumed.  
The second statement is a sparse eigenvalue condition useful for proving results about high-dimensional techniques.  In standard regression analysis where the number of covariates is small relative to the sample size, a common assumption used in establishing properties of conventional estimators of $\theta$ is that $G$ has full rank. In the high-dimensional setting,  $G$ will be singular if $p>n$ and may have an ill-behaved inverse even when $p \leq n$. However, good performance of many high-dimensional estimators only requires good behavior of certain moduli of continuity of $G$.   There are multiple formalizations and moduli of continuity that can be considered here; see  \cite{BickelRitovTsybakov2009}. This analysis focuses on a simple eigenvalue condition that was used in \cite{BellChenChernHans:nonGauss}.  Condition 1 could be shown to hold under more primitive conditions by adapting arguments found in \cite{BC-PostLASSO}, which build upon results in \cite{ZhangHuang2006} and \cite{RudelsonVershynin2008}; see also \cite{RudelsonZhou2011}.  Condition 1 is notably weaker than previously used irrepresentability conditions.  Irrepresentability conditions require that for certain sets $S$ and $k \notin S$, letting $x_{iS}$ be the subvector of $x_i$ with components $j\in S$, $\| \En[x_{iS}x_{iS}']^{-1}\En[x_{iS}x_{ik}] \|_1$ is strictly less than 1.  The normalization $\En[x_{ij}^2]=1$ is used to keep exposition concise and can be relaxed (and, e.g.,  is relaxed in Theorem 5).

The final statement in Condition 1 is a regularization condition similar to regularization conditions common in the analysis of Lasso.  The condition requires $ t^{1/2}$ to dominate a multiple of the $\| \En[x_i \eps_i]\|_\infty$.  This condition is stronger than that typically encountered with Lasso, because the multiple depends on the sparse eigenvalues of $G$.  To illustrate why such a condition is useful, let $\check x_{ij}$ denote $x_{ij}$ residualized away from previously selected regressors and renormalized.  Then even if  $\En[x_{ij} \eps_i] < t^{1/2}$,  $\En[\check x_{ij}\eps_i]$ can exceed $t^{1/2}$, resulting in more selections into the model.   Nevertheless, using the multiple $2 \varphi_{\min}(K_n)(G)^{-1}$, which stays bounded with $n$, is sufficient to ensure that $\hat s$ does not grow faster than $s_0$.  Furthermore, this requirement does not implicitly impose a $\beta$-min condition and does not implicitly impose irrepresentability.   \textcolor{black}{The requirements on $t$ can be relaxed if there is additional control on quantities of the form $\En[\check x_{ij}\eps_i]$.}  Relative to analogous Lasso bounds in \cite{BellChenChernHans:nonGauss}, Theorem 1 does not involve maximal sparse eigenvalues. This may become relevant if the components of $x_i$ arise from factor model structures.

From a practical standpoint, Condition 1 does, however, implicitly require the user to know more about the design of the data in choosing an appropriate $t$.  Choosing feasible thresholds which satisfy a similar condition to Condition 1 is considered in the next section, where \textcolor{black}{analysis of} TBFMS is developed.

\begin{theorem}
Consider a sequence of data sets $\mathscr D_n$ indexed by $n$ with parameters $\theta_0$ and threshold $t$ which satisfy Condition 1.  Suppose $\hat \theta$ is obtained by Algorithm 1.  Then there are bounds
$$
\En [(x_i'\theta_0 -\x_{ i}'\hat \theta)^2]^{1/2} = O\(\sqrt{ \frac{s_0 \log p }{n} } \ \),$$
$$ \ \hat s \leq O(s_0),$$
which hold with probability $1 - o(1)$ as $n \rightarrow \infty$. 
\end{theorem}

More explicitly, the implied $O$ constants and $o$ sequence in bounds for Theorem 2 are understood to depend only on the implied $O$ constants and $o$ sequences in Condition 1.

The theorem shows that Simple Forward Selection can obtain asymptotically the same convergence rates (specifically $\sqrt{s_0 \log p /n}$ for the quantities $\En [(x_i'\theta_0 -\x_{ i}'\hat \theta)^2]^{1/2} $) as other high-dimensional estimators like Lasso, provided an appropriate threshold $t$ is used.  In addition, it selects a set with cardinality commensurate with $s_0$.

Finally, two direct consequences of Theorem 2 are   bounds on the deviations $\| \hat \theta - \theta_0\|_1$ and $  \| \hat \theta - \theta_0\|_2$ of $\hat \theta$ from the underlying unknown parameter $\theta_0$.  Theorem 3 above shows that deviations of $\hat \theta$ from $\theta_0$ also achieve rates typically encountered in high-dimensional estimators like Lasso.

\begin{theorem} Consider a sequence of data sets $\mathscr D_n$ with parameters $\theta_0$ and thresholds $t$ which satisfy Condition 1.  Suppose $\hat \theta$ is obtained by Algorithm 1.  Then there are bounds  $$\| \theta_0 - \hat \theta \|_2 = O\(\sqrt{ \frac{s_0 \log p }{n}} \hspace{.5mm} \) \text{ and }\| \theta_0 - \hat \theta \|_1 = O\(\sqrt{\frac{ s_0^2 \log p }{n  }} \hspace{.5mm}  \)$$ which hold with probability $1-o(1)$ as $n \rightarrow \infty.$
\end{theorem}  

%\textcolor{black}{As a final remark, inspection of the proofs Theorems 1 - 3 reveals that the arguments can be adjusted slightly to give bounds on $\hat s$ without the need for the condition that $t^{1/2} \geq 2 \varphi_{\min}(K_n)(G)^{-1} \| \En[\eps_i x_i] \|_\infty$.  An alternative involves conditions which restrict the maximum of quantities of the form $\| \En[ z_i \eps_i] \|_\infty$ where $z_i$ arise from different orthogonalizations of the original covariates $x_i$.  The specific $z_i$ correspond to components of QR decompositions over different orderings of variables.  Additional details for this possible extension are discussed in the proof of Theorem 1. }

\section{Testing-Based Forward Model Selection}

The previous section presented results on convergence rates of Simple Forward Selection.  The results of Theorem 1 are useful in developing intuition and proof techniques for inductive variable selection algorithms.  However, in terms of practical implementation, Section 2 leaves the question of how to choose a threshold unanswered.   This section develops TBFMS in order to analyze feasible, data-driven ways to decide which covariates to select, and when to stop selecting.

\subsection{Framework}

The basic framework for this section is similar to the earlier one.  Again, the observed data is given by $\mathscr D_n = \{(x_i,y_i)\}_{i=1}^n$, is generated by $\Pr$, and $y_{i} = x_i ' \theta_0 + \varepsilon_i$ for a parameter $\theta_0$ which is sparse with $s_0$ non-zero components supported on $S_0$.  Define $\ell(\theta)$ and $\ell (S)$ as before.

Define the expected loss function $\mathscr E:\mathbb R^p \rightarrow \mathbb R$ by $$\mathscr E(\theta) = \Ep \left [ \En [(y_i - x_i'\theta)^2] \right ]$$ where $\Ep$ is expectation with respect to $\Pr$.  Note that $\mathscr E(\theta ) = \Ep \ell(\theta)$. 
Extend the definition of $\mathscr E$ to apply also as a map $\E: 2^{\{1,...,p\}} \rightarrow \mathbb R$ by  $\E(S) = \min_{\theta: \text{supp}(\theta) \subseteq S} \E(\theta).$  Similarly to before, 
for any $S$ define the incremental loss from the $j$th covariate by
$$\Delta_j \E(S) = \E(S \cup \{j\} ) -   \E(S).$$ 

Within the class of greedy algorithms, it would be preferable to consider a greedy algorithm which inductively selects the $j$th covariate to enter a working model if $\Delta_j\E(S)$ is large \textcolor{black}{in absolute value} and $\textcolor{black}{-}\Delta_j \E (S) \geq \textcolor{black}{-}\Delta_k \E (S)$ for each $k \neq j$.  However, because $\Delta_j \E(S)$ cannot generally be directly observed from the data, the idea that follows is to make use of statistical tests to gauge the magnitude of $\Delta_j \E(S)$.  
Consider a set of tests given by
$$T_{jS\alpha} \in \{0,1\} \text{ associated to } H_0: \Delta_{j} \E (S) =0 \text{ and level } \alpha>0.$$

\noindent Assume that the tests reject ($T_{jS\alpha} = 1$) for large values of a test statistic $W_{jS}$.  

The model selection procedure is as follows.  Start with an empty model (consisting of no covariates).  At each step, if the current model is $\hat S$, select one covariate such that $T_{j\hat S\alpha} = 1$, append it to $\hat S$, and continue to the next step;  if no covariates have $T_{j\hat S\alpha}=1$, then terminate the model selection procedure and return the current model.  If at any juncture, there are two indices $j,k$ (or more) such that $T_{j\hat S\alpha} = T_{k\hat S\alpha}=1$, the selection is made according to the larger value of $W_{j\hat S}, W_{k\hat S}$.   

The use of $T_{jS\alpha}$ to define the set of covariates \textit{eligible} for entry into the model, and $W_{jS}$ to select \textit{which} eligible covariate actually enters is conceptually important: it dissociates and highlights the two fundamental tasks of regularization and fitting.  

%Alternatively, additional tests $T_{jkS\alpha}$ associated to $H_0:\Delta_{j} \E(S) \geq \Delta_k \E(S)$ could be devised to break ties.  The test statistic approach is natural for breaking potential multi-way ties. 
To summarize, the algorithm for forward selection given the hypothesis tests $(T_{jS\alpha},W_{jS})$ is now given formally.

{

\

\noindent {{\textbf{Algorithm 2.} \textit{ Testing-Based Forward Model Selection}}}

 \textit{Initialize}  \noindent Set $\hat S = \varnothing$

 \textit{For $1 \leq k \leq p$} 

\quad \textit{If} $T_{j\hat S\alpha} = 1$ for some $j\notin \hat S$ 

\quad \quad \textit{Set} $\hat j \in \text{arg} \max   \{W_{j\hat S}: T_{j\hat S \alpha} =1  \}$

\quad \quad \textit{Update }  $\hat S = \hat S \cup \{\hat j\}$ 

\quad \textit{Else}  \textit{Break}

\textit{Set } $\hat \theta \in \arg \min_{\theta:\supp(\theta) \subset \hat S}  \En[ (y_i - x_i'\theta)^2]$.

}

\subsection{Formal Analysis}

This section formally states conditions on the hypothesis tests and conditions on the data before analyzing properties of {Algorithm 1}.  These conditions are measures of the quality of the given testing procedure and the regularity of the data.

\

\noindent \textbf{Condition 2} (\textit{Hypothesis Tests}).
There is an integer $K_{\text{test}} > s_0$ and constants $ \alpha$, $\delta_{\text{test}}$, $c_{\text{test}}$, $c_{\text{test}}'$, $c_{\text{test}}''>0$ such that each of the following conditions hold.

\begin{itemize}
\item[1.]The tests have power in the sense that   
$$\Pr \( \{ T_{jS\alpha} =1 \ \ \text{for every} \  j, |S| \leq K_{\text{test}} \ \text{such that } -\Delta_{j} \mathscr E(S) \geq  c_{\text{test}} \} \)\geq 1 - \frac{1}{3}\delta_{\text{test}}.$$

\item[2.] The tests control size in the sense that
$$\Pr \( \{ T_{jS\alpha} = 1 \text{ for some } j,|S| \leq K_{\text{test}} \ \text{such that }  -\Delta_{j} \E(S) \leq  c_{\text{test}}' \} \) \leq \alpha+ \frac{1}{3}\delta_{\text{test}}.$$

\item[3.] The tests are continuous in the sense that   
$$\Pr (\{ -\Delta_j \E(S) \geq - c_{\text{test}}'' \Delta_k \E(S)  \text{ for each } \ j,k,|S| \leq K_{\text{test}} \text{ such that } $$ $$T_{jS\alpha} =1 \ \text{and} \  W_{jS} \geq W_{kS}  \} ) \geq 1-  \frac{1}{3}\delta_{\text{test}}.$$
%$$\Pr ( \{W_{jS} \geq W_{kS} \text{ for each } \ j,k,|S| \leq K_{\text{test}} \text{ such that } $$ $$T_{jS\alpha} =1 \ \text{and} \ -\Delta_j \E(S) \geq - c_{\text{test}}'' \Delta_k \E(S)  \} ) \leq  \frac{1}{3}\delta_{\text{test}}.$$
%$$\Pr ( \{-\Delta_j \E(S) \geq - c_{\text{test}}'' \Delta_k \E(S)  \text{ for each } j,k,|S| \leq K_{\text{test}} \text{ such that } $$ $$ T_{jS\alpha}=1 \text{ and } W_{jS} \geq W_{kS}   \} ) \geq 1 -  \frac{1}{3}\delta_{\text{test}}.$$

\end{itemize}

\

 The constants $c_{\text{test}}$ and $c_{\text{test}}'$ measure quantities related to the size and power of the tests and provide a convenient language for subsequent discussion.  The constant $c_{\text{test}}''$ measures the extent to which the test statistics $W_{jS}$ reflect the actual magnitude of $\Delta_{j}\mathscr E(S)$.  Note again that the hypothesis tests are considered simply tools for model selection, which coincidentally have many properties in common with traditional inferential hypothesis tests.

\

\noindent \textbf{Condition 3} (\textit{Regularity}).
Normalizations $\Ep[\En[x_{ij}^2]] = 1$ hold for all $j$. The residuals decompose into $\varepsilon_i = \varepsilon_i^{\text{o}} + \varepsilon_i^{\text{a}}$ where $\Ep [ \En [ {\varepsilon_i^{\text{o}}}^2]] < \infty$, $\Ep [ \En [ \varepsilon_i^{\text{o}} x_{ij} ]] = 0$ for all $j$, and  $\Ep[\En[{\eps^{\text{a}}_i}^2]] \leq \frac{1}{2}\varphi_{\min}(K_{\text{test}})(\Ep[G])^{-1} c_{\text{test}}'$. Finally, $  ( \textcolor{black}{80 \times}\varphi_{\min}(K_{\text{test}})(\Ep[G])^{-\textcolor{black}{4}}{c_{\text{test}}''}^{-3}+1)s_0 < K_{\text{test}}$.

\

%which are called the orthogonal components and the sparse approximation components. 

Condition 3 imposes regularity conditions for the class of models considered in the following theorem.  First, $\eps_i$ is decomposed into an orthogonal component $\eps_i^{\text{o}}$ and an approximation component $\eps_i^{\text{a}}$, each of which exhibits a different kind of regularity.  The orthogonal component is orthogonal to the covariates in the population.  The approximation component need not be orthogonal to the covariates, but its magnitude must be suitably controlled by the sparse eigenvalues of $\Ep[G]$ and by the parameter $c_{\text{test}}'$, which is a detection threshold for the profile of hypothesis tests $T_{jS\alpha}$.   This decomposition allows for approximately sparse models similar to the framework of \cite{BellChenChernHans:nonGauss}.    The fact that $\eps_i^{\text{a}}$ need not be orthogonal to the covariates also allows this framework to overlay onto many problems in traditional nonparametric econometrics.

Condition 3 also imposes conditions on the relative values of the sparse eigenvalues of $\Ep[G]$, $c_{\text{test}}'', s_0,$ and $K_{\text{test}}$.  Note that $K_{\text{test}}$ measures the size of the set $S \subset \{1,...,p\}$ over which the hypothesis tests perform well, as defined by Condition 2.    Consequently, this condition requires that the hypothesis tests $T_{jS\alpha}$ perform sufficiently well over sets $S$, which must be larger when $\Ep[G]$ has small eigenvalues, when $c_{\text{test}}''$ is small, or when $s_0$ is large. 

There are a few cases where Condition 3 can be simplified.  If $p >n$, even though the empirical Gram matrix is necessarily rank deficient, the population Gram matrix may be full rank.    When $\Ep[G]$ is full rank, then $\lambda_{\min}(\Ep[G])$ may be used in place of $\varphi_{\min}(K_{\text{test}})(\Ep[G])$.    In addition, the condition on $\eps_{i}^{\text{a}}$ implicitly imposes constraints on $c_{\text{test}}'$ and  $\varphi_{\min}(K_{\text{test}})(\Ep[G])^{-1}$.  When there is no approximation error, this requirement is no longer needed.  %In the application below, $c_{\text{test}}''$ will satisfy $c_{\text{test}}'' = O(1)$.

Let \begin{align*}
c_{\text{T}}& = s_0 \varphi_{\min}(K_{\text{test}})(\Ep[G])^{-1}c_{\text{{test}}}  \\
c_{\text{T}\hspace{.5mm}}'& =  \textcolor{black}{80 \times }\varphi_{\min}(K_{\text{test}})(\Ep[G])^{-\textcolor{black}{4}}{c_{\text{{test}}}''}^{-3} \textcolor{white}{\Big |} \\
c_{\text{T}}''(\hat s)& =  \varphi_{\max}(\hat s + s_0)(G)^{1/2}   \varphi_{\min}(\hat s + s_0)(G)^{-1/2}\hat s^{1/2}  \hspace{.5mm} \| \En[x_i\eps_i] \|_\infty \\ & \ \ \ +  3\varphi_{\max}(\hat s + s_0   )(G)  \varphi_{\min}(\hat s + s_0)(G)^{-1/2}  {({\hat s + s_0})^{1/2} {c_{\text{test}}} ^{1/2}} \varphi_{\min}(K_{\text{test}})(\Ep[G])^{-1} .\end{align*}

\

\begin{theorem}

Consider  $\mathscr D_n \sim \Pr$ for some fixed $n$ and $\{T_{jS\alpha},W_{jS}\}$ such that Conditions 2 and 3 hold.  Suppose $\hat \theta$ is obtained by Algorithm 2. Then the bounds 
$$\E(\hat S) - \E(S_0) \leq c_{\text{\em T}}$$  
$$ \hat s \leq (c_{\text{\em T}}' +1)\hspace{.5mm} s_0$$
$$\En [ (\hspace{.5mm} x_i'\theta_0 - \x_{ i}'\hat \theta \hspace{1mm} )^2 ]^{1/2} \leq c_{\text{\em T}}''(\hat s\hspace{.5mm}) $$

\noindent hold with probability at least $1 - \alpha - \delta_{\text{\em test}} $.

\end{theorem}

Theorem 4 provides finite sample bounds on the performance of TBFMS.  % Theorem 4 works with the fact that $-\Delta_{j} \mathscr E(S) > c_{\text{test}}$ on the event implied by Condition 2, instead of the simpler fact $-\Delta_j \ell(S) > t$ as was used in Theorem 1. 
 In contrast to Theorem 1, Theorem 4 also addresses the possibility that if covariate $j$ is selected ahead of covariate $k$, it is not necessarily the case that $-\Delta_{j} \mathscr E(S) > -\Delta_{k} \mathscr E(S)$.  This is done by making use of the continuity constant $c_{\text{test}}''$ in Condition 2.

Theorem 4 can be used to derive asymptotic estimation rates by allowing the constants to change with $n$.  The next subsection provides an example to a linear model with heteroskedastic disturbances, where, under the stated regularity conditions, the prediction and estimation error attain the rate $O (s_0 \log p/n)$.  This convergence rate matches typical Lasso and Post-Lasso rates.  

The results aim to control the hybrid objectives, described in the introduction, of producing a good fit and returning a sparse set of variables.  One useful implication of bounds controlling both $\hat s$ and $\En[ (x_i'\theta_0 - \x_{ i}\hat \theta)^2]$ is that the results can be applied to constructing uniformly valid post-model selection inference procedures (see \cite{BCH-PLM}), in which for some applications, the prediction error bound alone is insufficient.  %Applications to inference are discussed in Appendix C.

%Though the statement and conditions of Theorem 4 are very simple, it is convenient to have an version which does not require a bound on the eigenvalues of the population Gram matrix $\Ep[G]$.  This is useful, for instance, when the covariates $\{x_i \}_{i=1}^n$ are considered fixed or when conditioning on covariates.

\section{Examples and Extensions}

This section describes an example application of Theorem 4.  The main theoretical application is an illustration with heteroskedastic data.  For this setting, a TBFMS procedure is constructed for which optimal convergence rates are proven.  %Next, several variants are also discussed, as well as investigated in a simulation study that follows.   

%An important aspect of the result in the previous section is that it explicitly allows the researcher to easily dissociate the \textit{regularization} component of high-dimensional estimation with the \textit{fitting} component.  For instance, in the heteroskedasticity example, test statistics can be redefined $\tilde W_{jS} = \Delta_{j}\ell(S)$ while $T_{jS\alpha}$ remain those defined by Definition 1.  This is a slightly ``greedier'' formulation.  Under condition 5, this modification can be shown to achieve the same rates of convergence (by showing that the ratio $ W_{jS} / \tilde W_{jS}$ is sufficiently well behaved.).

\subsection{Heteroskedastic Disturbances}

This section gives an example of the use of Theorem 4 by illustrating an application of model selection in the presence of heteroskedasticity in the disturbance terms $\eps_i$.  A TBFMS procedure is constructed based on the Heteroskedasticity-Consistent standard errors described in \cite{white:het}.  The conditions required for the application of Theorem 4 are verified under low-level conditions on data generating processes.    Other TBFMS procedures are possible, and these are discussed in the next section.  The analysis begins with a formulation stated in \cite{DK:TBFMS:PandP} (which does not derive nor claim any theoretical properties.)

For shorthand, write $x_{ijS}$ \textcolor{black}{(with $j \notin S$)} to be the vector with components $x_{ik}$ with $k = j$ or $k \in S$.  To construct the tests, begin with the least squares estimate of the regression $y_i$ on $x_{ijS}$.    
$$\hat \theta_{jS} = \En[ x_{ijS} x_{ijS}' ]^{-1}  \En[  x_{ijS}'y_i ]$$ 

\noindent Define  $\hat \varepsilon_{ijS} = y_i - x_{ijS}'\hat \theta_{jS}.$
One heteroskedasticity robust estimate of the sampling variance of $\hat \theta_{jS}$, proposed in \cite{white:het}, is given by the expression

$$\hat V_{jS}=\frac{1}{n} \En  [x_{ijS} x_{ijS}' ]^{-1} \Psi_{jS}^{\hat \varepsilon}  \En  [x_{ijS} x_{ijS}']^{-1}$$

\noindent where $$\Psi_{jS}^{\hat \varepsilon} = \En [ \hat \varepsilon_{ijS}^{ \ 2} x_{ijS} x_{ijS}' ] .$$

\noindent Define the test statistics

$$W_{jS}^{\mathsf{het}} = [\hat V_{jS}]_{jj}^{-1/2} \left |[\hat \theta_{jS}]_j \right |.$$

\noindent Reject $H_0$ for large values of $W_{jS}^{\mathsf{het}}$ defined relative to an appropriately chosen threshold.  
To define the threshold, first let $\eta_{jS}:= (1 \ , - \beta_{jS}')'$  where $\beta_{jS}$ is the coefficient vector from the least squares regression of $\{ x_{ij} \}_{i=1}^n$ on $\{ x_{ik} \}_{i=1, k\in S}^n$.   Then define

$$\hat \tau_{jS} = \frac{ \| \eta_{jS} ' \text{Diag}(\Psi_{jS}^{\hat \varepsilon})^{1/2} \|_1}{  \sqrt{  \eta_{jS}'\Psi_{jS}^{\hat \varepsilon} \eta_{jS} }        }.$$

\noindent The  $\hat \tau_{jS}$ will be helpful in addressing the fact that many different model selection paths are possible under different realizations of the data under $\Pr$.\footnote{There is an unfortunate misprint in a Papers and Proceedings version of this paper, \cite{DK:TBFMS:PandP}, in which the exponent $1/2$ is missing from the term $\text{Diag}(\Psi_{jS}^{\hat \varepsilon})$.}    Not taking this fact into account can potentially lead to false discoveries.   The next condition states precisely the hypothesis tests $T_{jS\alpha}$.

\begin{definition}[\textit{Hypothesis Tests for Heteroskedastic Disturbances}]  Let $c_\tau>1$ and $\alpha > 0$ be parameters.  Assign $W_{jS}=W_{jS}^{\mathsf{het}}$.   Assign
 $$T_{jS\alpha} = 1 \iff W_{jS} \geq c_{\tau}\hat \tau_{jS}\Phi^{-1}(1-\alpha/p).$$

 \end{definition}

The term $\Phi^{-1}(1-\alpha/p)$ can be informally thought of as a Bonferroni correction term that takes into account the fact that there are $p$ potential covariates.   The term $c_\tau\hat \tau_{jS}$ can be informally thought of as a correction term that can account for the fact that the set $S$ is random and can have many potential realizations.  The simulation study uses the settings $c_\tau = 1.01$ and $\alpha = .05$.

 \

\noindent \textbf{Condition 4} (\textit{Regularity for Data with Heteroskedasticity}).  Consider a sequence of data sets $\mathscr D_n = \{(x_i,y_i)\}_{i=1}^n \sim \Pr = \Pr_n$.  The observations $(x_i,y_i)$ are i.n.i.d. across $i$ and
$y_i = x_i ' \theta_0 + \varepsilon_i$ for some $\theta_0$ with $s_0 = o(n)$.  The residuals decompose into $\eps_i = \eps^{\text{o}}_i + \eps^{\text{a}}_i$ such that a.s.,  $\Ep[ \varepsilon_i^{\text{o}} | x_i] = 0$ and $\max_{i} | \eps_{i}^{\text{a}}| = O(n^{-1/2})$.  In addition, a.s., uniformly in $i$ and $n$, $\Ep [\varepsilon_i^4|x_i] $ are bounded above and $\Ep[\varepsilon_i^2|x_i]$ is bounded away from zero.   The covariates satisfy  $\max_{j \leq p} \En[x_{ij}^{12}] = O(1)$ with probability $1 - o(1)$.  There is a sequence $K_n$, where $s_0 = o(K_n)$, and bounds
$\varphi_{\min}(K_n)(G)^{-1} = O(1)$, $\varphi_{\max}(K_n)(G) = O(1)$, $ \varphi_{\min}(K_n)( \En [(\varepsilon_{i} x_{i})(\varepsilon_{i}x_{i})' ] )^{-1}  = O(1),$ and $\max_{|S| < K_n, j \notin S} \| \eta_{jS} \|_1= O(1)$,
which hold with probability $1-o(1)$.  The rate condition $K_n^4 \log^3 p/n = o(1)$ holds.
  
\

Condition 4, as before, gives conditions on the sparse eigenvalues, this time applying to both $G$ and to $\En[(\eps_i x_i)(\eps_i x_i')]$.  %The requirement that $\varphi_{\min}(K_n)(G)^{-1}$ stay bounded can be relaxed, though the cost is that the proof of Theorem 5 could not directly call upon Theorem 4 (since $c_{\text{T}}''$ involves a maximal sparse eigenvalue).  
In addition, Condition 4 assumes a bound on $\eta_{jS}$ that may be strong in some cases.  Previous results in \cite{Tropp:Greed}, \cite{Zhang:GreedyLeastSquares} assume the strict condition that $\max_{j \notin S_0} \| \eta_{jS_0} \|_1 < 1$, which is the genuine irrepresentability condition, in analysis of inductive variable selection algorithms.  Here, the requirement $< 1$ is replaced by the weaker requirement $=O(1)$.  Other authors, for instance \cite{MY2006}, use conditions analogous to $\max_{|S| \leq K_n, j \notin S} \| \eta_{jS} \|_1= O(1)$ in the context of learning high-dimensional graphs, and note that the relaxed requirement is satisfied by a much broader class of data-generating processes.  Analogous bounds on $\| \eta_{jS}\|_1$ were not required in Theorem 1, because its proof does not leverage bounds relating $\hat W_{jS}$ to the self-normalized sums $\En[x_{ij}' \eps_i] / \sqrt{ \En[ x_{ij}^2 \eps_i^2]}$, $j \leq p$.  Failure of the $O(1)$ bound would lead only to slightly slower convergence rates.
Condition 4 also states regularity conditions on moments of $\varepsilon_i$ and $x_i$, which are useful for proving laws of large numbers, central limit theorems, and moderate deviation bounds (see \cite{jing:etal}).   Finally, the rate condition controls relative sizes of $s_0,p,n$ because $s_0 < K_n$.

\begin{theorem}  Consider a sequence of data sets $\mathscr D_n \sim \Pr = \Pr_n$ which satisfies Condition 4.  Suppose that $c_{\tau}>1$ is fixed independent of $n$, and  that $\alpha = o(1)$ with $n\alpha \rightarrow \infty$.\footnote{Allowing $\alpha$ to be fixed is possible under more restrictive conditions on the approximation error terms $\eps_i^{\text{a}}$.  If $p > n$, then the rate $\log(p/ \alpha)$ becomes equivalent to simply $\log (p)$.} Let $\hat \theta$ be the estimate obtained from Algorithm 2 with tests defined by Definition 1.  Then there are bounds
$$\En[ (x_i'\theta_0 -x_i'\hat \theta)^2 ]^{1/2} = O\(\sqrt{\frac{s_0 \log (p/\alpha)}{n}} \ \) \ \ \ \text{and} \ \ \ \hat s =O (s_0),$$
which hold with probability at least $1- \alpha - o(1)$  as $n \rightarrow \infty$. 
\end{theorem}

%The theorem is proven by appealing to Theorem 4.  

%The condition $\| \eta_{jS} \|_1 = O(1)$ is potentially restrictive (see discussion above).  If instead the unrestrictive condition $\| \eta_{jS} \|_1= O( |S|^{1/2} )$ holds, then the following similar result can be shown. $\En[(x_i'\theta_0 -x_i'\hat \theta)^2]  = O_{\Pr}(s_0^2 \log p/n) $ and $\hat s = O(1)s_0$ with probability $1-o(1)$.  Note that the rate of convergence is slower by a factor of $s_0$.

\subsection{Additional TBFMS Formulations and Variants}

In general, the quality of statistical performance of a variant of TBFMS may depend on the structure of the data at hand through the size, power, and continuity properties of the tests, as articulated in Condition 2.  Theorem 4 is general and can thus potentially be applied for many different types of TBFMS procedures, depending on how $W_{jS}$ and $T_{jS\alpha}$ are defined in a particular setting.  Depending on the setting, some variants may have better size, power, and continuity properties than others. 
This section describes several variants of the TBFMS procedure defined in the previous section.   

The definition of the first variant considered is based on the observation that there is a much simpler formulation for the hypothesis tests that ignores the $c_{\tau} \hat \tau_{jS}$ terms.  This results in the following definition.

\begin{definition}[\textit{Simplified Hypothesis Tests for Heteroskedastic Disturbances}]   Let $\alpha > 0$ be parameters.  Assign
 $$T_{jS\alpha} = 1 \iff W_{js} \geq\Phi^{-1}(1-\alpha/p).$$
 \end{definition}

These tests are based on a simple Bonferroni-type correction.  Furthermore, though never previously formally justified, TBFMS using the simpler tests is natively programmed in some statistical software, including SPSS and Stata. It is unknown to the author at the time of this writing whether the same convergence rates can be attained using the simpler tests under the identical regularity as in Condition 4.  This option is explored in some finite sample settings in the simulation study that follows.  Evidence from the simulation study suggests that this option performs better than the more complex tests defined in Definition 1.  The tests in Definition 2 are not necessarily more conservative than those in Definition 1.

The next variant TBFMS procedure for heteroskedastic data illustrates an important aspect of the result of Theorem 4.
Namely, Theorem 4 explicitly allows the researcher to easily dissociate the \textit{regularization} component of high-dimensional estimation with the \textit{fitting} component.  For instance, the following formulation may be described as slightly greedier than Definition 2.

\begin{definition}[\textit{Fit-Streamlined Hypothesis Tests for Heteroskedastic Disturbances}]  Let $\alpha > 0$ be a parameter.  Assign $W_{jS}=\Delta_j \ell(S)$.   Assign
 $$T_{jS\alpha} = 1 \iff W_{jS}^{\mathsf{het}} \geq \Phi^{-1}(1-\alpha/p).$$
  
 \end{definition}

Under the more conservative tests, $T_{jS\alpha} = 1 \iff W_{jS}^{\mathsf{het}} \geq c_{\tau}\hat \tau_{jS} \Phi^{-1}(1-\alpha/p)$ for some $c_\tau>1$, the same convergence rates of this greedier TBFMS procedure under Condition 4 are proven in the same way as Theorem 5  (by showing that the ratios $ \Delta_j \ell(S) / W_{jS}^{\textsf{het}}$ are sufficiently well behaved).  For brevity, the details are omitted.   
  %\begin{corollary}
%Consider a sequence of data sets $\mathscr D_n \sim \Pr = \Pr_n$ which satisfies Condition 4.  Suppose that $c_{\tau}>1$ is fixed independent of $n$, and  that $\alpha = o(1)$ with $n\alpha \rightarrow \infty$. Let $\hat \theta$ be the estimate obtained from Algorithm 2 with tests defined by Definition 3.  Then there are bounds
%$$\En[ (x_i'\theta_0 -x_i'\hat \theta)^2 ]^{1/2} = O\(\sqrt{\frac{s_0 \log (p/\alpha)}{n}} \ \) \ \ \ \text{and} \ \ \ \hat s =O (s_0)$$
%which hold with probability at least $1- \alpha - o(1)$  as $n \rightarrow \infty$. 
 %\end{corollary}

When the data is approximately homoskedastic, the tests defined in Definition 1 may be too conservative and suffer in terms of power (noting that power is an explicit input into the bounds in Theorem 4).   In this case, tests of the following form can be considered using the homoskedastic-based test statistics $W_{jS}^{\mathsf{hom}} = [\hat V_{jS}^{\mathsf{hom}}]_{jj}^{-1/2} \left |[\hat \theta_{jS}]_j \right |$ with $\hat V_{jS}^{\mathsf{hom}} = \frac{1}{n} \En[\hat \eps_{ijS}^2] \En[x_{ijS}x_{ijS}']^{-1}$.

\begin{definition}[\textit{Simplified Hypothesis Tests for Homoskedastic Disturbances}]   Let $\alpha > 0$ be parameters.  Assign $W_{jS} = W_{jS}^{\textsf{hom}}$.  Assign
 $$T_{jS\alpha} = 1 \iff W_{js} \geq\Phi^{-1}(1-\alpha/p).$$
 \end{definition}

Evidence from the simulation study suggests that this option performs better than the more complex tests in Definition 1 when the data is homoskedastic, but not when heteroskedastic.  

\subsection{TBFMS with Baseline Covariates}

The initialization statement of Algorithm 2 can be modified trivially so that a baseline set $S_{\text{base}}$ of covariates are automatically included in $\hat S$ at the start of selection.   In this case, the initialization statement of Algorithm 2 is replaced with  

\

$\textit{Initialize.}  \text{ Set }\hat S = S_{\text{base}}.$ 

\

\noindent  Under this modification, a direct analogue of Theorem 4 holds.  It is proven using the same arguments.  {Sparsity bounds can be calculated as in the proof of Theorem 4 by separately tracking covariates $S_0 \setminus S_{\text{base}}$.  This requires an appropriate adjustment to Condition 3 in which $K_{\text{test}}$ must be larger, and in particular, bound a quantity depending on both $s_0$ and $| S_{\text{base}}|$.  The proof of the first and third statements of Theorem 4 do not depend on the initialization $S = \varnothing$. } Such modification is appropriate, for instance, in cases in which a researcher wishes to include a set of covariates into a model, but is unsure of which interactions to include.   In this case, TBFMS can be used to help identify relevant interaction terms.  This case is further explored in the empirical application.  Similarly, a constant term may be included automatically in the regression model.

\subsection{Additional Discussion of Potential Variants}

Analogous results potentially hold for dependent data using HAC-type standard errors (see \cite{newey:west}, \cite{andrews:hac1}.)  %The required central limit results for such an application are beyond the scope of this work, though using the moderate deviation results for dependent data of \cite{chen:biao:moderatedeviationsdependence} is a potential starting point.   
In addition, cluster-type standard errors for large-$T$-large-$n$ and fixed-$T$-large-$n$ panels can be used by adapting arguments from \cite{BCHK:Panel}.   Analogous results for homoskedastic disturbances can be derived as a corollary.  %The performance of homoskedastic variations on TBFMS are explored in the empirical example and in the simulations.

Another alternative is to consider generalized error rates.  The conditions set forth in Condition 2 require control of a notion resembling family-wise error rate uniformly over hypothesis tests $H_0: \Delta_j \mathscr E(S)=0$ for $j \leq p $ and $|S| < K_{\text{test}}$ for some integer $K_{\text{test}}$.  Other types of error rates like $k$-family-wise error rate, false discovery rate, or false discovery proportion are potentially possible as well.  In particular, the arguments in the proof of Theorem 4 would continue to be compatible with procedures that controlled an appropriate notion of false discovery proportion.  In order to keep exposition concise, these extensions are not considered here.

\section{Example: TBFMS for Asset-Based Poverty Mapping}

This section investigates the use of TBFMS to develop improved proxy-means tests in an application to poverty mapping in a Peruvian dataset covering years 2010--2011.   This analysis extends the original analysis in \cite{Hanna:Olken}, who estimated a predictive model of household consumption using the same data.  The data is from the Peruvian Encuesta Nacional de Hogares (ENAHO), maintained by the Instituto Nacional de Estadistica e Informatica (INEI), Peru.

Fighting poverty is a major priority for many developing countries and international organizations like the United Nations.  Strategies for combating poverty often require that governing bodies have accurate information about household level consumption, income, or other measures of welfare.  Methods for empirical identification of households and regions below a given poverty line is are considered in e.g. \cite{elbers2003micro}, among others.  %Once quality information about extent and location of poverty is available, instruments for fighting poverty, including conditional and unconditional cash transfers, can be implemented more efficiently (see e.g. \cite{RePEc:wbk:wbpubs:2597}).  

One method in which a government can obtain a signal about measures of welfare is termed a \textit{proxy-means test.} The implementation of a proxy-means test is usually based on large censuses of the population, in which government enumerators obtain information on easily observable and verifiable assets.  %For instance, typical observed assets may include ownership of items such as televisions and refrigerators, the type of material used in a household's roof, floor, and walls, etc.  
The government uses these assets to predict incomes or per-capita consumption or other measures of poverty or welfare by estimating a regression on a smaller sample with detailed measurement of consumption. The \textit{proxy-means score} is defined as the predicted income or consumption, which is calculated using the results from the predictive regression.
This method is widespread, and is implemented in several countries including Indonesia, Pakistan, Nigeria, Mexico, Philippines, Burkina Faso, Ecuador, Jamaica.  Improved predictions may be helpful to policy makers in deciding on strategies to eliminate poverty; see \cite{RePEc:wbk:wbpubs:2597}.    %Any means by which the predictive quality of the proxy-means score can be improved may lead to more efficiently allocated benefits.

%There are a variety of methods to predict income (or consumption), but they all share basic features: The government takes another dataset that was collected in a low-stakes context (for example, one collected not for targeting, but rather just for research purposes), and therefore for which households have no strong incentive to lie. In this dataset, the government observes the same asset variables as in the proxy- means census and also observes a measure of poverty, such as a household?s monthly income or per-capita expenditure. 

%The resulting proxy-means scores can be used to determine eligibility for benefits for individual households.  As a result, eligibility depends on predicted income instead of actual income.  

%\subsection{Estimation Methodology}

The data contain covariates which are indicator variables describing household-level assets and which are used to predict outcomes $y$ = Consumption (in $10^3$ Peruvian Soles) as well as $y=$ \textit{log} Consumption (in Peruvian Soles).   The set of 46305 household observations is split randomly (with equal probability) into a training sample of size 22674 and a testing sample of size 22704.  All estimation procedures considered are implemented on the training sample.  The indicators derive from factor variables describing a household's (1) water source, (2) drain infrastructure (3) wall material, (4) roof material, (5) floor material, (6) availability of electricity, (7) access to telephone, (8) education of head of household, (9) type of insurance, (10) crowding, (11) consumption of luxury items.  See \cite{Hanna:Olken} for more details.% Properties of the estimates are studied based on predicted consumption and predicted log consumption on the testing sample.  
%For additional details, see \cite{Hanna:Olken}.  

Here, TBFMS is used to determine which interactions of underlying covariates described above are helpful in developing an improved proxy-means test.\footnote{\cite{Hanna:Olken} also note that other estimation techniques may offer improvement.  The reference \cite{ml:poverty:mcbride:nichols} also investigates various cross validation and machine learning techniques toward this end.}  
%In this analysis, TBFMS is used to select possible interaction terms of the original set of 72 indicators, and investigate if selected terms may be used to improve significantly over the specification considered in \cite{Hanna:Olken}.\footnote{Note that there are several reasons to include original variables at the beginning of model selection.  For one, the original variables may often have been chosen by experts, familiar with what should be relevant from previous experience both empirical and theoretical.  In general, economic theory may not always provide enough information about functional forms to definitively answer the question of exactly which interaction terms of important variables should be included ex ante. }
For every (unordered) pair of indicator variables $A$ and $B$, three symmetric logic functions, 
$$\mathsf{and}(A,B), \ \mathsf{or}(A,B),\ \mathsf{xor}(A,B).$$
\noindent are generated.  Together with a constant term, these logic functions linearly span exactly the set of symmetric boolean functions on all pairs of indicators.  
In addition, for every unordered triple of indicator variables, $A,B$ and $C$, the function $$\mathsf{and}(A,B,C)$$
\noindent is generated. The final interaction expansion is based on the union of the generated functions. 
 %Though this number is not larger than the training sample size, the number of terms is large enough to create over-fitting concerns.  At the same time, the expansion is simple and intuitive.  Higher order expansions could also be considered.   %\footnote{There are 16 boolean functions on two binary variables.  These include the above five functions, their negations, two function which depend only on $A$, two function which depend only on $B$, and two functions which depend on neither $A$ nor $B$.}
%This paper uses the same data as \cite{Hanna:Olken} to investigate the performance of TBFMS as a predictor of consumption.  

\begin{table}[h] \caption
 {Asset-Based Poverty Estimation } 
\begin{tabular*}{\textwidth}{p{3.3cm} p{.7cm}  p{.9cm} p{.9cm} p{1cm} p{1cm} p{1.2cm} p{1.3cm} p{1.3cm} } 
\hline          \hline                                                                  \\ 
&	\textcolor{white}{\Big |}  \# Var	& Train MSE & Test MSE &  Rel(\%) $ \Delta$MSE &	 $\Delta $MSE	& $p$-val  $\Delta \text{MSE}$   &  \multicolumn{2}{c}{95\% CI $\Delta$MSE} \\   
 \\ \cline{2-9}     & \multicolumn{8}{c}{A. $y= \text{Consumption} \ (\times 10^3 \ \text{Peruvian Soles})$ \textcolor{white}{\Big |}  }\\ \cline{2-9} 
 \underline{Estimates}\\ 
\textit{OLS	}&	62	&	0.0780	&	0.0897	&	--	&		--	&	--	& \ \ \ 	--	& \ \ \ 	--	\\  
%\textit{TBFMS I  + Base} & 63  		&	0.0779	&	0.0897	&	0.0002		&	0.0000	&	 0.0455 &   ( 0.0000	&  0.0000 ) 		\\
%\textit{TBFMS II + Base} & 142 			&	0.0706	&  0.0862	& 3.8939	&	0.0035			&	0.0001 & ( 0.0017	&	0.0053 )		\\
%\textit{TBFMS III  + Base}&	103	&	0.0684	&	0.0842	&	6.1090		&	0.0055 	&4.5{\tiny\textrm E}-06	&	( 0.0031	&	0.0078 )	\\
%\textit{TBFMS IV + Base}&	103	&	0.0684	&	0.0842	&	6.1090		&	0.0055 	&4.5{\tiny\textrm E}-06	&	( 0.0031	&	0.0078 )	\\ \\
%\textit{TBFMS I	}&	79	&	0.0792	&	0.0932	&	-3.9264		&	-0.0035 	&1.6{\tiny\textrm E}-08	&	( -0.0047	&	-0.0023)	\\
%\textit{TBFMS II	 }&	169	&	0.0709	&	0.0872	&	2.7911		&	0.0025 &0.0014&		( 0.0010	&	0.0040 )	\\ 
%\textit{TBFMS III }	&	61	&	0.0687	&	0.0842	&	6.2080		&	-0.0056	&	3.9{\tiny\textrm E}-06	&	( 0.0032	&	0.0079 )	\\ 
%\textit{TBFMS IV }&	61	&	0.0687	&	0.0842	&	6.2080		&	-0.0056	&	3.9{\tiny\textrm E}-06	&	( 0.0032	&	0.0079 )	\\  \\
\textit{TBFMS I } &   	116	&	0.0727	&	0.0867	&	3.3503		&	0.0030	&	 2.3{\tiny \textrm{E}}-08  & ( 0.0020	& 0.0041 ) 		\\
\textit{TBFMS II  } & 159 			&	0.0697	&  0.0860	& 4.1264	&	0.0037			&	0.0002 & ( 0.0018	&	0.0056 )		\\
\textit{TBFMS III   }&	118	&	0.0654	&	0.0848	&	5.4865		&	0.0049 	&0.0002	&	( 0.0023	&	0.0076 )	\\
\textit{TBFMS IV   }&	118	&	0.0654	&	0.0848	&	5.4865		&	0.0049 	&0.0002	&	( 0.0023	&	0.0076 )	\\
%\textit{OLS	+ Interaction}&	7730	&	0.06150	&	0.0931	&	-3.7334 	& -0.0033 &	0.0588 & 	( -.0068	&	0.0001 )	\\ 
 \\ \cline{2-9}     & \multicolumn{8}{c}{B. $y= log \ \text{Consumption} \ (\text{Peruvian Soles} ) $ \textcolor{white}{\Big |}  }\\ \cline{2-9} 
 \underline{Estimates}\\
\textit{OLS}	&	62	&	0.1916	&	0.1910	&	--	&		--	&	--	&	\ \ \ --	&\ \ \ 	--	\\ 
%\textit{TBFMS I  + Base }& 63  		&	0.1915	&	0.1911	&	-0.0487		&	-0.0000	&	 0.3257 & (  -0.0003	&	0.0001 )		\\
%\textit{TBFMS II + Base }& 141 			&	0.1841	&	0.1887	&	1.2248		&	0.0023	&	0.0006 & ( 0.0010	&	0.0037 )		\\
%\textit{TBFMS III  + Base }&	85	&	0.1835	&	0.1866	&	2.3385		&	0.0045 &1.8{\tiny \textrm E}{-11}&		( 0.0032	&	0.0058 )	\\ 
%\textit{TBFMS IV + Base} &	85	&	0.1835	&	0.1866	&	2.3385		&	0.0045 &1.8{\tiny \textrm E}{-11}&		( 0.0032	&	0.0058 )	\\ \\
%\textit{TBFMS I}	&	72	&	0.1985	&	0.2016	&	-5.5025		&	-0.0105	&	5.4{\tiny\textrm E}-33	&	( -0.0122	&	-0.0088)	\\
%\textit{TBFMS II }&	128	&	0.1844	&	0.1898	&	0.6317		&	0.0012	&	0.1309	&	( -0.0004	&	0.0028 )	\\ 
%\textit{TBFMS III}	&	53	&	0.1839	&	0.1876	&	1.8024		&	0.0034	&	4.0{\tiny\textrm E}-07	&	( 0.0021	&	0.0048 )	\\ 
%\textit{TBFMS IV} &	53	&	0.1839	&	0.1876	&	1.8024		&	0.0034	&	4.0{\tiny\textrm E}-07	&	( 0.0021	&	0.0048 )	\\  \\
\textit{TBFMS I   }& 117  		&	0.1849	&	0.1884	&	1.3954		&	0.0027	&	 2.1{\tiny \textrm{E}}{-07} & ( 0.0017	&	0.0037 )		\\
\textit{TBFMS II  }& 140 			&	0.1829	&	0.1892	&	0.9669		&	0.0018	&	0.0109 & ( 0.0004	&	0.0033 )		\\
\textit{TBFMS III   }&	89	&	0.1818	&	0.1865	&	2.3809		&	0.0045 &1.8{\tiny \textrm E}{-10}&		( 0.0031	&	0.0060 )	\\ 
\textit{TBFMS IV  } &	89	&	0.1818	&	0.1865	&	2.3809		&	0.0045 &1.8{\tiny \textrm E}{-10}&		( 0.0031	&	0.0060 )	\\ \\
%\textit{OLS + Interaction}	&	7730 &	0.1639	&	0.2023	&		-5.8929	&	-0.0113	&	6.3{\tiny\textrm E}-10&	( -0.0148	&	-0.0077)	\\ 

 %\multicolumn{9}{l}{*TBFMS procedure selected from 4 candidate procedures based on training sample MSE.} \\  \multicolumn{9}{l}{  Candidate procedures: TBFMS I -- IV with baseline variables at selection start.}\\
 %\multicolumn{9}{l}{  Selected procedure:  TBFMS III with baseline variables.}\\
 % \multicolumn{9}{l}{  Full results for all procedures and additional procedures without baseline variables in Appendix.}\\
\hline

Data description& \multicolumn{8}{l}{  Training sample size : 22674, Test sample size : 22704} \\
& \multicolumn{8}{l}{  \# Baseline characteristics : 62, \# Interacted characteristics : 23964 } \\ 
& \multicolumn{8}{l}{  Years collected : 2010--2011 } \\
& \multicolumn{8}{l}{ Source : Instituto Nacional de Estadistica e Infomatica (INEI)} \\
& \multicolumn{8}{l}{Encuesta Nacional de Hogares (ENAHO);  also used in \cite{Hanna:Olken}}  \\ 
\multicolumn{8}{l}{ Outcome distribution summary (unconditional) } \\ 
& & &  Q05  & Q50 & Q95 & Mean & Stdev & Var \\ 
\multicolumn{3}{l}{$y$ = Consumption $(\times 10^3 \ \text{Peruvian Soles})$ } &     0.1015   & 0.3550 &   1.1582 &    0.4608 &   0.4136 & 0.1711 \\
\multicolumn{3}{l}{ $y$ = $log$ Consumption (Peruvian Soles) }& 4.6254  &  5.8766   & 7.0516   & 5.8566   & 0.7442  & 0.5539\\

\hline 
\hline
 \end{tabular*}	
 
 %\

\flushleft
Asset-based proxy-means estimation results.   Estimates are presented for \textit{TBFMS I--TBFMS IV}, \textit{OLS} as described in the text.  The \textit{OLS} estimator replicates analysis in  \cite{Hanna:Olken}.

%This table presents full estimation results for all TBFMS procedures considered in estimation of Poverty-predictions.

 \end{table}

% Ordinary least squares using baseline indicators; identical to \cite{Hanna:Olken} (2)   Algorithm 2 with the tests defined in Definition 1--Definition with $c_\tau=1.01$, $\alpha = .05$, with initial inclusion of original baseline indicator variables as described in the introduction of Section 4.    

Proxy-means tests are estimated with ordinary least squares (OLS) using un-interacted indicator variables (replicating  \cite{Hanna:Olken}) and with four TBFMS estimators adding interactions.  The  TBFMS estimators are \textit{TBFMS I--IV} and use Algorithm 2 together with tests of Definitions 1--4, with selection initialized, as described in Section 4.3, with baseline variables $S_{\text{base}}$ consisting of original (un-interacted) indicators and a constant term.  Tuning parameters are set to $c_\tau = 1.01, \alpha = .05$ in all cases.  Including baseline indicators is natural because they are ex-ante researcher chosen as relevant and because they are few in number relative to the sample size. 10 baseline indicators of the 72 contained in the original data were excluded due to multicollinearity (after having generated interactions), leaving 62 baseline indicators (which includes a constant term).  Generated interaction terms with training sample standard deviation $<0.03$ (corresponding to 20 or less 1s) were discarded.  Interactions are residualized away from baseline indicators on the training sample\footnote{Residualizing does not involve the outcome variable and can be viewed similarly to the process of demeaning covariates provided the number of baseline indicators is sufficiently small.  Residualizing introduces a small amount of dependence across observations as does demeaning.  The analysis in the proof of Theorem 5 can be adjusted to show that dependence arising from residualization away from a fixed, suitably small set of covariates has negligible effects.}.  Residualizing only affects the terms $\hat \tau_{jS}$ in \textit{TBFMS I}, with $\hat \tau_{jS}$ typically becoming smaller, thus resulting in more selections into the model. Residualizing does not affect the estimation \textit{TBFMS II---IV}.   The final number of generated interactions is 23964 which exceeds the training sample size of 22674.   

Table 1  reports test sample mean square error (MSE), training sample MSE, relative and absolute reduction in test sample MSE (compared to OLS), p-values and confidence sets for absolute reduction in MSE based two-sided paired t-tests.  Table 1 also reports summary statistics.

Each TBFMS estimate indicates significant improvement in test sample MSE relative to OLS at the 5\% level. Highest relative improvement are with \textit{TBFMS III} and \textit{TBFMS IV}, which selected the same models both for $y=$ \textit{log} Consumption and $y=$ Consumption.   For $y=$ Consumption, in the test sample, \textit{TBFMS III, IV} give a 5.4865\% relative improvement to \textit{OLS} in MSE; 0.0049 absolute improvement, $p$-value = .0002 based on 2-sided paired $t$-test, 95\% CI (0.0023, 0.0076).  For  \textit{log} Consumption, \textit{TBFMS III, IV} give a 2.3809\% relative improvement to OLS in MSE, absolute 0.0045 improvement; $p$-value = 1.8{\tiny \textrm{E}}-10 based on 2-sided paired $t$-test, 95\% CI (0.0031, 0.0060).  
  
% Table 4 suggests the power and size properties of the tests associated to \textit{TBFMS I} may not match those of \textit{TBFMS II} in finite samples.  In particular, \textit{TBFMS I} selects only one additional covariate when estimating on outcome log consumption, and zero additional covariates when estimating on outcome consumption.  \textit{OLS+Ineraction} gives negative relative MSE for the log consumption outcome (relative difference -5.89\%, aboslute difference $-0.01$, $p$-value = 6.3{\tiny \textrm{E}}-10), suggesting  overfitting.  

\section{{\em TBFMS} Simulation Studies}

The results in the previous sections suggest that estimation with TBFMS should produce quality estimates in large sample sizes. This section conducts two simulation studies to evaluate the finite sample performance of TBFMS relative to select other procedures commonly used in high-dimensional regression problems in finite samples.  

The simulation study compares the following estimators.

\begin{itemize}
\item[1.] \textit{TBFMS I}.  Algorithm 2 with tests defined in Definition 1 with $c_\tau=1.01$, $\alpha = .05$.  
\item[2.] \textit{TBFMS II}.   Algorithm 2 with simplified tests defined in Definition 2 with $\alpha = .05$.
\item[3.] \textit{TBFMS III}.   Algorithm 2 with streamlined tests defined in Definition 3 with $\alpha = .05$.  
\item[4.] \textit{TBFMS IV}.   Algorithm 2 with homoskedastic tests defined in Definition 4 with  $\alpha = .05$.
 \item[5.] \textit{Lasso-CV}.  Standard Lasso, with penalty parameter chosen using 10-fold cross validation.   $\hat \theta$ is the minimizer of the Lasso objective function and is not refit on the selected model.
 \item[6.] \textit{Post-Het-Lasso}.  Lasso implementation in \cite{BellChenChernHans:nonGauss}, which is designed specifically to handle heteroskedastic disturbances. This requires two tuning parameters to $c_\tau$ and $\alpha$ set to $c_{\tau} = 1.01$ and $\alpha = .05$. The estimate $\hat \theta$ is refit on the selected model. 
 \item[7.] \textit{Oracle} Least squares regression on model consisting of $S_0 = \{j: [\theta_0]_j \neq 0 \}$.  
 \end{itemize}

\subsection{Simulation I}

\comment{
\begin{table}[h!] \caption 
 {Simulation I : {\em TBFMS} Simulation Design}
\begin{tabular}{ll}
\hline \hline \\ 
\text{Data: }& $\mathscr D_n = (y_i, x_i)_{i=1}^n \ iid $ \\
\text{DGP: }& $y_i = x_i\theta_0 + \eps_i $\\
&$p = \text{dim}(x_i) = 2n, \ \ \ \theta_{0j} = b_0^{j-1}\textbf{1}_{j \leq s_0}$\\
&$x_{ij} \sim N(0,1), \ \text{with corr}(x_{ij},x_{ik}) = 0.5^{|j-k|} $\\ 
&$\varepsilon_i \sim \mathbb \sigma_i N(0,1), \ \sigma_i = \text{exp}{(\rho_0 \textstyle \sum \nolimits_{j=1}^p 0.75^{(p-j)}x_{ij})}$ \\ 
%Settings: & I.i \ :  $s_0 =6$ , $b_0 = -0.5 $ , $ \rho_0 =0.5$ ,$n \in \{100,200,300,400,500 \} $ \\  
Settings: &  $s_0 =6$ , $b_0 \in\{ -0.5, 0.5\} $ , $ \rho_0 \in \{0, 0.5\}$ ,$n \in \{100,500 \} $ \\  \\
\hline
\end{tabular}
\end{table}
}

The first simulation study, Simulation I, evaluates TBFMS relative to several alternative estimators on a series of performance metrics for high-dimensional sparse linear regression problems.   Simulation I considers data of the form $\mathscr D_n = \{(y_i,x_i)\}_{i=1}^n$ with $y_i = x_i'\theta_0 + \eps_i$.   Samples $\mathscr D_n$ are drawn from several data-generating processes reported in Table 2.  Each considered data generating process is characterized by parameters $s_0, \rho_0, b_0, n$.
The parameter $s_0$ is the sparsity.  The parameter $b_0$ controls the nature of the coefficient vector by $\theta_{0j} = b_0^{j-1}\textbf{1}_{j\leq s_0}$.  When $b_0 <0$, the coefficients $\theta_{0j}$ alternate sign in $j$ and when $|b_0| < 1$, the coefficients decay. %, and when $b_0 = 0.5$, the coefficients are all positive.  The two different settings for $b_0$ create different interplay between the Toeplitz correlation structure in the covariates and their corresponding coefficients.  
The parameter $\rho_0$ controls the presence of heteroskedasticity in the disturbance terms $\eps_i$.  The terms $\eps_i$ are homoskedastic when $\rho_0=0$ and heteroskedastic otherwise.  The dimensionality is always taken to be double the sample size so that $p = 2n$.   Each simulation design is replicated 1000 times.

The simulation results for Simulation I are reported in Table 3.  The results track various measures of estimation quality for the five estimators for fixed $s_0, b_0,$ and $\rho_0$, and for $n$.  Columns in Table 3 display (1) \textit{MPEN}---mean prediction error norm defined as $\En[ (x_{i}'\theta_0 - x_i'\hat \theta)^2]^{1/2}$, (2) \textit{RMSE}---root mean square estimation error defined as $\|\hat \theta_2 - \theta_0\|_2$,  (3) \textit{MNCS}---mean number of correctly selected covariates from $S_0$, (4) \textit{MSSS}---mean size of selected set of covariates in all cases averaged over simulation replications. 

Table 3 indicates that predictive and estimation performance improve for all estimates as $n$ increase from 100 to 500.  \textit{Lasso-CV} selects an increasing number of variables.  Outside the \textit{Oracle} estimator, for all $n$, \textit{TBFMS II} attains the best predictive and estimation performance in settings with heteroskedasticity.  
In Table 3, there is no single feasible estimator that dominates in every setting in terms of estimation error or prediction error.  However, in all settings, \textit{Lasso-CV} selects the most covariates (both in absolute terms and in terms of the number of correctly identified covariates), followed by \textit{TBFMS I - IV} and \textit{Post-Het-Lasso.}  % There are instead important instances when each estimator performs better.  In both settings with positive coefficients ($b_0 = 0.5$), Lasso-CV achieves the smallest estimation error with TBFMS I and TBFMS II having slightly higher estimation error.  In these settings, however, the prediction error is smallest with TBFMS II.  
With alternating coefficients  ($b_0=-0.5$), %however, 
\textit{TBFMS I} and \textit{TBFMS II} dominate \textit{Lasso-CV} and \textit{Post-Het-Lasso} on prediction error and estimation error.  This suggests that the performance of these estimators depends on the configuration of the signal $\theta_0$ relative to the correlation structure of the covariates.   Finally, the relative difference in performance across estimators is larger in the presence of heteroskedasticity.  In the presence of heteroskedasticity, the \textit{Post-Het-Lasso} exhibits faster improvement in estimation error and prediction error with increasing $n$, though it is still dominated by the other estimators.  Note that each of the techniques, \textit{TBFMS I} and \textit{Post-Het-Lasso}, are theoretically valid for sequences of data-generating processes with heteroskedasticity.  In addition, the properties of cross-validation with Lasso are only beginning to be understood (see \cite{CV:Lasso} for analysis of Lasso with cross-validation).  But it is seen from this simulation study that \textit{Lasso-CV} leads to selection of substantially more covariates selected.% to the extent that the effects of heteroskedasticity on the performance of the estimator are still not fully clear.  

 \comment{
 \begin{figure}[H]
\caption{Simulation Results}

\hrule
  \includegraphics[scale=.19]{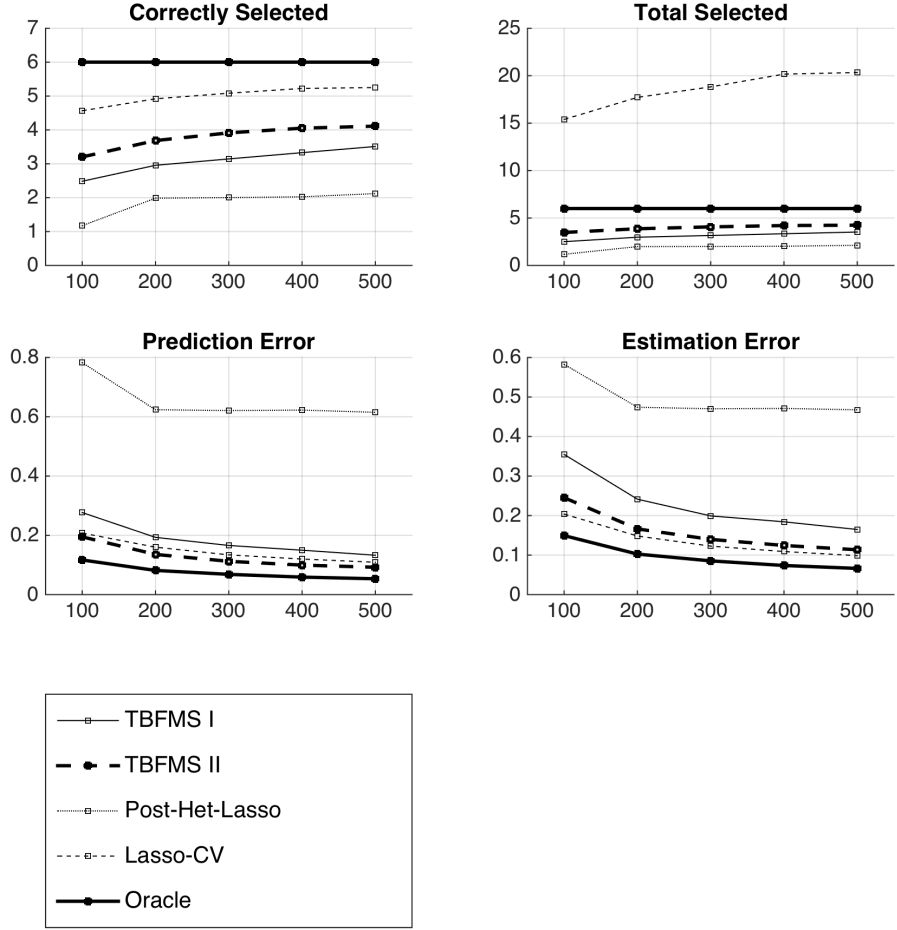}

\footnotesize
\flushleft
Simulation results \textit{TBFMS I-IV, Post-Het-Lasso, Lasso-CV}, and \textit{Oracle.}  Simulation design is described in Table 2, using settings 
$s_0 = 6$, $b_0 = -0.5$, $\rho_0 = 0.5$.
Plots based on 1000 simulation replications for every $n = 100,200,300,400,500$ indexed on the horizontal axis. Plots display (1) \textit{Correctly Selected}---number of correctly identified covariates from $S_0$, (2) \textit{Total Selected}---total number of selected covariates, (2) \textit{Prediction Error}---prediction error defined as $\En[ (x_{i}'\theta_0 - x_i'\hat \theta)^2]^{1/2}$, (3) \textit{Estimation Error}---estimation error defined as $\|\hat \theta_2 - \theta_0\|_2$, in all cases averaged over simulation replications.     

\

\hrule

 \end{figure}
}

%\bibliographystyle{plain}
%\bibliography{dkbib1}

%\pagebreak

%\appendix

\subsection{Simulation II}

This section performs a second simulation study, Simulation II, on the use of TBFMS in selecting control covariates for the estimation of the effect of a covariate of interest on an outcome from a large set of potential observable controls.

Simulation II considers data of the form $\mathscr D_n = \{(y_i,x_i,w_i)\}_{i=1}^n \sim {\Pr_n}$ where $y_i \in \mathbb R$ are outcome variables, $x_i\in \mathbb R$ are variables of interest, and $w_i\in \mathbb R^p$ are controls.   
In particular,
\begin{align*}
  y_{i}  &= x_{i}\beta_0 + w_i'\theta_0^{1} +  \eps_{i}\\
  x_{i}  &= w_i'\theta_0^{2}  + u_{i}
\end{align*}
for some parameters $\beta_0 \in \mathbb R, \theta_0^{1},\theta_0^{2}\in \mathbb R^p$ with $\Ep[\varepsilon_{i} | w_{i},x_{i}]= 0$ and $ \Ep[u_{i} | w_i] = 0$.  Here, the impact of the policy/treatment variable $x_i$ on the outcome $y_i$ is measured by the unknown parameter $\beta_0$ which is the target of inference.  The  $w_i$ are potentially important conditioning variables.  %Data are assumed independent but not necessarily identically distributed across $i$.
The confounding factors $w_{i}$ affect $x_i$ via the function $w_i'\theta_0^{2}$ and the outcome variable via the function $w_i'\theta_0^{1}$. Both of the parameters $\theta_0^{1}$ and $\theta_0^{2}$ are unknown.
As in Simulation I, specific data-generating processes used in Simulation II depend on parameters $s_0, b_0, \rho, n$ and are given by Table 2.

The structure of the TBFMS procedure applied to a specific linear regression problem ensures that any coefficient in that problem, unless reliably distinguishable from zero, is estimated to be exactly zero.  This property complicates inference after model selection in sparse models that may have a set of covariates with small but non-zero coefficients.  The use of TBFMS for a linear regression $y_i$ on $(x_i,z_i)$, possibly initializing the $\hat S$ to including $x_i$ at the start of model selection, may result in excluding important conditioning covariates, which may lead to non-negligible omitted variables bias of parameters of interest.  As a result, inference which does not take account the possibility of such model selection mistakes can be distorted.  This intuition is formally developed in %\citen{potscher} and 
\citen{leeb:potscher:pms}.  Several recent papers offer solutions to this problem; see, for example, \citen{ZhangZhang:CI}, \citen{BCH-PLM}, \citen{vdGBRD:AsymptoticConfidenceSets}; \citen{JM:ConfidenceIntervals}. %This section considers one technique for valid post-model selection inference, the Post-Double Selection technique (\cite{BCH-PLM}),  in conjunction with TBFMS.

%Inference about $\beta_0$ is impossible in this model without imposing further structure since $p > n$ elements in $w_{i}$ are allowed.  The additional structure is added by assuming that a sparsity condition applies to both $\theta_0^{\text{B}_1}$ and $\theta_0^{\text{B}_2}$.  Once these assumptions are in place, the exogeneity of $x_{i}$ may be taken as given after controlling linearly for a relatively small number, $s_0 < n$, of the variables in $w_{i}$ whose identities are \textit{a priori} unknown.  Specifically, impose the following restrictions.

To estimate $\beta_0$ in this environment, adopt the \textit{post-double-selection} method of \citen{BCH-PLM} in conjunction with TBFMS.
This method proceeds by first substituting to obtain predictive relationships for the outcome $y_{i}$ and the treatment $x_{i}$ in terms of only control variables:
\begin{align*}
 y_{i}  &=  w_{i}'\theta_0^{\text{RF}} +  v_{i}\\
 x_{i}  &=  w_{i}'\theta_0^{\text{FS}} +  u_{i}
\end{align*}
with $\theta_0^{\text{FS}} = \theta_0^{2}$, referred to as the \textit{first stage} (FS) coefficient, and $\theta_0^{\text{RF}} = \theta_0^{1} + \beta_0 \theta_0^{2}$, referred to as the \textit{reduced form} (RF) coefficient.  TBFMS is applied to each of the above two equations to select one set of variables that are useful for predicting $ y_{i}$ and another set of variables useful for predicting $x_i$.  Once this is done, the union of the selected sets will index the final set of control variables.
\citen{BCH-PLM} develop and discuss the post-double-selection method in detail.  %They note that including the union of the variables selected in each variable selection step helps address the issue that model selection is inherently prone to errors unless stringent assumptions are made.  %As noted by \citen{leeb:potscher:pms}, the possibility of model selection mistakes precludes the possibility of valid post-model-selection inference based on a single Lasso regression within a large class of interesting models.  
% The chief difficulty arises with covariates whose effects in (\ref{LinearFE1}) are small enough that the variables are likely to be missed if only (\ref{LinearFE1}) is considered but have large effects in (\ref{LinearFE2}).  The exclusion of such variables may lead to substantial omitted variables bias if they are excluded which is likely if variables are selected using only (\ref{LinearFE1}).\footnote{The argument is identical if only (\ref{LinearFE2}) is used for variable selection exchanging the roles of (\ref{LinearFE1}) and (\ref{LinearFE2}).  The argument also holds if considering only one of (\ref{PLMReducedForm}) or (\ref{PLMFirstStage}) for variable selection.} 
Using two model selection steps guards against distorted inference and guarantees that variables excluded in both model selection steps have a negligible contribution to omitted variables bias.

Analogously to Simulation Study I, Seven estimators are considered for estimation and are named \textit{TBFMS I --TBFMS IV, Lasso-CV, Post-Het-Lasso}, and \textit{Oracle}.  The estimators differ only in that they replace \textit{TBFMS I} with a different model selection technique in selecting covariates into the final estimated model as described above.  Final variance estimates $\hat V$ for $\hat \beta$ are based on HC3 standard errors \cite{white:het}. For each estimator and simulation setting, the bias and standard deviation of the point estimates, and coverage probability and average interval length are computed over 1000 simulation replications. 
Results are shown in Table 4.

The simulation results indicate that across the data-generating processes considered, TBFMS estimators largely achieve near-\textit{Oracle} coverage probabilities. \textit{TBFMS IV}-based estimates exhibit some size distortion in heteroskedastic settings.  In all settings, bias, standard deviation, and interval lengths of \textit{TBFMS I -- III} closely resemble the \textit{Oracle} estimator.   In some simulations (notably in Panel A) there is a large difference in coverage probabilities between the TBFMS estimates and the \textit{Post-Het-Lasso} estimate.  Though the \textit{Post-Het-Lasso}-based confidence sets are asymptotically uniformly valid, and are theoretically robust against model selection mistakes, finite sample model selection properties remain important.  % In this example, the signal-to-noise ratio in the first stage is low, which is a setting in which detection of important variables is more difficult for Post-Het-Lasso relative to TBFMS.  
Interestingly, in this case, using the relaxed penalty level with \textit{Lasso-CV} does not improve coverage probability.

\begin{table}[H]\caption
{Design for Simulation Studies I and II}
\begin{tabular}{lllllll}
\hline \hline  \\
 \multicolumn{2}{l}{\textbf{I. High-Dimensional Prediction}.  }& \multicolumn{2}{l}{\textbf{II.  High-Dimensional Controls}.  }\\  \\
\text{Data: }& $\mathscr D_n = (y_i, x_i)_{i=1}^n \ iid $  &  Data: &$\mathscr D_n = (y_i, x_i, w_i)_{i=1}^n \ iid $\\
\text{DGP: }& $y_i = x_i\theta_0 + \eps_i $ & DGP: &$y_i = x_i\beta_0 +w_i'\theta_0^{\text{1}}+ \eps_i $, \ \ $x_i = w_i'\theta_0^{2} + u_i $ \\
&$p = \text{dim}(x_i) = 2n$ &&$p = \text{dim}(w_i) = 2n $ \\ & $\theta_{0j} = b_0^{j-1}\textbf{1}_{j \leq s_0}$  & &  $\theta_{0j}^{1} = b_0^{j-1}\textbf{1}_{j \leq s_0}$, \ $ \theta_{0j}^{2} = \sin(j) \textbf{1}_{j \leq s_0}$ \\
&$x_{ij} \sim N(0,1), \ \text{corr}(x_{ij},x_{ik}) = 0.5^{|j-k|} $ &&$w_{ij} \sim \textrm{N}(0,1)$, \  $\text{corr}(w_{ij},w_{ik}) = 0.5^{|j-k|} $ \\
&$\varepsilon_i \sim \mathbb \sigma_i N(0,1)$ &&$(\varepsilon_i,u_i) \sim \mathbb \sigma_i \textrm{N}\(0,\( \begin{array}{cc} 1&0 \\ 0&1 \end{array}\) \)$ \\&  $ \ \sigma_i = \text{exp}{(\rho_0 \textstyle \sum \nolimits_{j=1}^p 0.75^{(p-j)}x_{ij})}$ &&$\sigma_i = \text{exp}{(\rho_0 \textstyle \sum \nolimits_{j=1}^p 0.75^{(p-j)}w_{ij})}$\\
\text{Settings: } &$s_0 = 6$, $ b_0 \in \{-0.5,0.5\}$,  $\rho_0 \in \{0,0.5\}$ & \text{Settings: } &$s_0 = 6$, $ b_0 \in \{-0.5,0.5\}$,  $\rho_0 \in \{0,0.5\}$ \\
& $n \in \{100,500\}$ & & $n \in \{100,500\} $\\ \\
 \hline 
\end{tabular}
\end{table}

\begin{table}[H] \caption 
 {Simulation I Results: Prediction in the Linear Model  } 
\begin{tabular*}{\textwidth}{p{2.2cm} p{.9cm}  p{.9cm} p{.9cm} p{.9cm} p{.6cm} p{.9cm}  p{.9cm}  p{.9cm} p{.9cm} } 
\hline          \hline                                                                  \\ 
&			&  \multicolumn{2}{c}{$n=100$}	&		& & 		& \multicolumn{2}{c}{$n=500$}		&		\\   \cline{2-5} \cline{7-10} \\  
&	\textcolor{white}{\Big |}MPEN		& RMSE		&	MNCS	&	MSSS	& & MPEN		& RMSE		&	MNCS	&	MSSS  \\   
  \\  \cline{2-10}     & \multicolumn{9}{c}{A. $\rho_0 = 0 $ : Homoskedastic,  $s_0 = 6$ : High Sparsity, $b_0 = -0.5$ : Alternating Coefficients   }\\ \cline{2-10}  \\ 
\textit{TBFMS I }  &    0.313  &  0.484  &  1.470  &  1.470	&&    0.146 & 0.218  & 2.599  & 2.599 \\
\textit{TBFMS II  } &   0.192  & 0.281  & 2.330 & 2.444	&&    0.094  & 0.132 &  3.403 &  3.478\\
\textit{TBFMS III } &   0.194  & 0.286  &  2.307  & 2.307	&& 0.094  & 0.133 &  3.432  & 3.432  \\
\textit{TBFMS IV }  &   0.191  &  0.281  & 2.301  &  2.358	&&   0.093  & 0.132 &  3.400  & 3.461\\
\textit{Post-Lasso}&   0.417  &  0.640  &  0.958  &  0.958	&&   0.401  & 0.617 &  1.000 & 1.000\\
\textit{Lasso-CV  }    &  0.285  &  0.468  &  2.850  &  21.463 &&      0.173  & 0.275  &  3.693  &40.617\\
\textit{Oracle   }    &  0.117  &  0.149  &  6.000  & 6.000	&&   0.053  & 0.066 &  6.000 & 6.000\\
 \\   \cline{2-10}     & \multicolumn{9}{c}{B. $\rho_0 = 0.5$ Heteroskedastic, $s_0 = 6$ : High Sparsity, $b_0 = -0.5$ : Alternating Coefficients   }\\ \cline{2-10}  \\ 
\textit{TBFMS I  }  &  0.507  &0.720  &  0.815  & 0.815	&    &  0.395  & 0.593  &1.077  & 1.077\\
\textit{TBFMS II }  &   0.478 &0.671  & 0.994  &1.012	&& 0.276  &0.401  &1.765  & 1.780\\
\textit{TBFMS III  }&   0.526  &0.734  & 0.824  & 0.824	&& 0.294  & 0.429 & 1.650  & 1.650\\
\textit{TBFMS IV }  &   0.529  & 0.726  & 0.896 & 0.973	&& 0.299  & 0.425 &1.713  & 1.807\\
\textit{Post-Het-Lasso }&    0.884  & 1.147  &  0.015  & 0.015	&&   0.581&0.846  & 0.867  & 0.867\\
\textit{Lasso-CV  }   &  0.594  & 0.849  &  1.242  & 8.863	&&   0.438  & 0.661  &1.670 &18.326\\
\textit{Oracle   }    & 0.379  & 0.482  & 6.000 & 6.000	&& 0.180 &0.222 &6.000 & 6.000\\
 \\   \cline{2-10}     & \multicolumn{9}{c}{C. $\rho_0 = 0 $ : Homoskedastic, $s_0 = 6$ : High Sparsity, $b_0 = 0.5$ : Positive Coefficients   }\\ \cline{2-10} \\  
\textit{TBFMS I  }  & 0.307  & 0.395  & 2.331  & 2.331	&&     0.147  &  0.186 & 3.396  & 3.396\\
\textit{TBFMS II  } &   0.193  & 0.244  & 3.094  & 3.193	&&  0.091  &  0.113  & 4.078  & 4.152\\
\textit{TBFMS III } &   0.193  & 0.244  & 3.064  & 3.064	&&  0.091  &  0.114  & 4.105  & 4.105\\
\textit{TBFMS IV }  &   0.191  & 0.241  &3.062  & 3.109&&0.091  &  0.114  & 4.068  &4.128\\
\textit{Post-Het-Lasso }&  0.782  & 0.583  &  1.174  & 1.174	&&  0.615  &0.468  & 2.121  &2.121\\
\textit{Lasso-CV }    &   0.207 & 0.204  & 4.570  & 15.392	&&     0.109  &  0.099  & 5.257  &  20.100\\
\textit{Oracle    }   &   0.117  & 0.149  & 6.000  & 6.000	&& 0.053  & 0.066  &6.000  &6.000\\
  \\  \cline{2-10}     & \multicolumn{9}{c}{D. $\rho_0 = 0.5 $ : Heteroskedastic, $s_0 = 6$ : High Sparsity, $b_0 = 0.5$ : Positive Coefficients  }\\ \cline{2-10} \\  
\textit{TBFMS I  }  &   0.665  & 0.789 & 1.274  & 1.274	&&      0.405   & 0.499  & 1.976 & 1.976\\
\textit{TBFMS II  } &   0.513  & 0.617  &1.759  &1.780	&&     0.285   & 0.346  & 2.580  &2.590\\
\textit{TBFMS III } &   0.577  &0.662  & 1.505  & 1.505	&&   0.310   & 0.369  &  2.400 & 2.400\\
\textit{TBFMS IV }  &   0.580  &  0.670  & 1.574  &1.656	&&     0.314   & 0.373  & 2.478  & 2.574\\
\textit{Post-Het-Lasso }&     1.314  & 1.014  &  0.288  & 0.288	&& 0.879   & 0.663  & 1.898  & 1.898\\
\textit{Lasso-CV }   &    0.570  & 0.522  & 3.125  & 11.661	&&  0.334   & 0.289  & 3.873  & 16.337\\
\textit{Oracle }     &    0.379  & 0.482  & 6.000  & 6.000	&&  0.180   & 0.222  & 6.000  & 6.000\\
 \\ \cline{1-10}	
 \end{tabular*}

 \
 
 \flushleft
%This table presents simulation results for the estimators \textit{TBFMS I--IV, Post-Het-Lasso, Lasso-CV,} and \textit{Oracle.}  The simulation design is described fully in Table 2. Tables are based on 1000 simulation replications for every $n = 100,500$.  Columns display (1) \textit{Correctly Selected}---number of correctly identified covariates from $S_0$, (2) \textit{Total Selected}---total number of selected covariates, (2) \textit{Prediction Error}---prediction error defined as $\En[ (x_{i}'\theta_0 - x_i'\hat \theta)^2]^{1/2}$, (3) \textit{Estimation Error}---estimation error defined as $\|\hat \theta_2 - \theta_0\|_2$, in all cases averaged over simulation replications.     

Simulation results for estimation in the design described in Table 2 with
 $s_0 =6$, 
 $b_0 \in \{-0.5,0.5 \}$,
$ \rho_0 \in \{0, 0.5\}$, and $n \in \{100,500 \} 
$.  Estimates are presented for the estimators, \textit{TBFMS I-- TBFMS IV, Post-Het-Lasso, Lasso-CV}, and \textit{Oracle} described in the text. Columns display (1) \textit{MPEN}---mean prediction error norm defined as $\En[ (x_{i}'\theta_0 - x_i'\hat \theta)^2]^{1/2}$, (2) \textit{RMSE}---root mean square estimation error defined as $\|\hat \theta_2 - \theta_0\|_2$,  (3) \textit{MNCS}---mean number of correctly selected covariates from $S_0$, (4) \textit{MSSS}---mean size of selected set of covariates in all cases averaged over simulation replications.  Figures are based on 1000 simulation replications.

 \end{table}

\begin{table}[H]\caption 
 {Simulation II Results: Control Selection in the Linear Model   } 
\begin{tabular*}{\textwidth}{p{2.2cm} p{.9cm}  p{.9cm} p{.9cm} p{.9cm} p{.6cm} p{.9cm}  p{.9cm}  p{.9cm} p{.9cm} } 
\hline          \hline                                                                  \\ 
&			&  \multicolumn{2}{c}{$n=100$}	&		& & 		& \multicolumn{2}{c}{$n=500$}		&		\\   \cline{2-5} \cline{7-10} \\  
&	\textcolor{white}{\Big |}  Bias		& StDev		&	Length	&	Cover	& & Bias		& StDev		&	Length	&	Cover  \\   
  \\  \cline{2-10}     & \multicolumn{9}{c}{D. $\rho_0 = 0 $ : Homoskedastic, $s_0 = 6$ : High Sparsity, $b_0 = -0.5$ : Alternating Coefficients  }\\ \cline{2-10} \\  
\textit{ TBFMS I     }     & -0.035 & 0.129 &  0.354 & 0.772 && -0.036 & 0.060 &  0.182 & 0.797 \\ 
\textit{ TBFMS II    }     & -0.009 & 0.114 &  0.428 & 0.934 && 0.000 & 0.047 &  0.178 & 0.941 \\ 
\textit{ TBFMS III    }    & -0.016 & 0.112 &  0.424 & 0.928 && -0.000 & 0.047 &  0.178 & 0.941 \\ 
\textit{ TBFMS IIV    }    & -0.011 & 0.112 &  0.425 & 0.929 && -0.000 & 0.047 &  0.178 & 0.940 \\ 
\textit{ Post-Het-Lasso}   & -0.190 & 0.057 &  0.210 & 0.068 && -0.193 & 0.025 &  0.094 & 0.000 \\ 
\textit{ Lasso CV       }  & -0.193 & 0.054 &  0.214 & 0.067 && -0.193 & 0.025 &  0.094 & 0.000 \\ 
\textit{ Oracle          } & 0.002 & 0.105 &  0.424 & 0.956 && -0.001 & 0.047 &  0.178 & 0.941 \\ 
 \\   \cline{2-10}     & \multicolumn{9}{c}{B. $\rho_0 = 0.5$ :  Heteroskedastic, $s_0 = 6$ : High Sparsity, $b_0 = -0.5$ :  Alternating Coefficients  }\\ \cline{2-10}  \\ 
 \textit{ TBFMS I      }    & 0.032 & 0.326 &  1.016 & 0.881 && 0.003 & 0.219 &  0.752 & 0.927 \\ 
\textit{ TBFMS II      }   & 0.010 & 0.373 &  1.179 & 0.901 && 0.005 & 0.231 &  0.806 & 0.929 \\ 
\textit{ TBFMS III     }   & 0.010 & 0.366 &  1.141 & 0.897 && 0.004 & 0.231 &  0.801 & 0.928 \\ 
\textit{ TBFMS IV     }   & 0.009 & 0.372 &  1.167 & 0.899 && 0.005 & 0.231 &  0.803 & 0.929 \\ 
\textit{ Post-Het-Lasso }  & -0.111 & 0.183 &  0.553 & 0.721 && -0.104 & 0.129 &  0.421 & 0.650 \\ 
\textit{ Lasso CV      }   & -0.116 & 0.209 &  0.637 & 0.745 && -0.104 & 0.130 &  0.422 & 0.651 \\ 
\textit{ Oracle         }  & 0.004 & 0.406 &  1.333 & 0.909 && 0.009 & 0.234 &  0.820 & 0.934 \\ 
  \\  \cline{2-10}     & \multicolumn{9}{c}{C. $\rho_0 = 0 $ : Homoskedastic, $s_0 = 6$ : High Sparsity, $b_0 = 0.5$ : Positive Coefficients  }\\ \cline{2-10} \\  
\textit{  TBFMS I      }    & 0.012 & 0.122 &  0.362 & 0.814 && -0.054 & 0.059 &  0.184 & 0.729 \\ 
\textit{ TBFMS II      }   & -0.012 & 0.113 &  0.426 & 0.929 && 0.000 & 0.044 &  0.179 & 0.960 \\ 
\textit{ TBFMS III     }   & -0.019 & 0.112 &  0.423 & 0.926 && -0.000 & 0.044 &  0.179 & 0.961 \\ 
\textit{ TBFMS IV     }   & -0.014 & 0.113 &  0.423 & 0.920 && 0.000 & 0.044 &  0.179 & 0.960 \\ 
\textit{ Post-Het-Lasso }  & -0.088 & 0.061 &  0.234 & 0.669 && -0.085 & 0.027 &  0.102 & 0.109 \\ 
\textit{ Lasso CV      }   & -0.086 & 0.059 &  0.234 & 0.680 && -0.085 & 0.027 &  0.102 & 0.109 \\ 
\textit{ Oracle       }    & 0.002 & 0.102 &  0.420 & 0.954 && -0.001 & 0.044 &  0.179 & 0.960 \\ 
  \\  \cline{2-10}     & \multicolumn{9}{c}{D. $\rho_0 = 0.5 $ : Heteroskedastic, $s_0 = 6$ : High Sparsity, $b_0 = 0.5$ : Positive Coefficients  }\\ \cline{2-10} \\  
\textit{  TBFMS I    }      & -0.032 & 0.334 &  1.020 & 0.857 && 0.009 & 0.250 &  0.763 & 0.918 \\ 
\textit{ TBFMS II    }     & -0.016 & 0.371 &  1.196 & 0.882 && -0.001 & 0.259 &  0.811 & 0.923 \\ 
\textit{ TBFMS III    }    & -0.010 & 0.362 &  1.153 & 0.884 && -0.001 & 0.257 &  0.807 & 0.920 \\ 
\textit{ TBFMS IV    }    & -0.007 & 0.363 &  1.176 & 0.885 && 0.000 & 0.259 &  0.809 & 0.920 \\ 
\textit{ Post-Het-Lasso  } & -0.161 & 0.216 &  0.596 & 0.638 && -0.044 & 0.154 &  0.427 & 0.859 \\ 
\textit{ Lasso CV      }   & -0.063 & 0.203 &  0.641 & 0.903 && -0.044 & 0.154 &  0.427 & 0.858 \\ 
\textit{ Oracle       }    & -0.019 & 0.397 &  1.352 & 0.895 && 0.003 & 0.263 &  0.826 & 0.924 \\ 
 \\ \cline{1-10}	
 \end{tabular*}

  \

 \flushleft
 Simulation results for estimation in the design described in Table 2 with
 $s_0 =6$, 
 $b_0 \in \{-0.5,0.5 \}$,
$ \rho_0 \in \{0, 0.5\}$, and $n \in \{100,500 \} 
$.  Estimates are presented for the estimators, \textit{TBFMS I-- TBFMS IV, Post-Het-Lasso, Lasso-CV}, and \textit{Oracle} described in the text. Columns display (1)\textit{Bias}---bias of the respective estimates for $\beta_0$, (2) 'textit{StDev}---standard deviation of the respective estimates for $\beta_0$, (3) \textit{Length}---length of confidence intervals for $\beta_0$, (4) \textit{Cover}---coverage probabilities of the respective 95\% confidence intervals for $\beta_0$.  Figures are based on 1000 simulation replications. 

\

\

\

\

\

 \end{table}

\section{Conclusion} This paper has considered TBFMS for high-dimensional sparse linear regression problems.  The procedure is shown to achieve estimation rates matching those of Lasso and Post-Lasso under a broad class of data-generating processes.  

\appendix
\section{Proofs}
This appendix proves Theorems 1 and 4.   As TBFMS is a greedy procedure which is not the resulting solution of a simple optimization problem, the proofs establishing the properties of TBFMS cannot refer to any global optimality conditions.   This fact limits the applicability of common m-estimation arguments or arguments for similar bounds for Lasso, and requires the development of certain different techniques.

In the course of the proofs,
several important results along the way are recorded as lemmas.  Lemmas which do not follow immediately from arguments in this section are proven in the online supplemental material, Supplement to ``Analysis of Testing-Based Forward Model Selection.''

In addition, Theorems 2, 3 and 5 are proven in the online supplemental material.

\subsection{Proof of the First Statement of Theorem 1}

The first statement of Theorem 1 is proven by creating an appropriate analogue of the \textit{basic inequality}\footnote{For Lasso estimation with penalty level $\lambda$, the basic inequality asserts that $\ell(\hat \theta) + \lambda \| \hat \theta \|_1 \leq \ell(\theta_0) + \lambda \| \theta_0 \|_1$.} from standard Lasso analysis. Specifically, the following lemma holds. \begin{lemma} $\ell(\hat S \cup S_0) \leq \ell( \theta_0) $.\end{lemma}  %To see this, note that the facts $\ell(\hat S) = \ell(\hat S \cup S_0) + [\ell(\hat S) - \ell(\hat S \cup S_0)]$ and $\ell(\hat S) = \ell(\thetahat)$ together imply $\ell(\hat S \cup S_0) \leq \ell(\theta_0)$.  %, thus giving an analogue to the \textit{basic inequality} used in Lasso arguments (see e.g. \cite{BickelRitovTsybakov2009}.) 

Lemma 1 holds by $\ell(\hat S \cup S_0) \leq \ell(S_0) \leq \ell(\theta_0)$.  Once the analogue basic inequality is noted, $\ell(\hat \theta)$ can be related to $\ell(\theta_0)$ with a bound that depends on $s_0, t,$ and $ \varphi_{\min}(\hat s + s_0)(G)^{-1}$.  The following is
Lemma 3.3 in \cite{Submodular:Spectral}.  

\begin{lemma}[\cite{Submodular:Spectral}] $\ell(\hat S) - \ell(\hat S \cup S_0) \leq \varphi_{\min}(\hat s + s_0)(G)^{-1} \sum_{j \in  S_0 \setminus \hat S} (-\Delta_j \ell(\hat S))$.  \end{lemma} 
Using the fact that $t$ is the threshold in Algorithm 1, and thus $-\Delta_j \ell(\hat S) \leq t$, gives the further bound $\ell(\hat S) - \ell(\hat S \cup S_0) \leq s_0t\varphi_{\min}(\hat s+ s_0)(G)^{-1}.$  
\noindent Applying the basic inequality along with the fact that $\ell(\hat S) = \ell(\hat \theta)$ gives
$\ell(\hat \theta) \leq \ell(\theta_0) + s_0t\varphi_{\min}(\hat s + s_0)(G)^{-1}.$ 
Expanding the quadratics $\ell(\hat \theta)$ and $\ell(\theta_0)$, and applying arguments analogous to those in Lemma 6 in \cite{BellChenChernHans:nonGauss}, gives
\begin{lemma} $\En[(x_i ' \theta_0 - x_i' \hat \theta)^2]^{1/2} \leq c_{\text{\em F}}(\hat s)$.\end{lemma} 
\noindent Lemma 3 is the first statement of Theorem 1.   

\subsection{Proof of First and Third Statements of Theorem 4}
Let $\mathscr T$ be the event implied by Condition 2.  Then  $\Pr(\mathscr T) \geq 1 - \alpha -3 \delta_{\text{test}}/3 = 1- \alpha - \delta_{\text{test}}$.  The rest of the proof works on the event $\mathscr T$.
If Algorithm 2 terminates at $K_{\text{test}}$ steps or earlier, then it terminates at a step with $-\Delta_j\E(\hat S)  \leq  c_{\text{test}}$
for every $j \notin \hat S$. Similarly to the proof of Theorem 1 above, Lemma 3.3 in \cite{Submodular:Spectral} yields
$ |\E(S_0) - \E(\hat S) |\leq  \varphi_{\min}(K_{\text{test}})(\Ep[G])^{-1} \sum_{j \in S_0 \setminus \hat S} -\Delta_j \E(S) \leq s_0 c_{\text{test}}  \varphi_{\min}(K_{\text{test}})(\Ep[G])^{-1} .$
It is shown in the next section that $\hat s \leq K_{\text{test}}$ on $\mathscr T$, completing the proof of the first statement of Theorem 4.  Next, \begin{lemma} $\En[(x_i ' \theta_0 - x_i' \hat \theta)^2]^{1/2} \leq c_{\text{\em T}}''(\hat s)$. \end{lemma} 
\noindent Lemma 4 is the third statement of Theorem 4.     

\subsection{Proof of Sparsity Bounds for Theorems 1 and 4}

The sparsity bounds in Theorems 1 and 4 are proven together.  In the case of Theorem 1, the covariates $x_j = (x_{1j},...,x_{nj})'$, outcome $y = (y_1,...,y_n)'$ as well as disturbances $\eps = (\eps_1,...,\eps_n)'$ are considered elements of the Hilbert space $\mathbb R^n$ with inner product $\langle a,b \rangle =n\times  \En[a_ib_i]$.  In the case of Theorem 4, $x_j, y,\eps$ are elements of the Hilbert space $L^2(\Omega, \mathbb R^n)$, of $\Pr$-square-integrable random vectors taking values in $\mathbb R^n$ with $\Omega$ an underlying probability space, with inner product $\langle a, b \rangle_{L^2(\Omega, \mathbb R^n)} = n \times \Ep \En[a_ib_i]$.   The notation $\mathsf H$ is used to denote the appropriate of these Hilbert spaces according to the cases of $\mathsf H = \mathbb R^n$ for Theorem 1 and $\H = L^2(\Omega, \mathbb R^n)$ for Theorem 4.  In addition, let 
$t_{\mathsf H} = t$ and $G_{\mathsf H} = G $ in the case of Theorem 1 and let $ t_{\H} = c_{\text{test}}'$ and $ G_\H = \Ep [G]$ in the case of Theorem 4.  In the case of Theorem 4, the arguments that follow hold on $\mathscr T$, the event defined by Condition 2 (see also the previous subsection).

\subsubsection{Two Orthogonalizations} 

 Let $v_k\in \H$, $k=1,...,s_0$ denote \textit{true covariates} which refer to $x_j$  for $j \in S_0$. The term \textit{false covariates} refers to those $x_j$ for which $j \notin S_0$.  
Consider \textcolor{black}{the step} after which there are exactly $m$ false covariates selected into the model.  These are denoted $w_1,...,w_m$, ordered according to the order they were selected, and indexed by the set $A = \{1,...,m\}$.    

Apply Gram-Schmidt orthogonalization to $v_1,...,v_{s_0}, \eps$ with the inner product from $\mathsf H$.  The ordering is according to selection into $\hat S$.  Any true covariates unselected at $m$ false covariate selections are temporarily ordered arbitrarily after the selected true covariates. $\eps$ is placed last.  This yields a new set of unit-norm elements $$v_1,...,v_{s_0}, \eps \mapsto \tilde v_1,...,\tilde v_{s_0}, \tilde \eps \in \mathsf H.$$    This orthogonalization also yields parameters $$\tilde \theta = (\tilde \theta_1,...\tilde \theta_{s_0})' \in \mathbb R^{s_0}, \ \tilde \theta_{\tilde \eps} \in \mathbb R \ \ \text{satisfying} \ \ y = \tilde v_1 \tilde \theta_1 +...+ \tilde v_{s_0} \tilde \theta_{s_0} + \tilde \theta_{\tilde \eps }\tilde \eps.$$
Reorder the unselected covariates $v_k$ such that for any unselected true covariate, $\tilde \theta_k > \tilde \theta_l$ whenever $l>k$.  No additional orthogonalization is performed.

Apply a separate orthogonalization to $w_1,...,w_m\in \H$.  These are orthogonalized by the Gram-Schmidt process according to order of inclusion into $\hat S$, with true covariates included (interspersed according to when they were selected between the $w_j$) in the orthogonlization process.   The resulting Gram-Schmidt-orthogonalized elements are renormalized to give $$ w_1,..., w_m ; v_1,...,v_{s_0} \mapsto \tilde w_1,...,\tilde w_m \in \H$$ such that the component of each $\tilde w_j$ orthogonal to the span of $\tilde v_1,...,\tilde v_{s_0}$ in $\mathsf H$ has unit norm.    This renormalization is possible whenever $\varphi_{\min}(m+s_0)(G_\H) > 0$.
Therefore, $\tilde w_j$ can be decomposed into $\tilde w_j = \tilde r_j + \tilde u_j$ which satisfy $\tilde r_j \in \text{span}( \tilde v_1,...,\tilde v_{s_0})$ and  $\tilde u_j \in \text{span}(\tilde v_1,...,\tilde v_{s_0})^\perp$ and $\| \tilde u_j \|_\H = 1$. 

 For ease of reading, the online supplemental material also presents additional descriptive notation and details about the orthogonalization constructions.% \textcolor{black}{If $\mathscr Q_{S_{\text{pre-}w_j}}w_j \in \text{span}{(\tilde V)}$ then set $c_j$ such that $\langle \tilde w_j , \tilde w_j \rangle= 1$.} \textcolor{black}{Note again that the process for defining $\tilde w_j$ is the Gram-Schmidt process, only this time with a different normalization}.

The two orthogonalizations are next related to each other by associating, for each $j \in A$, parameters $$j \mapsto ( \tilde \gamma_j \in \mathbb R^{s_0}, \tilde \gamma_{j\tilde \eps} \in \mathbb R)$$
which  are defined as 
%projection coefficients (ie. regression coefficients) from projecting $\tilde w_j$ on $\tilde v_1,...,\tilde v_{s_0}, \tilde \eps$.  %They therefore measure dependence between $w_j$ and appropriate constituents of $y$.    
%Because $\tilde v_1,...,\tilde v_{s_0}, \tilde \eps$ are orthonormal, 
$\tilde \gamma_{jk} = \langle \tilde w_j, \tilde v_k \rangle_{\mathsf H}, \tilde \gamma_{j \tilde \eps} = \langle \tilde w_j, \tilde \eps \rangle_\H$.  Assume without loss of generality that each component of $\tilde \theta$ is positive (the remainder of the proof does not depend on the sign assigned during orthogonalization).  Similarly, assume without loss of generality that $\tilde \gamma_{j}'\tilde \theta \geq 0$.

A large part of the following analysis is at the level of the parameters  $\tilde \theta, \tilde \theta_{\tilde \eps}, \tilde \gamma_j, \tilde \gamma_{j\tilde \eps}$.  Therefore, some remarks are helpful.  Aside from relating the two orthogonalizations, these parameters also encode information about incremental loss for various $j$ and sets $S$.
For example, in the case of Theorem 1, it can be shown that 
$-\Delta_j \ell(S)   = \frac{1}{n} \frac{1}{\| \tilde w_j \|_{\mathsf H}^2} (\tilde \theta ' \tilde \gamma_j + \tilde \theta_{\tilde \eps}\tilde \gamma_{j\tilde \eps})^2$ if $S$ is the set of all covariates selected into $\hat S$ before $w_j$.  Similarly, $-\Delta_k\ell(S) = \frac{1}{n} \tilde \theta_k^2$ if $S$ corresponds to the set $\{v_1,...,v_{k-1}\}$.  Relatedly, note that if $\tilde \gamma_j'\tilde \theta$ is large for sufficiently many $j$, then some dependence between $\tilde \gamma_j$ may be anticipated; see \cite{tao:blog:transitivity}, for a general discussion of partial transitivity of correlation.  This, however, heuristically creates tension with the fact that $\tilde \gamma_j$ arise from orthogonalized covariates. 

%Let $\mathscr M_{k}$ denote projection in $\mathsf H$ onto the space orthogonal to $\text{span}(\{v_1,...,v_k\})$.  Then note that $\tilde v_k = \frac{ \mathscr M_{k-1}v_k }{ \langle v_k,\mathscr M_{k-1}v_k\rangle^{1/2}} \ \text{for} \ k = 1,...,s_0$ and $\tilde \eps =  \frac{\mathscr M_{s_0 }\eps}{(\langle \eps,\mathscr M_{s_0}\eps \rangle)^{1/2}}.$  For more general sets $S$, let $\mathscr Q_{S}$ be projection onto the space orthogonal to $\text{span}(   \{   x_j, \ j \in S \}  )$.  For each selected covariate, $w_j$, set  $S_{\text{pre-}w_j}$ to be the set of (both true and false) covariates selected prior to $w_j$ and note that $\tilde w_j = c_j \mathscr Q_{\text{pre-}j} w_j$ for a suitable normalizing constant $c_j$.  

\subsubsection{Main Sparsity Bounds}

Divide the set $A$ of false covariates into two sets $A_1$ and $A_2$, with cardinalities $m_1$ and $m_2$, on the basis of the magnitude of $\tilde \gamma_{j \tilde \eps}$.  Set  $A_1 = \left \{ j : | \tilde \gamma_{j \tilde \eps} |^2 \leq  \frac{t_\H n}{ \textcolor{black}{ 3} \| \tilde \eps \|_{\mathsf H} }  \right \}, \ A_2 = A \setminus A_1$.  Note that large values of $\tilde \gamma_{j \tilde \eps}$ indicate higher dependence between orthogonalized versions of $w_j$ and $\eps$.  Bounds on the size of $A_1$ are given first. 

Let $A_{1k}$ be the set of $j\in A_1$ such that $j$ is selected prior to the $k$-th true selection, but not prior to any earlier true selections.
For $j \in A_1, k,l \in S_0$, let $C_1>0$ and $1 \geq C_2 > 0$ be constants which satisfy
$$\tilde \gamma_j ' \tilde \theta \geq \tilde \theta_k  C_1\text{ for  }j \in A_{1k},\text{ and }\tilde \theta_k \geq \tilde \theta_l  C_2 \text{ for }l >k.  $$

%Define the constant $$C_A = \sqrt{2} \frac{ \|  \En[ x_i \eps_i ]\|_{\infty}}{\sqrt{t}} \varphi_{\min}(m+s_0)(G)^{-1}.$$  
 %Let  $$G_\H = \begin{cases} G \text{ in the case of Theorem 1 } \\  \Ep G \text{ in the case of Theorem 4 } \end{cases} \ \ \text{and} \ \ \ \ t_\H = \begin{cases} t \text{ in the case of Theorem 1 } \\  c_{\text{test}}' \text{ in the case of Theorem 4. } \end{cases}$$

The two key constants $C_1$,$C_2$ encode information about \textit{relative} incremental loss values at various points of the forward selection procedure.    Lemma 5 calculates suitable $C_1,C_2$.

\begin{lemma} $\tilde \gamma_j ' \tilde \theta \geq \tilde \theta_k  \( \frac{1}{6}\varphi_{\min}(m+s_0)(G_\H)\)^{1/2}$ for $j \in A_{1k}$.  In addition, in the case of Theorem 1,  $\tilde \theta_k \geq \tilde \theta_l  \varphi_{\min}(m+s_0)(G_\H)^{1/2}$ for $l>k$.  In the case of Theorem 4, $\tilde \theta_k \geq \tilde \theta_l   \(c_{\text{\em test}}'' \varphi_{\min}(m+s_0)(G_\H)\)^{1/2}$ for $l>k$.   \end{lemma}

%\begin{lemma} For $j \in A_{1k}$, $\tilde \gamma_j ' \tilde \theta \geq \tilde \theta_k  \( \frac{1}{6}\varphi_{\min}(m+s_0)(G_\H)\)^{1/2}$, provided that $m + s_0 +1 <n$ in the case of Theorem 1.  In addition, for $l > k$, $\tilde \theta_k \geq \tilde \theta_l  \varphi_{\min}(m+s_0)(G_\H)^{1/2}$  in the case of Theorem 1 and $\tilde \theta_k \geq \tilde \theta_l   \(c_{\text{\em test}}'' \varphi_{\min}(m+s_0)(G_\H)\)^{1/2}$ in the case of Theorem 4.   \end{lemma}

Define the following two $s_0 \times s_0$ matrices.
%$$\Gamma = 
%\[\begin{array}{ccccccc} 
%\sum \limits_{j\in m_1} \tilde  \gamma_{j1} & \sum \limits_{j \in m_1} \tilde \gamma_{j2} & ... & \sum \limits_{j \in m_1} \tilde \gamma_{j\hat k} &...&  \sum \limits_{j \in m_1} \tilde \gamma_{j(s_0+1)}  \\ 
%\\ 0 & \sum \limits_{j\in m_2} \tilde \gamma_{j2} & ... & \sum \limits_{j \in m_2} \tilde \gamma_{j\hat k}&...&  \sum \limits_{j \in m_2} \tilde \gamma_{j(s_0+1)}  \\ 
%\\ \vdots & \vdots &  & \vdots && \vdots \\ 
%\\ 0 & 0 & ... &  \sum \limits_{j \in m_{\hat k+1}} \tilde \gamma_{j(\hat k+1)} &  ...&\sum \limits_{j \in m_{\hat k+1}} \tilde \gamma_{j(s_0+1)} \\
%\\ 0 & 0 & ... & 0 & ...& 0\\
%\vdots & \vdots &  & \vdots && \vdots \\ 
%\\ 0 & 0 & ... & 0 & ...& 0\\
%\end{array} \] $$ 
$$\Gamma: \  \Gamma_{kl} = \sum_{j \in A_{1k}} \tilde \gamma_{jl}, \ \ B \ : \  \text{symmetric} , \  B_{kl} = \tilde \theta_{l}/\tilde \theta_{k} \ \text{if} \ l\geq k.$$ 

%$$\Gamma = \[\begin{array}{cccc} \sum \limits_{j\in A_{11}} \tilde  \gamma_{j1} & \sum \limits_{j \in A_{11}} \tilde \gamma_{j2} & ... & \sum \limits_{j \in A_{11}} \tilde \gamma_{js_0}\\  \\ 0 & \sum \limits_{j\in A_{12}} \tilde \gamma_{j2} & ... & \sum \limits_{j \in A_{12}} \tilde \gamma_{j{s_0}} \\ \\ \vdots & \vdots & \ddots & \vdots \\ \\ 0 & 0 & ... &  \sum \limits_{j \in A_{1{s_0}}} \tilde \gamma_{j{s_0}} \end{array} \],   \ B =  \[\begin{array}{cccc} \frac{\tilde\theta_{1}}{\tilde\theta_{1}} & \frac{\tilde\theta_{2}}{\tilde\theta_{1}} &...& \frac{\tilde\theta_{{s_0}}}{\tilde\theta_{1}} \\  \\\frac{\tilde\theta_{2}}{\tilde\theta_{1}} & \frac{\tilde\theta_{2}}{\tilde\theta_{2}} & ... & \frac{\tilde\theta_{{s_0}}}{\tilde \theta_{2}} \\ \\ \vdots & \vdots & \ddots & \vdots \\ \\ \frac{\tilde\theta_{{s_0}}}{\tilde\theta_{1}} & \frac{\tilde\theta_{{s_0}}}{\tilde\theta_{2}} &...& \frac{\tilde\theta_{{s_0}}}{\tilde\theta_{{s_0}}}   \end{array} \] $$

\noindent Empty sums are taken to be 0.   In case $\tilde \theta_k = 0$ for some $k$, $\tilde \theta_l/\tilde \theta_k$ is defined to be $ 1$.  The above definitions of $\Gamma, B$ are useful because the diagonal elements of the product $\Gamma B$ satisfy the equality 
$[\Gamma B]_{kk} = \sum_{j \in A_{\textcolor{black}{1}k}} \tilde \gamma_j ' \tilde \theta/\tilde \theta_{k}.$  This follows from $\Gamma$ being upper triangular (by the orthogonalization construction) and from the fact that $A_{1k}$ is empty if $\tilde \theta_k=0$ (see remark in the proof of Lemma 5).

The definition of $C_1$ implies that
$[\Gamma B]_{kk} \geq C_1|A_{1k}|  \ \text{and subsequently} \  $ $$ \text{tr}(\Gamma B) \geq C_1m_1 .$$

The product $\Gamma B$ has a convenient decomposition.  Consider the set $$\mathscr G_{s_0} = \left \{ Z \in \mathbb R^{s_0 \times s_0} : \begin{array}{c} Z_{kl} = \langle \mathsf X_k, \mathsf Y_l \rangle_{\mathsf H_1} \text{ for some elements } \ \| \mathsf X_k \|_{\mathsf H_1} ,\|\mathsf Y_l \|_{\mathsf H_1} \leq1 \\  \text{in some $s_0$-dimensional real Hilbert space } \mathsf H_1 \end{array} \right \}$$

\begin{lemma} The matrix product $\Gamma B$ may be expressed as $\Gamma B = \Gamma C_3 \bar Z$ where $C_3$ is a constant which may be taken as $C_3 = C_2^{-2}$ and where $\bar Z'  \in \mathscr G_{s_0}$. \end{lemma} 
\textcolor{black}{This decomposition is helpful because of the following result due to Grothendieck. }  \begin{lemma}[\cite{Grothendieck:Resume}] $\sup_{ Z \in \mathscr G_{s_0}} \text{\em tr}(MZ) \leq K^{\mathbb R}_G \| M \|_{\infty \rightarrow 1}$.\end{lemma} 
\textcolor{black}{Here, $ K_G^{\mathbb R} $ is an absolute constant which is known to be less than 1.783.  Importantly, it does not depend on $s_0$.   The notation $\| \cdot \|_{\infty \rightarrow 1}$ indicates the operator norm for bounded linear operators $L^\infty \rightarrow L^1$.  When the matrix $M$ is $s_0 \times s_0$ dimensional, the implied $L^\infty, L^1$ spaces are $L^\infty (\{1,...,s_0\}), L^1 (\{1,...,s_0\})$ or equivalently $(\mathbb R^{s_0}, \| \cdot \|_{\infty})$, $(\mathbb R^{s_0}, \| \cdot \|_{1})$.   The form used here is that described in \cite{Vershynin:Le:GraphConcentration}, Equation 3.2}. 
Therefore,
  $$C_3^{-1} C_1m_1 \leq C_3^{-1} \text{tr}(\Gamma B) =  C_3^{-1} \text{tr} ( \Gamma C_3 \bar Z) =    \text{tr} (  \Gamma'\bar Z' ) \leq K^{\mathbb R}_G \| \Gamma'\|_{\infty \rightarrow 1}.$$%$ Grothendieck's inquality gives
%$ C_3^{-1}\text{tr}(\Gamma B)\leq K_G^{\mathbb R}  \| \Gamma' \|_{\infty \rightarrow 1}.  $
%\textcolor{black}{  Next, let $\mathsf{Sym}\bar B$ be a symmetrized version of $\bar B$ according to $\mathsf{Sym}\bar B_{kl} = \bar B_{lk}$ for $l>k$.  Because $\Gamma$ is upper triangular and $\bar B$ is symmetric,}
%\textcolor{black}{  Because $\Gamma$ is upper triangular,}
%Assume without loss of generality that $C_2 \leq 1$.

In light of this lower bound on $\| \Gamma '\|_{\infty \rightarrow 1}$, there exists $\nu \in \{-1,1\}^{s_0}$ such that $$\| \textcolor{black}{\Gamma' \nu}\|_1 \geq \(K_G^{\mathbb R} \)^{-1}{C_3}^{-1} C_1m_1.$$ 

On the other hand, $\| \Gamma' \nu \|_1$ may be upper bounded by a quantity that depends on $s_0$ by $$\| \Gamma' \nu \|_1 \leq s_0^{1/2} \| \Gamma' \nu \|_2.$$

A key property of the $\tilde \gamma_{j}$, which constitute $\Gamma$, is that they are approximately orthogonal to each other in the sense of the following lemma.  In particular, signed sums of $\tilde \gamma_j$ scale in norm like $m_1^{1/2}$ up to a factor depending on $\varphi_{\min}(m+s_0)(G_\H)$.

\begin{lemma} For any signs $e_j \in \{ -1, 1 \}$,
$\big \| \sum_{j \in A_1} e_j \tilde \gamma_j  \big \|_2 \leq  {m_1}^{1/2} \varphi_{\min}({m}+s_0)(G_{\mathsf H})^{-1/2}.$ \end{lemma}

%For this particular choice of $\nu$, it follows that $ \| \textcolor{black}{\Gamma ' \nu} \|_2 \geq s_0^{-1/2} \(K_G^{\mathbb R} \)^{-1}{C_3}^{-1} C_1 m_1.$
%\noindent 

Therefore, $ \| \textcolor{black}{\Gamma' \nu} \|_2 =   \| \sum_{k=1}^{s_0} \sum_{j \in A_{1k}} \nu_k \tilde \gamma_{j}  \|_2 \leq   m_1^{1/2} \varphi_{\min}({m}+s_0)(G_{\mathsf H})^{-1/2},$ 
%It can be shown for any choice of signs $e_j \in \{-1,1\}^{m_1}$ that $ \| \sum_{j=1}^{m_1}  e_j \tilde \gamma_{j}  \|^2_2 \leq  m_1 \varphi_{\min}(m+s_0)(G)^{-1}$ in the case of Theorem 1 and  $ \| \sum_{j=1}^{m_1}  e_j \tilde \gamma_{j}  \|^2_2 \leq  m_1 \varphi_{\min}(m+s_0)(\Ep G)^{-1}$ in the case of Theorem 4.  
which, when combined with the bound $\| \Gamma' \nu \|_1 \leq s_0^{1/2} \| \Gamma' \nu \|_2,$ immediately implies the following. \begin{lemma}  $m_1 \leq  \varphi_{\min}(m+s_0)(G_\H)^{-1} C_1^{-2} {C_3}^2 \(K_G^{\mathbb R} \)^2s_0.$\end{lemma}

Having controlled $m_1$, it is left to give a bound which controls $m_2$.  The following lemma is proven by showing that the orthogonalization process $w_j \mapsto \tilde w_j$ cannot create too many variables, $j$, with large $\tilde \gamma_{j\tilde \eps}$, given the relevant regularization condition is met.  

\begin{lemma} 
$m_2 \leq  \textcolor{black}3 (m_1 + s_0)$ provided ${t^{1/2}} \geq  2\varphi_{\min}(m+s_0)(G)^{-1} \| \En[x_i \eps_i]\|_{\infty}$ in the case of Theorem 1 and $\Ep[\En[\eps_i^{\text{a\em 2}}]] \leq \frac{1}2 \varphi_{\min}(m+s_0)(\Ep[G])^{-1} c_{\text{\em test}}'$ in the case of Theorem 4. 
\end{lemma}

The final lemma restates the sparsity bounds of Theorems 1 and 4. Its proof involves only assembling the previous arguments.  Recall that $m_1+m_2 = m$ and that $m$ is the number of false selections being considered.

\begin{lemma}
 In the case of Theorem 1, if ${t^{1/2}} \geq  2\varphi_{\min}(m+s_0)(G)^{-1} \| \En[x_i \eps_i]\|_{\infty}$ holds, then also $m \leq 80 \times \varphi_{\min}(m+s_0)(G)^{-4} s_0 $  holds.   In the case of Theorem 4, $\hat s \leq (80 \times \varphi_{\min}(m+s_0)(\Ep [G])^{-4} c_{\text{\em test}}''^{-3}+1)s_0$.
\end{lemma}

This completes the proof of Theorems 1 and 4.

\bibliographystyle{abbrvnat}
\bibliography{dkbib1}

\pagebreak

\setcounter{page}{1}
\section*{S.  Supplement to ``Analysis of Testing-Based Forward Model Selection''}

This supplement proves Theorems 2 and 3, supporting lemmas for Theorems 1 and 4, and Theorem 5.
\setcounter{section}{19}
\setcounter{subsection}{0}
\subsection{Proof of Theorems 2 and 3}
Theorem 2 follows by applying Theorem 1 in the following way.  If $\hat s$ grows faster than $s_0$, then there is $m< \hat s$ such that $s_0 < m < K_n$ and $m/s_0$ exceeds $c_{\text{F}}'(K_n) = O(1)$, giving a contradiction.  The first statement of the theorem follows from applying the bound on $\hat s$.  Theorem 3 follows by $\| \theta_0 - \thetahat \|_1 \leq \sqrt{\hat s + s_0} \| \theta_0 - \thetahat \|_2 \leq \sqrt{\hat s + s_0} \varphi_{\min}(\hat s + s_0)(G)^{-1} \En[(x_i'\theta_0 - x_i' \thetahat)^2]^{1/2}$.   

\subsection{Proof of Lemmas 3 and 4}

%\subsubsection{Proof of Lemma 1}
%Lemma 1 holds by $\ell(\hat S \cup S_0) \leq \ell(S_0) \leq \ell(\theta_0)$.  

%\subsubsection{Proof of Lemma 2}
%Lemma 2 is Lemma 3.3 of \cite{Submodular:Spectral}.   

\subsubsection{Proof of Lemma 3}

It was already shown that $\ell(\hat \theta) \leq \ell(\theta_0) + s_0t\varphi_{\min}(\hat s + s_0)(G)^{-1}.$ 
 Expanding the above two quadratics in $\ell(\cdot)$ gives
\begin{align*} \En [ (x_i'\theta_0-x_i'\hat \theta )^2] &\leq |2\En[\eps_i x_i'(\hat \theta - \theta_0)] |+ s_0t\varphi_{\min}(\hat s + s_0)(G)^{-1}\\
& \leq 2\| \En[\eps_i x_i] \|_\infty \|\theta_0 - \hat \theta \|_1 + s_0t\varphi_{\min}(\hat s+s_0)(G)^{-1}.
\end{align*}
%If $\En [ (x_i'\theta_0-x_i'\hat \theta )^2] \geq s_0t\varphi_{\min}(s_0)(G)^{-1}$, then
To bound $\| \theta_0 - \hat \theta \|_1$: 
\begin{eqnarray*}
&\| \theta_0 - \thetahat \|_1 &\leq \sqrt{\hat s + s_0} \| \theta_0 - \thetahat \|_2\\ 
&&\leq \sqrt{\hat s + s_0} \varphi_{\min}(\hat s + s_0)(G)^{\textcolor{black}{-1/2}} \En[(x_i'\theta_0 - x_i' \thetahat)^2]^{1/2} .  
\end{eqnarray*}

 \textcolor{black}{If  $\En[(x_i'\theta_0 - x_i' \thetahat)^2]^{1/2} =0$, then the first conclusion of Theorem 1 holds.  Otherwise,} combining the above bounds and dividing by $\En[(x_i'\theta_0 - x_i' \thetahat)^2]^{1/2} $ gives 
\begin{align*}\En [ (x_i'\theta-x_i'\hat \theta )^2]^{1/2} &\leq    2\| \En[\eps_i x_i] \|_\infty \sqrt{\hat s + s_0} \varphi_{\min}(\hat s + s_0)(G)^{\textcolor{black}{-1/2}} \\ & + \frac{s_0t\varphi_{\min}(\hat s + s_0)(G)^{-1}}{\En[(x_i'\theta_0 - x_i' \thetahat)^2]^{1/2}  }.
\end{align*}
Finally, either $\En[(x_i'\theta_0 - x_i' \thetahat)^2]^{1/2} \leq \sqrt{  s_0t\varphi_{\min}(\hat s + s_0)(G)^{\textcolor{black}{-1/2}} }$, in which case Lemma 3 holds,  or alternatively $\En[(x_i'\theta_0 - x_i' \thetahat)^2]^{1/2} > \sqrt{  s_0t\varphi_{\min}(\hat s + s_0)(G)^{\textcolor{black}{-1/2}} }$, in which case 
\begin{align*}\En [ (x_i'\theta-x_i'\hat \theta )^2]^{1/2} &\leq    2\| \En[\eps_i x_i] \|_\infty \sqrt{\hat s + s_0} \varphi_{\min}(\hat s + s_0)(G)^{\textcolor{black}{-1/2}} \\ & + \sqrt{s_0t\varphi_{\min}(\hat s + s_0)(G)^{-1}}.  
\end{align*}

\subsubsection{Proof of Lemma 4}

For any $S$ define $\theta_S^*$ to be the minimizer of $\mathscr E(S)$.  For any $S$ define also $d_{S} =\theta_{ S}^* - \theta_{S_0\cup S}^*$.  Finally, let $\delta_{0,S} = \theta_0 - \theta_{S_0\cup S}^*$.  Note that
$\E( S)  - \E(S_0 \cup S)= d_{ S} '\Ep[G]d_{ S}.$
By arguments in the earlier sections, 
$d_{\hat S} '\Ep[G ]d_{\hat S} \leq s_0  c_{\text{test}} \varphi_{\min}(K_{\text{test}})(\Ep[G])^{-1}.$
But $d_{\hat S} '\Ep[ G]d_{\hat S} \geq \varphi_{\min}(K_{\text{test}})(\Ep[G])\| d_{\hat S}  \|_2^2.$
So $ \|d_{\hat S} \|_2 \leq \sqrt{s_0c_{\text{test}}} \varphi_{\min}(K_{\text{test}})(\Ep[G])^{-1}.$ In addition, $\delta_{0,S}$ is bounded by $$\| \delta_{0,S}\|_2  = \| \Ep[\En[x_{iS_0\cup S}'\eps_i] ]\|_2$$ $$ \leq (|S| + s_0)^{1/2} \max_{j} |\Ep [\En [ x_{ij} \eps_i^{\text{a}}] ]| \leq \frac{1}{2} \sqrt{(|S| +  s_0) c_{\text{test}}} \varphi_{\min}(K_{\text{test}})(\Ep[G])^{-1}$$  
where the last bound comes from Cauchy-Schwarz
 (passing to $\Ep [\En[ x_{ij}^2]]^{1/2} \Ep[\En[ \eps_{i}^{\text{a}2}]]^{1/2} $) 
 along with the assumed condition on $\eps_{i}^{\text{a}}$ and the fact that $c_{\text{test}}' \leq c_{\text{test}}$.  Next,
\begin{align*}\hat \theta &= G_{\hat S}^{-1} \En [x_{i\hat S} (x_{i\hat S}'\theta_{\hat S}^* + \varepsilon_i - x_{i\hat S \cup S_0}' d_{\hat S} + x_{i \hat S\cup S_0}'\delta_{0,\hat S} )] \\
& = \theta_{ \hat S}^* +  G_{\hat S}^{-1} \En  \left [x_{i\hat S}\varepsilon_i \right ]+G_{\hat S}^{-1} \En  \left [x_{i\hat S} x_{i \hat S\cup S_0}'(-d_{\hat S}+ \delta_{0,\hat S} ) \right ] \\
 \Rightarrow \| \thetahat - \theta_{\hat S}^* \|_2 &\leq \varphi_{\min}(\hat s)( G)^{-1/2}  \left  \|  \En [x_{i\hat S} \varepsilon_{i} ]\right \|_2+\left \|G_{\hat S}^{-1}    \En \left [x_{i \hat S} x_{i \hat S \cup S_0} (-d_{\hat S}+ \delta_{0,\hat S} )  \right ]\right \|_2 \\
& \leq \varphi_{\min}(\hat s)( G)^{-1/2}\hat s^{1/2} \|  \En [x_{i} \varepsilon_{i}  ]\|_{\infty}  \\ & \ \ \ +\varphi_{\min}(\hat s)(G)^{-1/2} \varphi_{\max}(\hat s + s_0)(G)^{1/2} (\| d_{\hat S} \|_2 + \| \delta_{0,\hat S} \|_2).
\end{align*}
Finally,
\begin{align*} ( \En [(x_i' \thetahat - x_i' \theta_0)^2] )^{1/2}  & \leq \varphi_{\max}(s_0 + \hat s \hspace{.5mm}  )(G)^{1/2}\|  \thetahat - \theta_0\|_2\\
& \leq   \varphi_{\max}(s_0 + \hat s \hspace{.5mm}  )(G)^{1/2}( \|  \thetahat - \theta_{\hat S}^*\|_2 + \| \delta_0 \|_2 + \| d_{\hat S} \|_2 )\ \\
& \leq    \varphi_{\max}(s_0 + \hat s)(G)^{1/2}   \varphi_{\min}(s_0 + \hat s)(G)^{-1/2}\hat s^{1/2}  \hspace{.5mm} \| \En[x_i\eps_i] \|_\infty \\ & \ \ \ +  \varphi_{\max}(s_0 + \hat s \hspace{.5mm}  )(G)^{1/2} (\frac{3}{2} + \frac{3}{2}  \varphi_{\max}(s_0 + \hat s \hspace{.5mm}  )(G)^{1/2} \varphi_{\min}(\hat s + s_0)(G)^{-1/2} ) \\ & \ \ \ \ \ \ \times \sqrt{(\hat s + s_0) c_{\text{test}}}\varphi_{\min}(K_{\text{test}})(\Ep[G])^{-1} \\
& \leq    \varphi_{\max}(s_0 + \hat s)(G)^{1/2}   \varphi_{\min}(s_0 + \hat s)(G)^{-1/2}\hat s^{1/2}  \hspace{.5mm} \| \En[x_i\eps_i] \|_\infty \\ & \ \ \ +  3\varphi_{\max}(s_0 + \hat s \hspace{.5mm}  )(G)  \varphi_{\min}(\hat s + s_0)(G)^{-1/2}  \sqrt{(\hat s + s_0) c_{\text{test}}}\varphi_{\min}(K_{\text{test}})(\Ep[G])^{-1} . 
\end{align*}  \vspace{-1cm}

\subsection{Proof of Supporting Lemmas for Sparsity Bounds for Theorems 1 and 4.}

\subsubsection{Additional notation}

Additional notation is used for the proof of the lemmas which follow.  
The inner product from $\H$ is hereafter denoted simply with $\langle\hspace{.5mm}  \cdot \hspace{.5mm} , \hspace{.5mm} \cdot  \hspace{.5mm} \rangle_\H= \langle\hspace{.5mm}  \cdot \hspace{.5mm} , \hspace{.5mm} \cdot  \hspace{.5mm} \rangle$.   The symbol  $'$  is kept for use for transposition of finite dimensional real matrices and vectors derived from certain elements of $\mathsf H$ defined below.  For $a, b \in L^2(\Omega : \mathbb R^n)$, $a'b$ is defined pointwise (over $\Omega$) and thus defines a random variable $\Omega \rightarrow \mathbb R$ and $\langle a, b \rangle = \Ep [a'b]$.   In the case of Theorem 1, $a'b = \langle a,b\rangle$. 

Let $V = [v_1,...,v_{s_0}]$ with the understanding that $V$ and similar quantities are formally defined as linear mappings $\mathbb R^{s_0} \rightarrow \mathsf H$.   Then $y = V \theta_0 + \eps$ is well defined for both Theorems 1 and 4.

  Let  $\mathscr M_{k}$ denote projection in $\mathsf H$ onto the space orthogonal to $\text{span}(\{v_1,...,v_k\})$. Then note that $\tilde v_k = \frac{ \mathscr M_{k-1}v_k }{ \langle v_k,\mathscr M_{k-1}v_k\rangle^{1/2}} \ \text{for} \ k = 1,...,s_0.$  
In addition,  $\tilde \eps =  \frac{\mathscr M_{s_0 }\eps}{(\langle \eps,\mathscr M_{s_0}\eps \rangle)^{1/2}}.$ 
%Let $\tilde V_{\text{temp}} = [\tilde v_1,...,\tilde v_{s_0}]$, ordered according to the temporary order.  \textcolor{black}{Note that the process obtaining the new elements $\tilde v_k ,\tilde \eps$ is the Gram-Schmidt process.  As a result, it also possible to define $\tilde V_{\text{temp}}$ through a thin QR decomposition $V = \tilde V_{\text{temp}} \mathsf  R_{\text{temp}}$.  } 
%Let $\hat k$ denote the index of the final true covariate selected into the model when the $m$-th false covariate is selected.  The variables $\tilde v_1,...,\tilde v_{\hat k}$ maintain their original order.  The unselected true covariates $\tilde v_{\hat k+1},..., \tilde v_{s_0}$ are reordered in such a way that under the new ordering, $\textcolor{black}{ \tilde \theta_{k} \geq \tilde \theta_{l}}$ whenever $l > k$ \textcolor{black}{for $\tilde \theta$ defined by  $\tilde V \tilde \theta + \tilde \theta_{\tilde \eps} \tilde \eps = y$ and $ \tilde V = [\tilde v_1,...,\tilde v_{s_0}]$ consistent with the new ordering.}   No new orthogonalization is done.  
%For notational convenience, indexing of true covariates will be done with the letter $k$ (so that the true covariates are $w_k$ with $k$ ranges from 1,...,$s_0$.)  Indexing of false covariates will be done with letter $j$ (so that the false selected covariates are $w_j$ with $j$ in $1,...,m$.)  The original covariates $x_j$ are also indexed by $j$ this time in the range $j=1,...,p$.
For more general sets $S$, let $\mathscr Q_{S}$ be projection onto the space orthogonal to $\text{span}(   \{   x_j, \ j \in S \}  )$.  For each selected covariate, $w_j$, set  $S_{\text{pre-}w_j}$ to be the set of (both true and false) covariates selected prior to $w_j$. 

\subsubsection{Proof of Lemma 5} It is needed to calculate $C_1, C_2$ such that  $\tilde \gamma_j ' \tilde \theta \geq \tilde \theta_k C_1$ for $j \in A_{1k}$ and $ \tilde \theta_k \geq  \tilde \theta_l C_2$ for $l >k$.  
% It was noted in the earlier text that taking $\tilde \theta_k \neq 0$ is without loss of generality.  This is addressed at the end of this section.  
Define $$\Delta_j\ell^{\hspace{.5mm}\H}(S)=  \begin{cases} \Delta_j\ell(S) \text{ in the case of Theorem 1} \\ \Delta_j\mathscr E(S) \text{ in the case of Theorem 4} \end{cases} $$

Also recall that $t_\H = t$ in the case of Theorem 1 and $t^\H = c_{\text{test}}'$ in the case of Theorem 4.  Note that $c_{\text{test}''}$ is not defined in the context of Theorem 1.  In the case of Theorem 1, during the proof of this lemma, $c_{\text{test}}''$ is taken to be equal to 1.

A simple derivation can be made to show that

$$-\Delta_j \ell^{\hspace{.5mm}\H}(S_{\text{pre-}w_j}) =  \frac{1}{n} \langle y, \tilde w_j \rangle  (\langle \tilde w_j, \tilde w_j\rangle )^{-1}\langle  \tilde w_j, y \rangle  = \frac{1}{n} \frac{1}{\| \tilde w_j \|_{\mathsf H}^2} (\tilde \theta ' \tilde \gamma_j + \tilde \theta_{\tilde \eps}\tilde \gamma_{j\tilde \eps})^2.$$
Note the slight abuse of notation in $-\Delta_j\ell^{\hspace{.5mm}\H}(S_{\text{pre-}w_j})$ signifying change in loss under inclusion of $w_j$ rather than $x_j$. Next, 
$$ (\tilde \theta ' \tilde \gamma_j + \tilde \theta_{\tilde \eps}\tilde \gamma_{j\tilde \eps})^2 \leq 2(\tilde \theta ' \tilde \gamma_j )^2 + 2 ( \tilde \theta_{\tilde \eps}\tilde \gamma_{j\tilde \eps})^2.$$
Since $\tilde \theta_{\tilde \eps} = $\textcolor{black}{ $\langle \tilde \eps,y\rangle $ $= \langle \eps , \mathscr M_{s_0} y \rangle / \langle \eps, \mathscr M_{s_0} \eps\rangle^{1/2} $} $=   \langle \eps, \mathscr M_{s_0} \eps\rangle^{1/2},$  $\| \tilde w_j \|_{\mathsf H}^2 \geq 1$, and $j \in A_1$ it follows that $$\frac{1}{n}\frac{1}{\| \tilde w_j \|_{\mathsf H}^2} ( \tilde \theta_{\tilde \eps}\tilde \gamma_{j\tilde \eps})^2 \leq \frac{1}{n}  \frac{1}{\| \tilde w_j \|_{\mathsf H}^2}  \tilde \theta_{\tilde \eps}^2\( \frac{t_\H^{1/2}n^{1/2}}{ (\textcolor{black}{ 3 }\langle \eps,\mathscr M_{s_0} \eps \rangle )^{1/2}}  \)^2 \leq \frac{t_\H}{\textcolor{black}{3}}.$$
This implies
$$ \frac{1}{2}( -\Delta_j \ell^{\hspace{.5mm} \H}(S_{\text{pre-}w_j})  )\leq  \frac{1}{n} \frac{1}{\| \tilde w_j \|_{\mathsf H}^2} (\tilde \theta'  \tilde \gamma_j)^2 + \frac{t_\H}{\textcolor{black}{3}}.$$
By the condition that the false $j$ is selected, it holds that $ -\Delta_j \ell^{\hspace{.5mm} \H} (S_{\text{pre-}w_j})  > t_\H$ and so $ \frac{1}{\textcolor{black}{3}}(-\Delta_j \ell^{\hspace{.5mm} \H}(S_{\text{pre-}w_j}) )> \frac{t_\H}{\textcolor{black}{ 3}}$ which implies that \textcolor{black}{$ -\frac{t_\H}{3} > \frac{1}{3} \Delta_j \ell^{\hspace{.5mm} \H}(S_{\text{pre-}w_j})$ and}
%$$\frac{\tilde \gamma_j ' \tilde \theta }{1 + \tilde \gamma_j' \tilde \gamma_j } > \frac{(\tilde v_j '\mathscr Q_j \tilde V \tilde \theta)^2}{\tilde v_j '\mathscr Q_j \tilde v_j}> \frac{(\tilde v_j '\mathscr Q_j \tilde V \tilde \theta)^2}{\tilde v_j '\tilde v_j}.$$
$$ \frac{1}{2}(-\Delta_j \ell^\H(S_{\text{pre-}w_j})) - \frac{t_\H}{\textcolor{black}{ 3}}  \geq \frac{1}{\textcolor{black}{ 6}}(-\Delta_j \ell^\H(S_{\text{pre-}w_j})).$$
Finally, this yields that
$$\frac{1}{n\| \tilde w_j \|_{\mathsf H}^2} (\tilde \gamma_j ' \tilde \theta)^2  \geq  \frac{1}{\textcolor{black}{ 6}}(- \Delta_j \ell^\H(S_{\text{pre-}w_j})).$$
By the fact that $w_j$ was selected ahead of $v_k$ it holds that 
$$- \Delta_j \ell^{\hspace{.5mm} \H}(S_{\text{pre-}w_j}) \geq - \Delta_k \ell^\H(S_{\text{pre-}w_j})c_{\text{test}}''.$$

%Then 
%Therefore, further bound the righthand side.  Let $\tilde z_k$ be the projection of $\tilde v_k$ onto the space orthogonal to all previously selected (true and false) covariates.  Then \textcolor{black}{$\tilde z_k = \mathscr Q_{\text{pre-}j}\tilde v_k$} and
%$$ - \Delta_k \ell(S_{\text{pre-}w_j}) \geq \frac{1}{n} \tilde z_k'\tilde z_k \tilde \theta_k^2.$$
%{$$ - \Delta_k \ell(S_{\text{pre-}w_j}) = - \Delta_k \ell(S_{\text{pre-}w_j} \cup  \{ \tilde v_1,...\tilde v_{k-1} \} ) = - \Delta_{\tilde v_k}\ell(S_{\text{pre-}w_j}  \cup \{ \tilde v_1,...\tilde v_{k-1} \})  .$$}
\textcolor{black}{
Next,} to lower bound $-\Delta_k \ell^{\hspace{.5mm}\H}(S_{\text{pre-}w_j})$, define a perturbed version of $\ell^{\hspace{.5mm}\H}$.  Let $\xi \in \mathbb \H$.  Let $\ell_{y + \xi}^\H$ be defined analogously to $\ell^{\hspace{.5mm}\H}$ except with the role of $y$ in $\ell$ played by $y+\xi$ in $\ell^{\hspace{.5mm}\H}_{y+ \xi}$.   Choose $\xi$ such that $\langle \xi, w_j  \rangle= 0$ for $j=1,...,m$, $\langle \xi, v_k \rangle = 0$ for $v_k = 1,...,s_0$ and $\langle \xi , \eps \rangle= 0$.  In the case of Theorem 1, $\xi\neq 0$ exists provided $m+s_0+1 <n$.  If not, then $\H$ can be enlarged appropriately to allow $\xi$ to exist, for example, with the inclusion $\iota:H \rightarrow H \oplus \mathbb R$, $x\mapsto (x,0)$, $\xi = (0,1)$.  
%, in which case $\ell^{\H \oplus \mathbb R}$ agrees with $\ell^{\H}$ all of the arguments that follow continue to hold and provide the needed bounds.  
Then due to the orthogonality of $\xi$ to $w_j$ and $v_k$ and $\eps$, it follows that 
$$ - \Delta_k \ell^\H(S_{\text{pre-}j}) = - \Delta_k \ell^\H_{y + \xi} (S_{\text{pre-}j}) $$
with the right hand side possibly defined on an enlarged $\H$ as described above.  

Next, the following reduction holds.  
$$ - \Delta_k \ell^\H_{y + \xi} (S_{\text{pre-}w_j})  \geq  - \Delta_{k} \ell^\H_{y + \xi} (S_{\text{pre-}w_j } \cup \{  \tilde v_{k+1} \tilde \theta_{k+1} + ... + \tilde v_{s_0}\tilde \theta_{s_0} + \tilde \eps +  \xi \}) .$$
$$  =  - \Delta_{\tilde v_k}  \ell^\H_{y + \xi} (S_{\text{pre-}w_j } \cup \{  \tilde v_{k+1} \tilde \theta_{k+1} + ... + \tilde v_{s_0}\tilde \theta_{s_0} + \tilde \eps +  \xi \}) .$$
 \textcolor{black}{Let $\overset{  \leftrightarrow   _\xi}{\mathscr M_{k} }$ be projection on the corresponding orthogonal space to the span of the vectors listed in  $S_{\text{pre-}w_j }\cup \{  \tilde v_{k+1} \tilde \theta_{k+1} + ... + \tilde v_{s_0}\tilde \theta_{s_0} + \tilde \eps +  \xi \}$.}  \textcolor{black}{(The accent $\overset{\leftrightarrow}\cdot$ meant to emphasize that covariates selected before and after $v_k$ (or not at all) are considered.)  Then the above term is further reduced by}

\textcolor{black}{$$  =   \frac{1}{n} \frac{ \langle (y + \xi) ,  \overset{  \leftrightarrow   _\xi}{\mathscr M_{k} }\tilde v_k \rangle^2}{ \langle \tilde v_k, \overset{  \leftrightarrow   _\xi}{\mathscr M_{k}  } \tilde v_k \rangle}= \frac{1}{n} \frac{\langle \tilde \theta_k \tilde v_k   , \overset{  \leftrightarrow   _\xi}{\mathscr M_{k} } \tilde v_k\rangle^2}{ \langle \tilde v_k , \overset{  \leftrightarrow   _\xi}{\mathscr M_{k} } \tilde v_k \rangle}   = \frac{1}{n} \tilde \theta_k^2   \langle \tilde v_k , \overset{  \leftrightarrow   _\xi}{\mathscr M_{k}  } \tilde v_k \rangle.$$}

%Note that $ \tilde v_k   '  \overset{\leftrightarrow} {\mathscr M}_{k} \tilde v_k = \tilde z_k' \tilde z_k$.  This follows from a general result for projections onto complements described below in this paragraph.  Consider a real Hilbert space $\mathsf V$.  (This level of generality will be helpful in the proof of Theorem 4).  Let $\mathsf A, \mathsf B$ be closed subspaces. Let $\mathsf v \in \mathsf B^\perp$.  Let $\mathsf{AB} = \text{span}({\mathsf A \cup \mathsf B})$.  Let $\mathscr P_{\mathsf A}, \mathscr M_{\mathsf A}$ denote projections onto $\mathsf A$ and $\mathsf A^\perp$ and define $\mathscr P_{\mathsf B}$ etc analogously.  Then $\mathscr M_{\mathsf{AB}}\mathsf v = \mathscr M_{\mathsf A} \mathsf v$

%Let $\mathsf A = \text{span}( \tilde v_1,...,\tilde v_{k-1}, S_{\text{pre-}j})$ and let  $\mathsf B = \text{span}( \tilde v_{k+1},...,\tilde v_{s_0}, \tilde \eps )$.  Let $\mathsf{AB} = \text{span}(\mathsf A \cup \mathsf B)$. Let 

%Assume that $\tilde \theta_k^2 > Ct$.  As otherwise there is nothing to prove.  
%Then $$ \tilde \theta_k \geq C\varphi_{\min}^{-1/2}(G) \| \En[x_i \eps_i \|_\infty$$

\textcolor{black}{
Then seek a lower bound on $\frac{1}{n}\langle  \tilde v_k,  \overset{  \leftrightarrow   _\xi}{\mathscr M_{k} } \tilde v_k \rangle$.  Note that for some vector $\eta_k$ it holds that $\tilde v_k = {\langle v_k , \mathscr M_{k-1} v_k\rangle ^{-1/2}} v_k - [v_1,...,v_{k-1}]\eta_k$.  Then $ \langle \tilde v_k, \overset{  \leftrightarrow   _\xi}{\mathscr M_{k}} \tilde v_k  \rangle =  \langle{v_k , \mathscr M_{k-1} v_k}\rangle^{-1}  \langle v_k,\overset{  \leftrightarrow   _\xi}{\mathscr M_{k} }  v_k\rangle$.  Let $H = [V  \ \ W ]$.  Let $\tilde y_k = \tilde v_{k+1} \tilde \theta_{k+1} + ... + \tilde v_{s_0}\tilde \theta_{s_0} + \tilde \eps$.  
A lower bound on the term $ \langle v_k,\overset{  \leftrightarrow   _\xi}{\mathscr M_{k} }  v_k \rangle$ follows from a lower bound on the eigenvalues of the below matrix for any $c>0$:
$$\langle  v_k, \overset{  \leftrightarrow   _\xi}{\mathscr M_{k} }  v_k \rangle \geq  \lambda_{\min} \( \langle \[ H \ \  ( \tilde y_k +   \xi)c  \], \[ H  \  \ ( \tilde y_k +  \xi )c  \]  \rangle\)$$  
That is, it is enough to bound the spectrum of $n G_{c, \xi}$ defined by $$ G_{c,\xi} = \frac{1}{n}  \[\begin{array}{cc} \langle H, H\rangle &   c\langle \tilde y_k +  \xi,  H \rangle  \\ \langle H,\tilde y_k +  \xi\rangle c    & c^2  \langle \tilde y_k +  \xi, \tilde y_k +  \xi\rangle  \end{array} \].$$
Using the fact that $\xi$ is orthogonal to $H$ and $\eps$, $G_{c,\xi}$ reduces to 
$$G_{c,\xi} =  \frac{1}{n} \[\begin{array}{cc} \langle H, H \rangle &   c\langle \tilde y_k, H \rangle  \\ \langle H, \tilde y_k \rangle c    & c^2 \langle \tilde y_k,\tilde y_k \rangle+ c^2 \langle \xi,\xi \rangle    \end{array} \].$$
As a result of the above reductions, for each $c, \xi$, 
$$  - \Delta_k \ell^{\hspace{.5mm}\H}(S_{\text{pre-}w_j}) \geq \frac{1}{n}\langle v_k, \mathscr M_{k-1} v_k\rangle ^{-1} n   \lambda_{\min}(G_{c, \xi})  \tilde \theta_k^2.$$
And therefore, $$  - \Delta_k \ell^{\hspace{.5mm}\H}(S_{\text{pre-}w_j}) \geq   \frac{1}{n}\langle v_k, \mathscr M_{k-1} v_k\rangle ^{-1} n   \tilde \theta_k^2 \ \  \underset{{  {\Big \{ }  \begin{array}{c} c \rightarrow 0 \\  \frac{1}{n} \langle\xi,\xi\rangle = c^{-2} \end{array} \ {\Big  \}}  }} \lim \ \ \  \lambda_{\min}(G_{c, \xi}). $$
By continuity of eigenvalues for symmetric matrices, passing to the limit gives
$$  - \Delta_k \ell^\H(S_{\text{pre-}w_j}) \geq   \frac{1}{n}\langle v_k, \mathscr M_{k-1} v_k\rangle ^{-1} n   \tilde \theta_k^2 \ \  \lambda_{\min}\( \frac{1}{n} \[\begin{array}{cc} \langle H, H\rangle &   0   \\ 0  &  1   \end{array} \] \) $$
$$ \geq \frac{1}{n}\langle v_k, \mathscr M_{k-1} v_k /n \rangle^{-1}  \tilde \theta_k^2  \varphi_{\min}(m+s_0)(G_{\mathsf H}) \geq  \frac{1}{n} \cdot 1 \cdot \tilde \theta_k^2  \varphi_{\min}(m+s_0)(G_\H). $$
%
%
%These can be bounded by calculating
%
%
%$$ \min_{\alpha, \beta} \  \frac{1}{ \| \alpha \|_2^2 + \beta^2}    \[\begin{array}{cc}  \alpha & \beta \end{array} \]     \[\begin{array}{cc} \En[ x_{Si}x_{Si}'] &  \En [x_{Si} \eps_i] \hspace{.5mm}  \En[\eps_i^2]^{-1/2}  \\ \En [ \eps_ix_{Si}' ]  \hspace{.5mm}  \En[\eps_i^2]^{-1/2}  &1  \end{array} \] \[\begin{array}{c}  \alpha \\  \beta \end{array} \]     $$
%If $\alpha = 0$ and $\beta \neq 0$ then the above expression reduces to $1$.  Assume $\| \alpha \|_2 = 1$.  Then the above expression is bounded by
%$$ \geq\frac{  \varphi_{\min}(m+s_0)(G) - 2 \alpha' \En[ x_{iS} \eps_i] \En[\eps_i^2]^{-1/2} \beta + \beta^2}{ 1 + \beta^2 }$$
%$$ \geq {  \varphi_{\min}(m+s_0)(G) - 2 \alpha' \En[ x_{iS} \eps_i] \En[\eps_i^2]^{-1/2} \beta + \beta^2}$$
%$$ \geq {  \varphi_{\min}(m+s_0)(G) - 2 \| \alpha \|_1 \| \En[ x_{i} \eps_i] \|_\infty \En[\eps_i^2]^{-1/2} \beta + \beta^2}$$
%$$ \geq {  \varphi_{\min}(m+s_0)(G) - 2  \sqrt{m+s_0} \| \En[ x_{i} \eps_i] \|_\infty \En[\eps_i^2]^{-1/2} \beta + \beta^2}$$
This gives $$\frac{1}{n \| \tilde w_j \|_\H^2}(\tilde \gamma_j ' \tilde \theta)^2 \geq c_{\text{test}}'' \frac{1}{6} \frac{1}{n}\varphi_{\min}(m+s_0)(G_{\mathsf H}) \tilde \theta_k^2.$$}
Using the fact that $\| \tilde w_j\|_{\H} \geq 1$ implies that 
$$(\tilde \gamma_{j} ' \tilde \theta )^2 \geq   \tilde \theta_k^2 \hspace{.7mm} \textcolor{black}{c_{\text{test}}'' \frac{1}{6}}  \varphi_{\min}(m+s_0)(G_{\mathsf H}).$$
%\textcolor{black}{If $\tilde \theta_{\tilde \eps} = 0$, it is seen to be sufficient to bound $\lambda_{\min} \(  \[\tilde V \  \ \frac{1}{\sqrt n}  W\]' \[\tilde V \ \ \frac{1}{\sqrt n} W\] \) $ after noting that $  - \Delta_{\tilde v_k} \ell (S_{\text{pre-}w_j }\cup \{ \tilde v_1,...\tilde v_{k-1} \} \cup \{  \tilde v_{k+1},...\tilde v_{s_0} \} \cup \{  \tilde \eps  \}) =   - \Delta_{\tilde v_k} \ell (S_{\text{pre-}w_j }\cup \{ \tilde v_1,...\tilde v_{k-1} \} \cup \{  \tilde v_{k+1},...\tilde v_{s_0} \} ) $ and using similar reductions as above.  In this case, the final conclusion $\frac{1}{n\tilde w_j' \tilde w_j}(\tilde \gamma_j ' \tilde \theta)^2 \geq \frac{1}{4} \frac{1}{n}\varphi_{\min}(m+s_0)(G)\varphi_{\max}(m+s_0)(G)^{-1} \tilde \theta_k^2$ also holds.} 
Now suppose no true variables remain when $j$ is selected.  Then $\langle \tilde w_j, \tilde w_j \rangle= \langle \tilde u_j, \tilde u_j \rangle= 1$ and 
$$-\Delta_{j} \ell^{\hspace{.5mm}\H}(S_{\text{pre-}w_j}) =\frac{1}{n} \tilde \gamma_{j\tilde \eps}^2 \tilde \theta_{\tilde \eps}^2 \geq t_\H.$$
Note that $\tilde \theta_{\textcolor{black}{\tilde \eps}}$ is given by $\tilde \theta_{\tilde \eps} =$ $  \langle \eps, \mathscr M_{s_0} \eps\rangle^{1/2}. $ 
Therefore, $ \tilde \gamma_{j\tilde \eps}^2 \geq t_\H \frac{n}{\langle \eps, \mathscr M_{s_0} \eps\rangle } .$
This implies that $j \in A_2$.  Therefore, set $$C_1^2 = c_{\text{test}}''\textcolor{black}{ \frac{1}{ 6}} \varphi_{\min}(m+s_0)(G_{\mathsf H}).$$

Next, construct $C_2$.  
For each selected true covariate, $v_k$, set  $S_{\text{pre-}v_k}$ to be the set of (both true and false) covariates selected prior to $v_k$.   Note that $$\textcolor{black}{\frac{1}{n}}\tilde \theta_k^2  = - \Delta_k \ell^{\hspace{.5mm} \H}(\{v_1,...,v_{k-1} \}) \geq - \Delta_k \ell^{\hspace{.5mm} \H}(S_{\text{pre-}v_k}) $$  since $\{v_1,...,v_{k-1} \} \subseteq S_{\text{pre}-v_k}$.  In addition, if $v_k$ is selected before $v_l$ (or $v_l$ is not selected), then  
%$$  - \Delta_k \ell^{\hspace{.5mm} \H}(S_{\text{pre-}v_k}) \geq - \Delta_l \ell^{\hspace{.5mm} \H}(S_{\text{pre-}v_k}) \geq \varphi_{\min}(m+s_0)(G_\H) \textcolor{black}{ \frac{1}{n}\tilde \theta_l^2}$$ in the case of Theorem 1 and 
$$  - \Delta_k \ell^{\hspace{.5mm} \H}(S_{\text{pre-}v_k}) \geq  c_{\text{test}}'' ( - \Delta_l \ell^{\hspace{.5mm} \H}(S_{\text{pre-}v_k}) ) \geq c_{\text{test}} '' \varphi_{\min}(m+s_0)(G_\H) \textcolor{black}{ \frac{1}{n}\tilde \theta_l^2}.$$
% in the case of Theorem 4.  
Therefore, taking $$C_2^2 =  c_{\text{test}}'' \varphi_{\min}(m+s_0)(G_\H)$$
%Therefore, taking $$C_2^2 = \begin{cases} \varphi_{\min}(m+s_0)(G_\H) \ \text{in the case of Theorem 1}  \\ c_{\text{test}}'' \varphi_{\min}(m+s_0)(G_\H) \ \text{in the case of Theorem 4} \end{cases}$$
implies that $\tilde \theta_k \geq \tilde \theta_l  C_2$ for any $l > k$.

As a final remark, consider the case that $\tilde \theta_k = 0$.  Then $\tilde \theta_l = 0$ for $l > k$.  Then if $j \in A_{1k}$, it follows that $\tilde \gamma_j'\tilde \theta = 0$.  Therefore, using reasoning as above, $-\Delta_j \ell^\H(S_{\text{pre-}j}) = \frac{1}{n} \frac{1}{\| \tilde w_j\|_\H^2}(\tilde \theta_{\tilde \eps} \tilde \gamma_{j \tilde \eps})^2 \leq \frac{t_\H}{3}$.  But this is impossible, because being selected into the model requires $  -\Delta_j \ell^\H(S_{\text{pre-}j}) > t_\H$.  Therefore, $A_{1k}$ is empty if $\tilde \theta_k = 0$.
\subsubsection{Proof of Lemma 6}

The desired element $\bar Z$ of $\mathscr G_{s_0}$ is constructed as the covariance matrix of certain real, mean-zero,  random vectors 
 $$ \mathsf X = \( \mathsf X_k \)_{k=1}^{s_0} , \ \mathsf Y = \( \mathsf Y_l \)_{l=1}^{s_0}.$$ 
\textcolor{black}{  The random variables $\mathsf X_k, \mathsf Y_l$ constituting $\mathsf X, \mathsf Y$ are defined as follows. Let $\beta_k = \tilde \theta_{k}/\tilde \theta_{k-1}$ for $ k = 2,...,s_0$.  Then note that the components of $B$ can be expressed $B_{kl} = \prod_{q=k+1}^l \beta_{q}$ for $k < l$, and extended symmetrically for components $l < k$.}
 
Decompose the elements of the sequence $\beta_k$ into $$\beta_k = \beta_k^a \beta_k^b$$
in such a way that for all $l\geq k\geq 2$, $$C_2 \leq \prod_{q=k}^{l} \beta_{q}^a \leq C_2^{-1} $$ and for all $k \geq 2$, $$ 0 \leq \beta_k^b \leq 1.$$ Induction establishes the existence of such a decomposition with the additional property that: $\beta_k^a > \beta_k $ only if there is $q \leq k$ such that $\beta_q^a \cdot ... \cdot \beta_{k}^a  = C_2$.    The case $s_0 = 2$ follows by taking $\beta_2^a = \max \{ C_2, \beta_2 \}$ and noting that $\beta_2 = \tilde \theta_2/\tilde \theta_1 \leq C_2^{-1}$.   Assume the complete induction hypothesis that the decomposition exists for sequences with $s_0 =2,...,s$ for some $s$.   Consider a sequence $\beta_2,...,\beta_{s+1}$.  Apply the decomposition to obtain $\beta_k = \beta_k^a \beta_k^b$ for $k \leq s$.  The existence of the decomposition fails at $k=s+1$ only if $\beta_{s+1} > 1$ and there is an index $j$ such that $\beta_j^a  \cdot ... \cdot \beta_s^a \cdot  \beta_{s+1} > C_2^{-1}$.   Then there must be an index $o \geq j$ such that $\beta_o^a > \beta_o$ as otherwise this contradicts $\tilde \theta_{s+1} / \tilde \theta_{j-1} \leq C_2^{-1}.$  If there are multiple such indeces $o$, then consider the largest one.  There must then also be an index $q$ such that $\beta_q^a \cdot ... \cdot \beta_o^a =  C_2$.  There are two cases to consider: $q < j$ and $q \geq j$.  Consider the first case $q < j$.  In this case, the above conclusions can be visualized by: 
$$\phantom{  \beta_q^a \cdot ... \cdot \beta_{j-1}^a } \overbrace{ \phantom{\beta_j^a \cdot ... \cdot \beta_o^a \beta_{o+1} \cdot ... \cdot \beta_{s+1}.} }^ {> C_2^{-1}}$$
\vspace{-12mm}
$${\underbrace{ \underbrace{ \beta_q^a \cdot ... \cdot \beta_{j-1}^a  \beta_j^a \cdot ... \cdot \beta_o^a }_{= C_2} \underbrace{\beta_{o+1} \cdot ... \cdot \beta_{s+1} }_{\leq C_2^{-1}}}_{\leq 1}}.$$
 These imply that $\beta_q^a \cdot ... \cdot \beta_{j-1}^a < C_2$ which contradicts the inductive hypothesis.  The case $q \geq j$ is similar.  This completes the inductive argument and therefore establishes the decomposition $\beta_k = \beta_k^a \beta_k^b$, $k = 2,...,s_0$, for all $s_0$.

Using the fact that $\beta_k^b \leq 1$ for all $k$ allows the definition of the following autoregressive process. 
Let $\mathsf U_1 \sim \N(0,1)$ and let $\mathsf W_1 = \mathsf U_1$.  
Define $\mathsf U_k \sim \N(0,1)$ independently of previous random variables.  Define $\mathsf W_k$ inductively as $$\mathsf W_{k} = \beta_k^b \cdot \mathsf W_{k-1} + \sqrt{ 1 -  ({\beta_k^b})^2} \cdot \mathsf U_{k }.$$   
Note that $\Ep[\mathsf W_k^2] = 1$ and $\Ep[\mathsf W_k \mathsf W_l] = \prod_{q=k+1}^{l} \beta_{q}^b$ if $k<l$.  Then set $\mathsf X_k, \mathsf Y_l$ as follows: 
\begin{align*}
\mathsf X_k &= C_2 \mathsf W_k  \( \prod_{q=2}^{k} \beta_q^a  \)^{-1/2} \( \prod_{q=k+1}^{s_0} \beta_q^a \)^{1/2}  \\
\mathsf Y_l &= C_2 \mathsf W_l  \( \prod_{q=l+1}^{s_0} \beta_q^a  \)^{-1/2} \( \prod_{q=2}^{l} \beta_q^a \)^{1/2}  .
\end{align*}
By construction,  $$\Ep[\mathsf X_k \mathsf Y_l]  = C_2^2 B_{kl}, \ \text{for}\  k \leq l. $$

Next, note that $\Ep[ \mathsf X_k^2 ] \leq 1 $ and  $\Ep[ \mathsf Y_l^2 ] \leq 1 .$   This then implies (taking $\mathsf H_1$ to be the span of $\mathsf U_1,...\mathsf U_{s_0}$ within the set of square integrable random variables)  that  both $$ \Ep [\mathsf X \mathsf Y']  \in \mathscr G_{s_0} \text{ and } \Ep [\mathsf X \mathsf Y']'  \in \mathscr G_{s_0}.$$ 

Take $\bar Z =  \Ep [\mathsf X \mathsf Y']' $.   Let $C_3 = { C_2^{-2}}$.  Note $\Gamma$ is upper triangular due to the way $\tilde \gamma_j$ are defined.  Because $\Gamma$ is upper triangular, only lower triangular components of $\Ep[\mathsf X \mathsf Y']'$ matter for computing the product  $\Gamma C_3 \bar Z$.  Using this fact and the above calculations gives the desired factorization $$\hspace{5.2cm} \Gamma B = \Gamma C_3 \bar Z = \Gamma C_3 \Ep[ \mathsf X \mathsf Y']'. \hspace{4.6cm}   $$

%\subsubsection{Proof of Lemma 7}

%Lemma 6 is Grothendieck's inquality.  

\subsubsection{Proof of Lemma 8}

%Associate to each false covariate $\tilde w_j$, a vector $\tilde \gamma_j \in \mathbb R^{s_0}$, defined as the solution in $\mathbb R^{s_0}$ to the following equation $$\tilde V  \tilde \gamma_j=\tilde r_j.$$  Set $\tilde \gamma_{j\tilde \eps} = \langle \tilde \eps, \tilde w_j \rangle$.   Assume without loss of generality that each component of $\tilde \theta$ is positive (since otherwise, the true covariates can just be multiplied by $-1$.)  Also assume without loss of generality that $\tilde \gamma_{j}'\tilde \theta \geq 0$.  

%Denote the components of $\tilde \gamma_j$ by $\tilde \gamma_j = [\tilde \gamma_{j1}, \tilde \gamma_{j2},...,\tilde \gamma_{j{s_0}}, \tilde \gamma_{j \tilde \eps}]'$ and the components of $\tilde \theta$ by $[\tilde \theta_1,...,\tilde \theta_{s_0}, \tilde \theta_{\tilde \eps}]'$.  For notation purposes, when necessary, the following are considered identical: $\tilde \gamma_{j \tilde \eps} = \tilde \gamma_{j (s_0+1)}, \tilde \theta_{\tilde \eps} = \tilde \theta_{s_0+1}$.

Collect the $m_1$ false selections into
 $\tilde W = [\tilde w_{j_1},...,\tilde w_{j_{m_1}}]$.   Set $\tilde R = [\tilde r_{j_1},...,\tilde r_{j_{m_1}}], \tilde U = [\tilde u_{j_1},...,\tilde u_{j_{m_1}}]$.  Decompose $\tilde W= \tilde R+ \tilde U$. Then $\langle \tilde W,\tilde W \rangle= \langle \tilde  R, \tilde R \rangle  +  \langle \tilde U,  \tilde U \rangle$.  Here, the objects $\langle \tilde W,\tilde W \rangle$, $\langle \tilde R, \tilde R \rangle$ and $\langle \tilde U,  \tilde U \rangle$ etc are formally defined as $m_1 \times m_1$ real matrices with $k,l$ entry given by $\langle \tilde w_k, \tilde w_l \rangle$ $\langle \tilde r_k, \tilde r_l \rangle$ $\langle \tilde u_k, \tilde u_l \rangle$ etc (which, note, are genuine inner products on $\mathsf H$).
 
 \textcolor{black}{Next, by the above normalization}  \textcolor{black}{$diag(\langle \tilde U,  \tilde U \rangle) = I$ if $\langle \tilde u_j , \tilde u_j \rangle = 1$ for all $j \in A_1$.  Recall that this normalization is possible provided $\varphi_{\min}(m+s_0)(G_{\mathsf H}) > 0$.}   Since $diag( \langle \tilde U,  \tilde U \rangle ) = I$, it follows that the average inner product between the $\tilde u_j$, given by $$\bar \rho = \frac{1}{m_1(m_1-1)} \sum_{j \neq l \in A_1} \langle \tilde u_j ,\tilde u_l \rangle,$$ must be bounded below by $$\bar \rho \geq -\frac{1}{{m_1}-1}$$ due to the positive definiteness of $\langle \tilde U,  \tilde U \rangle $.   \textcolor{black}{(This can be checked as a direct consequence of the fact that $1_{m_1 \times 1} ' \langle \tilde U,  \tilde U \rangle {1_{m_1\times 1}} \geq 0$)}. This implies an upper bound on the average off-diagonal term in $\langle \tilde R , \tilde R \rangle $ since $ \langle \tilde W, \tilde W \rangle$ is a diagonal matrix.  Since $\tilde v_k$ are orthonormal, the sum of all the elements of $ \langle \tilde R, \tilde R\rangle $ is given by $ \| \sum_{j\in A_1} \tilde \gamma_j \|_2^2$.  Since $  \| \sum_{j\in A_1} \tilde \gamma_j   \|^2_2 = \sum_{j\in A_1} \| \tilde \gamma_j \|^2_2+ \sum_{j \neq l \in A_1} \tilde \gamma_j'\tilde \gamma_l$ and 
since $ \langle \tilde W,  \tilde W \rangle $ is a diagonal matrix, it must be the case that $$\frac{1}{m_1(m_1-1)}\sum_{j \neq l \in A_1} \tilde \gamma_j'\tilde \gamma_l = -\bar \rho.$$  Therefore,  
$$\textcolor{black}{-}\bar \rho = \frac{1}{{m_1}(m_1-1)} \( \Big \| \sum_{j\in A_1} \tilde \gamma_j  \Big \|^2_2 - \sum_{j \in A_1} \| \tilde \gamma_j\|^2_2 \)  \leq \frac{1}{m_1-1}.$$
This implies that $$ \Big \| \sum_{j \in A_1}  \tilde \gamma_j  \Big  \|^2_2 \leq {m_1}+ \sum_{j\in A_1} \| \tilde \gamma_j\|^2_2 . $$
%Note that $$\| \tilde \gamma_j\|_2^2  \leq  1 + \varphi_{\min}(\hat s + s_0)(G)^{-1}$$ which follows from the fact that $\tilde \gamma_j$ is a regression coefficient of $\tilde r_j$, a vector such that $\tilde r_j'\tilde r_j = 1 - \tilde u'\tilde u \leq 1 - \varphi_{\min}(\hat s + s_0)(G)^{-1}$, on an orthonormal set $\tilde V$. 
Next, bound $\max_{j \in A_1} \| \tilde \gamma_j \|_2^2.$

Note $\| \tilde \gamma_{j} \|_2^2 =  \| \tilde r_j\|_{\mathsf H}^2$ since $\tilde V$ is orthonormal.   \textcolor{black}{ Note that $ \| \tilde w_j \|_{\mathsf H}^2$ is upper bounded by $\varphi_{\min}(m+s_0)(G)^{-1}$.  To see this, note that $ \| \tilde w_j \|_{\mathsf H}^2 = \| c_j \mathscr Q_{\text{pre-}j} w_j \|_{\mathsf H}^2 \leq c_j^2 \| w_j \|_{\mathsf H}^2 = c_j^2 n$ where $c_j$ is the normalizing constant such that $\tilde w_j = c_j \mathscr Q_{\text{pre-}j}$.  At the same time, $c_j^2$ satisfies $ \| \mathscr M_{s_0} \mathscr Q_{\text{pre-}j} w_j \|_{\mathsf H}^2 = c_j^{-2}$ whenever $w_j \notin \text{span}(\tilde V)$.  Note also that  $ \| \mathscr M_{s_0} \mathscr Q_{\text{pre-}j} w_j \|_{\mathsf H}^2  \geq \| \mathscr Q_{S_0 \cup \text{pre-}j} w_j \|_{\mathsf H}^2$ where the notation $\mathscr Q_{S_{0} \cup \text{pre-}j}$ denotes projection onto the space orthogonal to covariates indexed in $S_0$ or selected before $w_j$.   To see this, consider an arbitrary Hilbert space $\check{\mathsf H}$,  projections onto closed subspaces $1,2, 12 = \text{span}(1 \cup 2)$ , $\mathscr P_1, \mathscr P_2, \mathscr P_{1 2}$, projections onto the respective orthogonal complements $\mathscr Q_{1}, \mathscr Q_{2}, \mathscr Q_{12 },$  and any vector $w$.  Then 
$w = \mathscr Q_{1 2} w + \mathscr P_{12} w.$
Then 
$\mathscr Q_2 \mathscr Q_1 w   = \mathscr Q_2 \mathscr Q_1 \mathscr Q_{12} w + \mathscr Q_2 \mathscr Q_1 \mathscr P_{12} w
		=\mathscr Q_{12} w +\mathscr Q_2 \mathscr  Q_1 \mathscr  P_{12} w.$
Note that the inner product between the above two terms vanishes: $\langle  \mathscr Q_{12}w,   \mathscr Q_2 \mathscr Q_1 \mathscr P_{12} w \rangle_{\check{\mathsf H}} = \langle w, \mathscr Q_{12} \mathscr P_{12} w \rangle_{\check{\mathsf H}} = \langle w, 0 w \rangle_{\check{\mathsf H}} = 0. $
Then by Pythagorean Theorem,
$\| \mathscr Q_2 \mathscr Q_1 w \|_{\check{\mathsf H}}^2 = \| \mathscr Q_{12} w \|_{\check{\mathsf H}}^2 + \| \mathscr Q_2  \mathscr Q_1 \mathscr P_{12} w \|_{\check{\mathsf H}}^2  \geq  \| \mathscr Q_{12} w \|_{\check{\mathsf H}}^2.$
So $\| \mathscr Q_{12} w \|_{\check{\mathsf H}} \leq \| \mathscr Q_2 \mathscr Q_1 w \|_{\check{\mathsf H}}.$
}
\textcolor{black}{
 Therefore, the quantity $\|\mathscr Q_{S_0 \cup\text{pre-}j}w_j\|_\H^2$ is lower bounded by $n\varphi_{\min}(m+s_0)(G_{\mathsf H})$.  As a result, $c_j^2 \leq \varphi_{\min}(m+s_0)(G_{\mathsf H})^{-1}$, giving the desired bound on $\|\tilde w_j\|_{\mathsf H}^2$.}   Therefore, $ \| \tilde r_j\|_{\mathsf H}^2  = \| \tilde w_j \|_{\mathsf H}^2 - 1 \leq \varphi_{\min}(m+s_0)(G_\H)^{-1} - 1$.  It follows that 
$$\max_{j \in A_1} \| \tilde \gamma_j \|_2^2 \leq  \varphi_{\min}(m+s_0)(G_\H)^{-1} -1. $$
 This then implies that $$\Big \| \sum_{j \in A_1}  \tilde \gamma_j  \Big \|^2_2 \leq {m_1} \varphi_{\min}({m}+s_0)(G_{\mathsf H})^{-1}.$$
The same argument as above also shows that for any choice $e_j \in \{ -1, 1 \}$ of signs,  it holds that 
$$\Big \| \sum_{j \in A_1} e_j \tilde \gamma_j  \Big \|^2_2 \leq  {m_1} \varphi_{\min}({m}+s_0)(G_{\mathsf H})^{-1}.$$
(In more detail, take $\tilde W_e = [\tilde w_{j_1}e_{j_1},...,\tilde w_{j_{m_1}} e_{j_{m_1}}],$ etc. and rerun the same argument.)

\subsubsection{Proof of Lemma 10}   In this \textcolor{black}{proof}, the number of elements of $A_2$ is bounded. 
Recall that the criteria for $j\in A_2$ is that $|\tilde \gamma_{j \tilde \eps} |>  \frac{t_\H^{1/2}n^{1/2}}{(\textcolor{black}{3}\langle \eps , \mathscr M_{s_0} \eps \rangle )^{1/2}}  $.  
Note also that $\tilde \gamma_{j\tilde \eps}$ is found by the coefficient in the expression
$$ \tilde \gamma_{j\tilde \eps} =   \langle  \tilde \eps, \tilde w_j \rangle=  \langle\eps,  \frac{1}{ \langle \eps , \mathscr M_{s_0} \eps \rangle ^{1/2}}   \mathscr M_{s_0} \tilde  w_j \rangle.$$
Next, let $H$ be $H = [v_1,...,v_{s_0}, w_1,...,w_m]$.  Note that $$\frac{1 }{ \langle \eps , \mathscr M_{s_0} \eps \rangle ^{1/2}}   \mathscr M_{s_0} \tilde  w_j   \in \text{span}(H),$$
%which implies that the above expression is unchanged when premultiplied by $H(H'H)^{-1}H'$.  
Therefore, $$ \tilde \gamma_{j\tilde \eps} = \langle \eps , H \rangle  \langle H,H\rangle ^{-1} \langle H, \frac{1}{(\langle \eps , \mathscr M_{s_0} \eps \rangle )^{1/2} }  \mathscr M_{s_0} \tilde  w_j \rangle. $$
Let $\mu_j$ be the \textcolor{black}{sign}  $+1$ for each $j \in A_2 $ such that $\tilde \gamma_{j\tilde \eps} > 0$ and $-1$ for each $j\in A_2$ such that $\tilde \gamma_{j\tilde \eps} < 0$. By the fact that $j \in A_2$, $\tilde \gamma_{j\tilde \eps}\mu_j > \frac{t_\H^{1/2} n^{1/2} }{(3\langle \eps , \mathscr M_{s_0} \eps \rangle )^{1/2}},$  summing over $j \in A_2$ gives $$\sum_{j \in A_{2}} \langle \eps, H \rangle  \langle H,H\rangle ^{-1} \langle H,  \frac{1}{( \langle \eps , \mathscr M_{s_0} \eps \rangle )^{1/2}} \mathscr M_{s_0} \tilde  w_j  \mu_j \rangle> m_2 \frac{t_\H^{1/2}n^{1/2}}{ (\textcolor{black}{3}\langle \eps , \mathscr M_{s_0} \eps \rangle  )^{1/2}}.$$
This implies that 
$$\Big \|  \langle H,H\rangle ^{-1} \langle H,   \frac{1 }{( \langle \eps , \mathscr M_{s_0} \eps \rangle )^{1/2}}  \sum_{j \in A_{2}} \mathscr M_{s_0} \tilde  w_j  \mu_j \rangle  \Big \|_1 \| \langle \eps, H \rangle \|_\infty > m_2  \frac{t_\H^{1/2}n^{1/2} }{ (\textcolor{black}{ 3}\langle \eps , \mathscr M_{s_0} \eps \rangle  )^{1/2}} $$
Which further implies that 
$$\sqrt{m+s_0} \Big \|  \langle H,H\rangle ^{-1}\langle  H,  \frac{1}{( \langle \eps , \mathscr M_{s_0} \eps \rangle )^{1/2}}  \sum_{j \in A_{2}} \mathscr M_{s_0} \tilde  w_j  \mu_j    \rangle \Big \|_2 \| \langle \eps, H \rangle  \|_\infty > m_2 \frac{t_\H^{1/2}n^{1/2}}{(\textcolor{black}{ 3}\langle \eps , \mathscr M_{s_0} \eps \rangle )^{1/2}} $$
Next,  further upper bound the $\| \cdot \|_2$ term on the left side above by 
$$ \Big \|  \langle H,H\rangle ^{-1} \langle H,  \frac{1}{( \langle \eps , \mathscr M_{s_0} \eps \rangle )^{1/2}} \rangle   \sum_{j \in A_{2}} \mathscr M_{s_0} \tilde  w_j  \mu_j \rangle   \Big \|_2  $$ $$\leq  \frac{n^{-1/2} }{( \langle \eps , \mathscr M_{s_0} \eps \rangle )^{1/2}}   \varphi_{\min}(s_0+m)(G_\H)^{-1/2}   \|  \mathscr M_{s_0} \sum_{j \in A_{2}} \tilde  w_j  \mu_j  \|_\H.$$  
Next, by the fact that $\mathscr M_{s_0}$ is a projection (hence non-expansive)  and $\tilde w_j$ are mutually orthogonal, 
%$$ \leq    \frac{t^{-1/2}}{( \eps ' \mathscr M_{s_0} \eps)^{1/2}} \sigma_{\max}((H'H)^{-1}H')  \| \sum_{j \in A_{2}} \tilde  w_j  \mu_j  \|_2 $$
$$ \leq \frac{n^{-1/2}}{(\langle \eps , \mathscr M_{s_0} \eps \rangle )^{1/2}}   \varphi_{\min}(s_0+m)(G_\H)^{-1/2} \sqrt{ \sum_{j \in A_{2}}  \| \tilde  w_j  \mu_j  \|_\H^2} .$$
Earlier, it was shown that $\max_j \| \tilde w_j\|_{\mathsf H}^2 \leq \varphi_{\min}(s_0 + m)(G_\H)^{-1}$.   Therefore, putting the above inequalities together, 
$$\frac{n^{-1/2} }{(\langle \eps , \mathscr M_{s_0} \eps \rangle )^{1/2}}  \sqrt{m+s_0} \varphi_{\min}(m+s_0)(G_\H)^{-1} \sqrt{m_2}\| \langle \eps, H\rangle \|_{\infty} > m_2  \frac{t_\H^{1/2}n^{1/2}}{ (\textcolor{black}{ 3}\langle \eps , \mathscr M_{s_0} \eps \rangle  )^{1/2}}  .$$
This implies that 
$$m_2 < \frac{1}{n^{\textcolor{black}{2}}}\frac{\textcolor{black}{3}}{t_\H} (\langle \eps, \mathscr M_{s_0} \eps \rangle)(m+s_0)  \frac{ \| \langle \eps, H \rangle \|_{\infty}^2}{ \eps ' \mathscr M_{s_0} \eps} \varphi_{\min}(m+s_0)(G_\H)^{-2}$$ 
In the case of Theorem 1, this is further bounded by $$\leq \textcolor{black}{ 3}(m +s_0)\frac{ \|  \En[ x_i \eps_i ]\|_{\infty}^2}{t} \varphi_{\min}(m+s_0)(G)^{-2}.$$
under the assumed condition that $t^{1/2}  \geq 2 \| \En[x_i \eps_i ]\|_{\infty} \varphi_{\min}(m+s_0)(G)^{-1}$,  it follows that
$$m_2 \leq \textcolor{black}{\frac{3}{4} }(m+s_0).$$
Similarly, the condition of Theorem 4 that $\Ep[\En[\eps^{\text{a}2}_i]] \leq \frac{1}{2} \varphi_{\min}(\Ep[G])^{-1}c_{\text{test}}')$ yields $m_2 \leq \textcolor{black}{\frac{3}{4} }(m+s_0)$ in the same way.
Finally, substituting $m = m_1 + m_2$ gives 
$$\hspace{6.2cm}m_2 \leq  \textcolor{black}{3} m_1 + \textcolor{black}{3} s_0. \hspace{5.2cm}  $$

\subsubsection{Proof of Lemma 11}  
Combining 
$m_1 \leq   \varphi_{\min}(m+s_0)(G_\H)^{-1}C_1^{-2} {C_3}^2 \(K_G^{\mathbb R} \)^2s_0 $
and
$m_2 \leq  \textcolor{black}3 (m_1 + s_0)$
gives 
$$m \leq  \[4  \varphi_{\min}(m+s_0)(G_\H)^{-1} C_1^{-2} {C_3}^2 \({K_G^{\mathbb R} }\)^2   
 + 3 \]s_0.$$  In addition, in the case of Theorem 1,
% \begin{align*} C_1 &= \frac{1}{2} \varphi_{\min}(m + s_0)(G), \\
%C_2 &=  \varphi_{\min}(m + s_0)(G), \\  
%C_3 &= (1 +  \varphi_{\min}(m + s_0)(G)^{-1/2} +  \varphi_{\min}(m + s_0)(G)^{-1}), 
%\end{align*}
 \textcolor{black}{$ C_1^2 = \frac{1}{6} \varphi_{\min}(m + s_0)(G_\H), \
C_2^2 =  \varphi_{\min}(m + s_0)(G_\H) , \  
C_3^2  = (C_2^{-2})^2 = \varphi_{\min}(m + s_0)(G_\H)^{-2}, \
C_1^{-2}C_3^2  = 6 \varphi_{\min}(m+s_0)(G_\H)^{-3} 
$
and $K_G^{\mathbb R} < 1.783$.  Therefore, $m  \leq (3 + \textcolor{black}{24} \times 1.783^2  \times  \varphi_{\min}(m+s_0)(G_\H)^{-4} \textcolor{black}{ }  ) s_0. $}
%Since $C_3^2 \leq 9 \varphi_{\min}(s_0+m)(G)^{-2}$, the expression above can be simplified at the expense of a slightly less tight constant, so that
%$$m \leq [ 1 + 8 \times 1.783^2  \times  \varphi_{\min}(m+s_0)(G)^{-5}   ]s_0.$$
\noindent Because $\varphi_{\min}(m+s_0)(G_\H)^{-1} \geq 1$ and $ 3 + 24 \times 1.783^2 = 79.2981 < 80$, it holds that $$m \leq 80 \times \varphi_{\min}(m+s_0)(G_\H)^{-4} s_0.$$
This bound holds for each positive integer $m$ of wrong selections, provided ${t^{1/2}} \geq  2\varphi_{\min}(m+s_0)(G)^{-1} \| \En[x_i \eps_i]\|_{\infty}$.  This concludes the proof of the sparsity bound for Theorem 1.   
Using similar reasoning in the case of Theorem 4, on the event $\mathscr T$,  it follows that $m \leq  80 \times \varphi_{\min}(m+s_0)(G_\H)^{-4} c_{\text{test}}''^{-3}s_0$ provided $\Ep[\En[\eps_i^{\text{a2}}]] \leq \frac{1}{2} \varphi_{\min}(m+s_0)(\Ep[G])^{-1} c_{\text{test}}'$.  Setting $m = K_{\text{test}}-s_0$ contradicts Condition 2 by   $K_{\text{test}} \leq 80 \times \varphi_{\min}(K_{\text{test}})(\Ep[G])^{-4}{c_{\text{test}}''}^{-3} + s_0 < K_{\text{test}}$.  Therefore, $m < K_{\text{test}}-s_0$ and thus
%no covariates beyond $K_{\text{test}}-s_0$ steps can be selected, which immediately results in 
$$\hat s \leq (80 \times \varphi_{\min}(K_{\text{test}})(G_\H)^{-4} c_{\text{test}}''^{-3}+1)s_0,$$ completing the proof of the sparsity bound for Theorem 4.

%\section{Proof of Theorem 4}

%The proof is similar to the proof of Theorem 1.  The only modification necessary is in Section 6.   The set $A_2$ can be further split into $A_{21} = \{ j \in A_2: \En[w_{ij}\eps_i] > \frac{1}{2}\varphi_{\min}(m+s_0)(G), \  A_{22} = A_2 \setminus A_{21}$.  On the set $A_{21}$, the same logic as in the proof of Theorem 1 in Section 6 show that $|A_{21}| \leq \frac{1}{2}(m_1+ |A_{21}| + |A_{22}|+s_0)$.  Then it follows that $|A_{21} | \leq |A_{22}| + m_1 + s_0$.   Since $|A_{22}|$ is bounded by $q(m)$ by assumption, the theorem follows.

\subsection{Proof of Theorem 5}

The strategy is to apply Theorem 4 using the conditional distribution $\Pr_x$ for $\mathscr D_n$, conditional on $x$.  The unconditional result is then shown to follow.
Let $\E_x(S) = \Ep[\ell(S)|x]$.   In addition, \textcolor{black}{for $j \notin S$} let $\theta_{jS}^{*|x} = (x_{jS}' x_{jS})^{-1} x_{jS}' \Ep[ x_{jS}' (x \theta_0 + \eps^{\text{a}})|x]$ so that $[\theta^{*|x}_{jS}]_j = (x_{j}' {\mathscr Q}_S x_j)^{-1} \Ep[x_j' {\mathscr Q}_S (x \theta_0 + \eps^{\text{a}}) | x]$.
\textcolor{black}{Throughout the proof of Theorem 5, use an abuse of notation by writing $\hat V_{jS} = [\hat V_{jS}]_{jj}$.  }
Let 
$$ \hat Z_{jS} =  \hat V_{jS}^{-1/2}([\hat \theta_{jS}]_j - [ \theta_{jS}^{*|x}]_j).$$ Let $t_{\alpha} = \Phi^{-1}(1-\alpha/p)$.  Let $\mathscr A$ be the event given by $$ \mathscr A = \left \{ | \hat Z_{jS}| \leq \left (\frac{1 + c_{\tau}}{2}\right ) \hat \tau_{jS} t_{\alpha} \ \text{for all} \ j , |S| < K_n \right \}.$$
Note that $-\Delta_j \mathscr E_x(S) = [\theta_{jS}^{*|x}]_j^2 A_{jS}$ for $A_{jS}$ defined by $A_{jS} = [G_{jS}^{-1}]_{jj}$.  

The next lemma states size, power, and continuity properties of the tests of Definition 1.
\begin{lemma}The following implications are valid on $\mathscr A$ for all $j,|S|<K_n$:

\begin{itemize}
\item[1.] $ \ T_{jS\alpha} = 1 \text{ \ if \ } -\Delta_{j} \mathscr E_x(S) \geq A_{jS} \hat V_{jS} (\textcolor{black}{2c_{\tau}})^2 \hat \tau_{jS}^2  t_\alpha^2  .$

\item[2.] $ -\Delta_{j} \mathscr E_x(S) \geq  A_{jS} \hat V_{jS} \left ( \frac{1-c_\tau}{2} \right )^2 \hat \tau_{jS}^2  t_\alpha^2  \text{ \ if \ }T_{jS\alpha}=1.$

\item[3.]  $  -\Delta_k\E_x(S) \leq \frac{\hat V_{kS}A_{kS}   }{\hat V_{jS}A_{jS}  }  \(1+ \frac{1 + c_\tau}{ c_\tau - 1}\( 1+\frac{\hat \tau_{kS}   }{  \hat \tau_{jS}   }\)  \)^2  (-\Delta_j\E_x(S))
$ if $ T_{jS\alpha}=1$,  $W_{jS} \geq W_{kS}$.
\end{itemize}
\end{lemma}

\

Next define a sequence of sets $\mathscr X = \mathscr X_n$ which will be shown to have the property that  both $\Pr(x \in \mathscr X) \rightarrow 1$ and 
$$\Pr^{\mathscr X}( \mathscr A) = \text{ess }  \inf_{x \in \mathscr X} \Pr(\mathscr A | x) \rightarrow 1.$$

In addition, there will be constants $\tilde c_{\text{test}}, \tilde c_{\text{test}}', c_{\text{test}}'' >0$ which are independent of $n$ and the realization of $x$, such that for $\ctest= \frac{1}{n} \tilde c_{\text{test}} ,$ $\ctest' = \frac{1}{n} \tilde c_{\text{test}}'$ and for the set $\mathscr B$ defined by
$$ \mathscr B = \begin{cases} 
\ \ 1. \ \ A_{jS}\hat V_{jS} (\textcolor{black}{2c_\tau})^2 \hat \tau_{jS}^2t_{\alpha}^2 \leq \ctest \\ \\
\ \ 2. \ \ A_{jS}\hat V_{jS} \left ( \frac{1-c_\tau}{2} \right )^2 \hat \tau_{jS}^2t_{\alpha}^2  \geq \ctest'  \ \ \ \ \ \  \ \ \ \ \ \ \ \ \ \ \ \ \ |S| < K_n \\ \\
\ \ 3. \ \ \frac{A_{kS} \hat V_{kS}}{A_{jS} \hat V_{jS}} \(1+ \frac{1+c_\tau }{1-c_{\tau}  } \( 1+\frac{\hat \tau_{kS}   }{  \hat \tau_{jS}   }\) \)^2 \geq \ctest''\end{cases} $$
it holds that $\Pr^{\mathscr X}(\mathscr B) \rightarrow 1$.

Define sets $\mathscr X = \mathscr X_{n}$ as follows.  Set
$\mathscr X = \mathscr X_1 \cap \mathscr X_2 \cap \mathscr X_3 \cap \mathscr X_4$ with
\begin{itemize} 
\item[] $\mathscr X_1 =\{ x: \max_{j\leq p } \En[ x_{ij}^{12} ] = O(1) \} $
\item[] $\mathscr X_2 = \{ x: \varphi_{\min}(K_n)(G)^{-1} = O(1)\} $
\item[] $\mathscr X_3 = \{ x:\textcolor{black}{\max_{j, |S| < K_n}} \| \eta_{jS} \|_1 = O(1) \} $
\item[] $\mathscr X_4 = \{ x : \Pr( \varphi_{\min}(K_n)( \En[ \eps_i^2 x_ix_i' ] )^{-1} = O(1) |x )= 1-o(1) \}.$
\end{itemize}

Note that $\Pr(\mathscr X_1), \Pr(\mathscr X_2), \Pr(\mathscr X_3) \rightarrow 1$ by assumption.  In addition, failure of $\Pr( \mathscr X_4 ) \rightarrow 1$ would contradict the unconditional statement in Condition 4 that $\Pr( \varphi_{\min}(K_n)( \En[ \eps_i^2 x_ix_i' ] )^{-1} = O(1) )= 1-o(1) $.  Therefore, $\Pr( \mathscr X) \rightarrow 1$.

The next two sections prove the following two lemmas.

\begin{lemma}$\Pr^{\mathscr X}( \mathscr A) \rightarrow 1$.\end{lemma}
\begin{lemma} $\Pr^{\mathscr X}( \mathscr B) \rightarrow 1$ for some $c_{\text{\em test}}, c_{\text{\em test}}',c_{\text{\em test}}''$ as described in the definition of $\mathscr B$ above.  \end{lemma}
%Note that these sections allow the application of Theorem 4 conditionally on $x$.  The final section then uses this fact, and then concludes the proof by showing that $\theta_0$ is bounded to $\theta_{S_0}^{*|x}$.

The previous results show that for each $n$, Theorem 4 can be applied conditionally on $x$ with $c_{\text{test}}, c_{\text{test}}',c_{\text{test}}''$ defined above, with $K_{\text{test}}=K_n-1$, and with $1-\alpha -\delta_{\text{test}} = \Pr^{\mathscr X}(\mathscr A \cap \mathscr B)$.  Note that renormalizing the covariates to satisfy $\En[x_{ij}^2] = 1$ does not affect $\mathscr E_x(S)$ and therefore does not affect the conclusions above.  Moreover, on $\mathscr X$, renormalizing does not affect boundedness of sparse eigenvalues of $G$.  The unconditional result is shown as follows.  By Theorem 4,
$$\Pr^{\mathscr X}\(\En [ (x_i' \theta_{S_0}^{*|x} - x_i \hat \theta)^2 ]^{1/2} \leq O(\sqrt{s_0 \log p /n})\) \rightarrow 1.$$  
%This follows from assumptions on $\eps^{\text{a}}$.  
Note that $\theta_{S_0}^{*|x} - \theta_0 = (x_{S_0}' x_{S_0})^{-1} x_{S_0}' \Ep[ \eps^{\text{a}}|x].$  As a result, $$\| \theta_0 - \theta_{S_0}^{*|x} \|_2 \leq \varphi_{\min}(s_0)(G)^{-1/2} \| \En[x_{is_0} \Ep[ \eps_i^{\text{a}}|x]] \|_2 \leq  \varphi_{\min}(s_0)(G)^{-1/2} \sqrt{s_0} \| \En[x_{ij} \Ep[ \eps_i^{\text{a}}|x]] \|_\infty .$$  By the assumed rate conditions, sparse eigenvalue conditions, and by $\max_i \Ep[\eps_i^{\text{a}}] = O(n^{-1/2}),$ the bound on $\| \theta_0 - \theta_{S_0}^{*|x} \|_2 $ implies further that $\Pr^{\mathscr X}\(\En [ (x_i' \theta_{S_0}^{*|x} - x_i \theta_0)^2 ]^{1/2} \leq O(\sqrt{s_0 \log p /n})\) \rightarrow 1.$   Theorem 5 follows by triangle inequality.  

\subsection{Proof of Supporting Lemmas for Theorem 5}

\subsubsection{Proof of Lemma 12}

For this proof, work on $\mathscr A$ and suppose $|S| < K_n$.  To prove the first statement, suppose that  $ -\Delta_{j} \mathscr E_x(S) \geq A_{jS} \hat V_{jS} (\textcolor{black}{2c_{\tau}} )^2 \hat \tau_{jS}^2  t_\alpha^2 .$  
Then 
\begin{align*} [\theta_{jS}^{*|x}]_j ^2 A_{jS} & \geq A_{jS} \hat V_{jS} (\textcolor{black}{2c_{\tau} })^2 \hat \tau_{jS}^2  t_\alpha^2 \\
| [\theta_{jS}^{*|x}]_j | & \geq  \hat V_{jS}^{1/2} (\textcolor{black}{2c_{\tau}} ) \hat \tau_{jS}  t_\alpha \\
| [\hat \theta_{jS}]_j | & \geq  \hat V_{jS}^{1/2} (\textcolor{black}{2c_{\tau}} ) \hat \tau_{jS}  t_\alpha - | [\theta_{jS}^{*|x}]_j  -   [\hat \theta_{jS}]_j |\\
| [\hat \theta_{jS}]_j | & \geq  \hat V_{jS}^{1/2} (\textcolor{black}{2c_{\tau}} ) \hat \tau_{jS}  t_\alpha  -\hat V_{jS}^{1/2} \( \frac{1 + c_{\tau}}{2} \) \hat \tau_{jS} t_{\alpha} \\
| [\hat \theta_{jS}]_j |&  \geq  \hat V_{jS}^{1/2} c_{\tau}  \hat \tau_{jS}  t_\alpha 
\end{align*}
which implies $T_{jS\alpha}=1$. 

%\subsubsection{Size Calculations}

Next, prove the second statement.  By construction, if $T_{jS\alpha} = 1 $ then  $| \hat V_{jS}^{-1/2} [ \hat \theta_{jS}]_j | \geq c_\tau \hat \tau_{jS} t_\alpha$, which is equivalent to 
$$  | [\hat \theta_{jS}]_j |\geq c_{\tau} \hat \tau_{jS} t_{\alpha}\hat V_{jS}^{1/2}.$$
Note that $| [\hat \theta_{jS}]_j - [ \theta_{jS}^{*|x}]_j |  \leq \hat V_{jS}^{1/2} \( \frac{1 + c_{\tau}}{2} \) \hat \tau_{jS} t_{\alpha} .$
Then  $T_{jS\alpha} = 1 \Rightarrow$
\begin{align*}  | [ \theta_{jS}^{*|x}]_j| \geq  c_\tau \hat \tau_{jS} t_\alpha \hat V_{jS}^{1/2} -  \hat V_{jS}^{1/2} \( \frac{1 + c_{\tau}}{2} \)\hat \tau_{jS} t_{\alpha} = \hat V_{jS}^{1/2} \hat \tau_{jS} t_{\alpha} \( \frac{c_{\tau} - 1}{2} \)
\end{align*}
Thefore 
 $ -\Delta_{j} \E_x(S)  \geq  A_{jS} \hat V_{jS} \hat \tau_{jS}^2 t_{\alpha}^2 \( \frac{c_{\tau} - 1}{2} \)^2 . $

Finally, prove the third statement.  Note that $  {W_{kS}} \leq {W_{jS}} $ implies $ { \hat V_{kS}^{-1/2} | [\hat \theta_{kS}]_k | } \leq {\hat V_{jS}^{-1/2} | [\hat \theta_{jS}]_j |} $. Then

$$    {  \hat V_{kS}^{-1/2} | [ \theta_{kS}^{*|x}]_k |  - \( \frac{1 + c_{\tau}}{2}\) \hat \tau_{kS} t_{\alpha} } \leq { \hat V_{jS}^{-1/2} | [ \theta_{jS}^{*|x}]_k |  + \( \frac{1 + c_{\tau}}{2}\) \hat \tau_{jS} t_{\alpha}}   $$

$$ \Rightarrow   {  \hat V_{kS}^{-1/2} | [ \theta_{kS}^{*|x}]_k |  }\leq { \hat V_{jS}^{-1/2} | [ \theta_{jS}^{*|x}]_j |  + \( \frac{1 + c_{\tau}}{2}\)(\hat \tau_{kS} +  \hat \tau_{jS} )t_{\alpha}}   $$
$$ \Rightarrow  {  \hat V_{kS}^{-1/2} A_{kS}^{-1/2} (-\Delta_k\E_x(S))^{1/2}  }\leq { \hat V_{jS}^{-1/2} A_{jS}^{-1/2} (-\Delta_j\E_x(S))^{1/2}  + \( \frac{1 + c_{\tau}}{2}\) (\hat \tau_{kS} + \hat \tau_{jS} )t_{\alpha}}   $$
$$ = { \hat V_{jS}^{-1/2} A_{jS}^{-1/2} (-\Delta_j\E_x(S))^{1/2}  + \( \frac{1 + c_{\tau}}{2}\) (\hat \tau_{kS} + \hat \tau_{jS} )t_{\alpha}} \( \frac{   A_{jS} \hat V_{jS} \left ( \frac{1-c_\tau}{2} \right )^2 \hat \tau_{jS}^2  t_\alpha^2   } {  A_{jS} \hat V_{jS} \left ( \frac{1-c_\tau}{2} \right )^2 \hat \tau_{jS}^2  t_\alpha^2} \)^{1/2} . $$

Using the fact that $-\Delta_{j} \mathscr E_x(S) \geq  A_{jS} \hat V_{jS} \left ( \frac{1-c_\tau}{2} \right )^2 \hat \tau_{jS}^2  t_\alpha^2$ (because $T_{jS\alpha}=1$), gives that the previous expression is bounded by

$$ \leq  \hat V_{jS}^{-1/2} A_{jS}^{-1/2} (-\Delta_j\E_x(S))^{1/2}  + \frac{ \( \frac{1 + c_{\tau}}{2}\) (\hat \tau_{kS} + \hat \tau_{jS} )t_{\alpha}} {\( A_{jS} \hat V_{jS} \left ( \frac{1-c_\tau}{2} \right )^2 \hat \tau_{jS}^2  t_\alpha^2 \)^{1/2}   } (-\Delta_j \mathscr E_x(S))^{1/2} $$

$$ =  \hat V_{jS}^{-1/2} A_{jS}^{-1/2} \(1 +   \frac{ 1 + c_\tau}{c_\tau-1} \frac{\hat \tau_{kS} + \hat \tau_{jS}    }{  \hat \tau_{jS}   }   \)(-\Delta_j\E_x(S))^{1/2}.  $$

%$$ {  \hat V_{jS}^{-1/2} A_{jS}^{-1/2} (-\Delta_j\E_x(S))^{1/2}  }\leq { \hat V_{kS}^{-1/2} A_{kS}^{-1/2} (-\Delta_k\E_x(S))^{1/2}  + \( \frac{1 + c_{\tau}}{2}\) (\hat \tau_{jS} + \hat \tau_{kS} )t_{\alpha}}   $$

%$\( \frac{1+c_{\tau}}{2} \) \hat \tau_{kS} t_{\alpha} \leq \hat V_{kS}^{-1/2} A_{kS}^{-1/2} (-\Delta_k\E_x(S))^{1/2} \( \frac{1 + c_\tau}{ c_\tau - 1} \)$. 

%Then
%$$   {  \hat V_{jS}^{-1/2} A_{jS}^{-1/2} (-\Delta_j\E_x(S))^{1/2}  } \leq { \hat V_{kS}^{-1/2} A_{kS}^{-1/2} (-\Delta_k\E_x(S))^{1/2}  +   2  \hat V_{kS}^{-1/2} A_{kS}^{-1/2} (-\Delta_k \mathscr E_x(S))^{1/2} \( \frac{1 + c_\tau}{ c_\tau - 1} \)  }  . $$

This gives
$  -\Delta_k\E_x(S) \leq \frac{\hat V_{kS}A_{kS}   }{\hat V_{jS}A_{jS}  }  \(1+ \frac{1 + c_\tau}{ c_\tau - 1}\( 1+\frac{\hat \tau_{kS}   }{  \hat \tau_{jS}   }\)   \)^2  (-\Delta_j\E_x(S)).
$

\subsubsection{Proof of Lemma 13}

% It is sufficient to prove that $\Pr(\mathscr A) \rightarrow 1$ unconditionally since this forces the existence of the desired set $\mathscr X$.

%In addition, let $\theta^{*|x,\text{a}}_{jS} = (x_{j} M_S x_j)^{-1}x_j M_S (x \theta_0 + \eps^{\text{a}})$ and 
%Define $\delta_{\text{a}} = \max_{j,|S| < K}| [\theta^{*|x,\text{a}}_{jS} ]_j- [\theta^{*|x}_{jS} ]_j|$.  Note that $$\delta_{\text{a}}= \max_{j,|S| < K} |(x_j' M_S x_j)^{-1}x_j'M_S ( \eps^{\text{a}} - \Ep[\eps^{\text{a}}|x])|$$ and can be bounded by Cauchy-Schwarz inequality to give by
%$$ \max_{j, |S| <K} (x_j' M_S x_j)^{-1/2} (( \eps^{\text{a}} - \Ep[ \eps^{\text{a}} |x]) ' ( \eps^{\text{a}} - \Ep[ \eps^{\text{a}} |x] ) )^{1/2}$$  
%which by assumption is $O( \frac{1}{\sqrt n}) O( 1 )$ with probability $1- o(1)$.

Note that 
\begin{align*}\hat Z_{jS} &=\hat V_{jS}^{-1/2}([\hat \theta_{jS}]_j - [\theta^{*|x}_{jS}]_j) \\
&= \hat V_{jS}^{-1/2} (x_j'{\mathscr Q}_Sx_j)^{-1} x_j'{\mathscr Q}_S ( \varepsilon - \Ep[\eps|x]  )\\
& = (  (x_j'{\mathscr Q}_Sx_j)^{-1}  \En[ \hat \eps_{ijS}^2 [{\mathscr Q}_Sx_{jS}]_i^2]  (x_j'{\mathscr Q}_Sx_j)^{-1} )^{-1/2}  (x_j'{\mathscr Q}_Sx_j)^{-1} x_j'{\mathscr Q}_S ( \varepsilon - \Ep[\eps|x]  ) \\
&=  \En[ \hat \eps_{ijS}^2 [{\mathscr Q}_Sx_{jS}]_i^2]^{-1/2}  x_j' {\mathscr Q}_S ( \varepsilon - \Ep[\eps|x]  )\\
&=  \En[ \hat \eps_{ijS}^2(\eta_{jS}'x_{ijS})^2]^{-1/2}  \eta_{jS}'x_{jS} ( \varepsilon - \Ep[\eps|x]  ).\\
&=  \En[ \hat \eps_{ijS}^2(\eta_{jS}'x_{ijS})^2]^{-1/2}  \eta_{jS}'x_{jS} ( \eps^{\text{o}} + \varepsilon^{\text{a}} - \Ep[\eps^{\text{a}}|x]  ).
\end{align*}
Let $\ddot \eps = \eps^{\text{o}} + \eps^{\text{a}} - \Ep[ \eps^{\text{a}}|x]$.
  Define the \textit{Regularization Event} by
$$ \mathscr R= \left \{\frac{ | \sum_{i=1}^n x_{ik} \ddot \varepsilon_i |}{\sqrt{\sum_{i=1}^n x_{ik}^2 {\ddot \eps_i}^2}} \leq t_{\alpha}  \ \ \text{for every} \ \ k \leq p \ \right \}$$ 

\noindent In addition, define the \textit{Variability Domination Event} by 
$$\mathscr V = \left \{ \sum_{i=1}^n x_{ik}^2 \ddot \varepsilon_i^2  \leq \( \frac{1+c_{\tau}}{2} \)^2 \sum_{i=1}^n x_{ik}^2 \hat \varepsilon_{ijS}^2 \ \ \text{for every} \ \ k \in jS, \ \text{for every} \  |S| < K_n \right \} $$
 The definition of the Regularization Event and the Variability Domination Event are useful because $$ \mathscr R \cap \mathscr V\Rightarrow \mathscr A.$$
To see this, note that on $\mathscr  R$, the following inequality holds for any conformable vector $\nu$:
$$
\(\sum_{i=1}^n \sum_{k \in jS} \nu_k x_{ik}\ddot \varepsilon_i \)^2  \leq \( t_\alpha \sum_{k \in jS} |\nu_k |\sqrt{ \sum_{i=1}^n x_{ik}^2\ddot  \varepsilon_i^2} \hspace{1mm} \)^2 
$$
Furthermore, on $\mathscr V$, the previous expression can be further bounded by
$$
\leq \( \frac{1+c_{\tau}}{2} \)^2 \( t_\alpha \sum_{k \in jS} |\nu_k |\sqrt{ \sum_{i=1}^n x_{ik}^2 \hat \varepsilon_{ijS}^2} \hspace{1mm} \)^2$$
$$ =  \( \frac{1+c_{\tau}}{2} \)^2  \frac{ \( t_\alpha \sum_{k \in jS}| \nu_k| \sqrt{ \sum_{i=1}^n x_{ik}^2 \hat \varepsilon_{ijS}^2} \)^2 }{ \sum_{i=1}^n \(\sum_{k \in jS} \nu_k x_{ik} \)^2 \hat \varepsilon_{ijS}^2    } \sum_{i=1}^n \(\sum_{k \in jS} \nu_k x_{ik} \)^2 \hat \varepsilon_{ijS}^2 $$
$$
=  \( \frac{1+c_{\tau}}{2} \)^2 t_\alpha^2 \frac{  \| \nu' \text{Diag}( \Psi_{jS}^{\hat \varepsilon})^{1/2} \|_1 ^2} { \nu ' \Psi_{jS}^{\hat \varepsilon} \nu} \sum_{i=1}^n \(\sum_{k \in jS} \nu_k x_{ik} \)^2 \hat \varepsilon_{ijS}^2. 
$$
Specializing to the case that $\nu = \eta_{jS}$, and using $\hat \tau_{jS} =  \frac{  \| \nu' \text{Diag}( \Psi_{jS}^{\hat \varepsilon})^{1/2} \|_1} { \sqrt{\nu ' \Psi_{jS}^{\hat \varepsilon} \nu}}  $ gives that $$|\hat Z_{jS} | \leq \frac{1+c_\tau}{2} \hat \tau_{jS} t_\alpha \ \ \text{on} \  \mathscr R \cap \mathscr V.$$

It is therefore sufficient to prove that $\mathscr R$ and $\mathscr V$ have probability $\rightarrow 1$ under $\Pr^{\mathscr X}$.  $\Pr^{\mathscr X}(\mathscr R) \rightarrow 1$ follows immediately from the moderate deviation bounds for self-normalized sums given in \cite{jing:etal}.  For details on the application of this result, see \cite{BellChenChernHans:nonGauss}.

Therefore, it is only left to show that $\Pr^{\mathscr X}(\mathscr V) \rightarrow 1$.  Define $\varepsilon_{ijS} = y_i  - x_{ijS}'\theta_{jS}^{*|x}.$ Furthermore, define $\xi_{ijS}$ through the decomposition $ \varepsilon_{ijS} = \ddot \varepsilon_i + \xi_{ijS}$.  Let $\eps_{jS}$ and $\xi_{jS}$  be the respective stacked versions.  Let $\tilde c_\tau = ((1 + c_{\tau})/2)^2$.
Then 
\begin{align*} \tilde c_\tau \sum_{i=1}^n x_{ik}^2 \hat \varepsilon_{ijS}^2 &= \tilde c_{\tau} \[ \sum_{i=1}^n x_{ik}^2 (\hat \varepsilon_{ijS}^2 -  \varepsilon_{ijS}^2) + \sum_{i=1}^n x_{ik}^2\ddot  \varepsilon_{i}^2 + 2\sum_{i=1}^n x_{ik}^2 \ddot \varepsilon_i \xi_{ijS} + \sum_{i=1}^n x_{ik}^2 \xi_{ijS}^2 \] \\ & \geq \tilde c_{\tau} \[ \sum_{i=1}^n x_{ik}^2 (\hat \varepsilon_{ijS}^2 -  \varepsilon_{ijS}^2)  + \sum_{i=1}^n x_{ik}^2 \ddot \varepsilon_{i}^2 + 2\sum_{i=1}^n x_{ik}^2\ddot \varepsilon_i \xi_{ijS}  \]  \\ 
& = \sum_{i=1}^n x_{ik}^2 \ddot \varepsilon_{i}^2  +  \tilde c_{\tau}   \sum_{i=1}^n x_{ik}^2 (\hat \varepsilon_{ijS}^2 - \varepsilon_{ijS}^2)  + \frac{(\tilde c_{\tau} - 1)}{2} \sum_{i=1}^n x_{ik}^2 \ddot \varepsilon_i^2      \\ 
& \ \ \ \ \ \ +2\tilde c_\tau \sum_{i=1}^n x_{ik}^2\ddot \varepsilon_i \xi_{ijS}  +  \frac{(\tilde c_{\tau} - 1)}{2} \sum_{i=1}^n x_{ik}^2 \ddot \varepsilon_i^2  .    \end{align*}

Define the two events 
$$\mathscr V' = \left \{ \tilde c_{\tau}   \En[ x_{ik}^2 (\hat \varepsilon_{ijS}^2 - \varepsilon_{ijS}^2) ] + \frac{(\tilde c_{\tau} - 1)}{2} \En[x_{ik}^2 \ddot \varepsilon_i^2 ] \geq 0 \ \text{for all} \ j,k \leq p , |S|<K_n \right \}$$
$$ \mathscr V'' = \left \{2\tilde c_\tau\En [ x_{ik}^2\ddot \varepsilon_i \xi_{ijS} ] +  \frac{(\tilde c_{\tau} - 1)}{2} \En [ x_{ik}^2 \ddot \varepsilon_i^2 ] \geq 0    \ \text{for all} \ j,k \leq p , |S|<K_n \right \} $$

Therefore $\mathscr V' \cap \mathscr V''  \Rightarrow \mathscr V$.  

%$\Pr^{\mathscr X}(\mathscr V') \rightarrow 1$ by Chebyshev inquality noting that $\En[x_{ij}^4] = O(1)$ on $\mathscr X_1$ and by boundedness of $\Ep[ \ddot \eps^4 | x]$.  Next turn to $\mathscr V''$.  By boundedness of the conditional variance of $\ddot \varepsilon_i$, 

%$$ \text{var}\( \En[ x_{ik}^2\ddot \varepsilon_i \xi_{ijS} ]  |x\) \leq  O(1) \frac{1}{n} \En[ x_{ik}^4 \xi_{ijS}^2]  $$

%$$ \leq \frac{1}{n} \( \En[ x_{ik}^8 ] \)^{1/2} \(\En [ \xi_{ijS}^4 ]\)^{1/2}  $$

%By definition, $\( \En[ x_{ik}^8 ] \)^{1/2}\leq O(1)$ on $\mathscr X_1$.  
%In addition, $\(\En[\xi_{ijS}^4 ] \)^{1/2} \leq O(1)K_n^4 $ on  $\mathscr X_1 \cap \mathscr X_3$.
%To see this, note $\xi_{jS} = M_{jS}x\theta_0 = \sum_{l=1}^{s_0} M_{jS} x_l\theta_{0,l} = \sum_{l=1}^{s_0} \eta_{l,(jS)} x_{ljS} = \tilde \eta_{jS} x_{S_0 \cup jS}$
%for some new linear combination $\tilde \eta_{jS}$.  Note that $\| \tilde\eta_{jS} \|_1 \leq s_0 O(1)$.
%Then $$ (\tilde \eta_{jS} x_{i,S_0 \cup jS} )^4 \leq \( \| \tilde \eta_{jS} \|_2 \)^4 \( \sum_{l \in S_0 \cup jS} x_{il}^2 \)^2  $$

%$$ \leq \( \sqrt{|S_0 \cup jS|} \| \tilde \eta_{jS} \|_1 \)^4 |S_0 \cup jS| \sum_{l \in S_0 \cup jS} x_{il}^4   $$

%$$\leq (s_0 + K_n)^2 s_0^4 O(1) (s_0 + K_n) \sum_{l \in S_0 \cup jS} x_{il}^4 $$

%Therefore, 

%$$\( \En [\xi_{ijS}^4 ]\)^{1/2}  \leq \sqrt{(s_0 + K_n)^2 s_0^4 O(1) (s_0 + K_n)^2  \max_{k \leq p} \En[ x_{ik}^4 ]} \ \leq O(1)K_n^4$$

%Since $K_n^4/n \rightarrow 0$, this gives $\Pr^{\mathscr X}(\mathscr V'') \rightarrow 1$.  

\

%First, note that on $\mathscr X_4$, $\Pr( \sum_{i=1}^n x_{ik}^2 \eps_i^2 \geq Cn|x) \rightarrow 1$ for some $C>0$ uniformly over $x \in \mathscr X_4$.

%Next, $\Pr( \sum_{i=1}^n x_{ik}^2 (\ddot \eps_i^2 - \eps_i^2 \geq Cn|x) \rightarrow 1$

Note that $\En [x_{ik}^2\ddot \eps_i^2] \geq \frac{1}{2} \En [ x_{ik}^2 \eps_i^2] - \En[ x_{ik}^2 \Ep[\eps_i^\text{a}|x] ] \geq  \frac{1}{2} \En[ x_{ik}^2 \eps_i^2 ]- \max_{i\leq n}  \Ep[\eps_i^{\text{a}2}|x]^{1/2} \En[ x_{ik}^4]^{1/2}$. 
This is bounded below with  $\Pr^{\mathscr X}\rightarrow 1$ by a positive constant independent of $n$.
Therefore, to show that $\Pr^{\mathscr X}(\mathscr V')\rightarrow 1$, $\Pr^{\mathscr X}(\mathscr V'')\rightarrow 1$, it suffices to show 
$ \En [x_{ik}^2 (\hat \varepsilon_{ijS}^2 - \varepsilon_{ijS}^2)]$ 
 and $\En[ x_{ik}^2 \ddot \eps_{i} \xi_{ijS}]$, respectively, are suitably smaller order.  
 
 First consider $ \En [x_{ik}^2 (\hat \varepsilon_{ijS}^2 - \varepsilon_{ijS}^2)]$.  It is convenient to bound the slightly more general sum $\En[x_{ik}x_{il} (\hat \varepsilon_{ijS}^2 - \varepsilon_{ijS}^2)]$, because this will show up again.
$$\En[x_{ik}x_{il}(\hat \varepsilon_{ijS}^2 - \varepsilon_{ijS}^2) ] = 2 \En \[ x_{ik}x_{il} \varepsilon_{ijS}x_{ijS}'(\theta_{jS}^{*|x} - \hat \theta_{jS}) \] +  \En\[x_{ik}x_{il}(x_{ijS}'(\theta_{jS}^{*|x} - \hat \theta_{jS}))^2 \] $$
$$ \leq 2 \| \En [ x_{ik}x_{il}  \varepsilon_{ijS}x_{ijS}' ] \|_2 \| \theta_{jS}^{*|x} - \hat \theta_{jS} \|_2  + \lambda_{\max}\En [ x_{ik}x_{il} x_{ijS}x_{ijS}']  \| \theta_{jS}^{*|x} - \hat \theta_{jS}) \|_2^2 $$
Standard reasoning gives that $ \| \theta_{jS}^{*|x} - \hat \theta_{jS} \|_2 \leq  \varphi_{\min}(K_n)(G)^{-1/2}\sqrt{K_n} \| \En x_{ijS} \varepsilon_{ijS} \|_\infty$.  Therefore, the bound continues
$$ \leq 2 \| \En [ x_{ik}x_{il}  \varepsilon_{ijS}x_{ijS}' ] \|_2\varphi_{\min}(K_n)(G)^{-1/2}\sqrt{K_n} \| \En x_{ijS} \varepsilon_{ijS} \|_\infty $$ $$ + \lambda_{\max}\En [ x_{ik}x_{il} x_{ijS}x_{ijS}']  \varphi_{\min}(K_n)(G)^{-1}K_n \| \En x_{ijS} \varepsilon_{ijS} \|_\infty^2. $$
Note that $\lambda_{\max} \En [ x_{ik}x_{il} x_{ijS}x_{ijS}'] \leq K_n \max_{j\leq p} \En[x_{ij}^4]$.
$$ \leq 2 \| \En [ x_{ik}x_{il} \varepsilon_{ijS}x_{ijS}' ] \|_2\varphi_{\min}(K_n)(G)^{-1/2}\sqrt{K_n} \| \En x_{ijS} \varepsilon_{ijS} \|_\infty $$ $$ +  K_n^2 \max_{j\leq p} \En[x_{ij}^4] \varphi_{\min}(K_n)(G)^{-1} \| \En x_{ijS} \varepsilon_{ijS} \|_\infty^2. $$
An application of Cauchy-Schwarz to the top line gives
$$ \leq 2 \sqrt{K_n} \max_{j}\En [ x_{ik}^4 ]^{1/2} \max_{j,S} \En[\varepsilon_{ijS}^2x_{ij}^2 ]^{1/2}\varphi_{\min}(K_n)(G)^{-1/2}\sqrt{K_n} \| \En x_{ijS} \varepsilon_{ijS} \|_\infty $$ $$ +  K_n^2 \max_{j\leq p} \En[x_{ij}^4] \varphi_{\min}(K_n)(G)^{-1} \| \En x_{ijS} \varepsilon_{ijS} \|_\infty^2. $$
Next, $\| \En x_{ijS} \varepsilon_{ijS} \|_\infty$ and $ \En[\varepsilon_{ijS}^2x_{ij}^2 ]^{1/2}$ are bounded using $\varepsilon_{ijS} =  \eps_i - \Ep[\eps_i|x]+ \xi_{ijS} $.
Note that by construction $\| \En[ x_{ijS} \xi_{ijS}] \|_\infty=0$.  Then
$$\| \En[ x_{ijS} \eps_{ijS} ]\|_\infty \leq \| \En[x_{i} \eps_i ] \|_{\infty}+ \| \En[x_{i} \Ep[\eps^{\text{a}}_i |x ] \|_{\infty}$$
$$\leq \| \En[x_{i} \eps_i ] \|_{\infty}+ \max_{j \leq p}\En[x_{ij}^2]^{1/2} \En[ \Ep[\eps_{i}^{\text{a}}|x]^2]^{1/2}=O( \sqrt{\log p/n} ) $$ with $ \Pr^{\mathscr X} \rightarrow 1$.
Next,  $$\En [\eps_{ijS}^2x_{ij}^2] \leq 3\En [\eps_{i}^2x_{ij}^2]  + 3\En [ \Ep[\eps_{i}^{\text{a}2}|x]x_{ij}^2]  +3 \En [\xi_{ijS}^2x_{ij}^2] $$  
$$ \leq 3\En [\eps_{i}^2x_{ij}^2]  + 3\En [ x_{ij}^2 ] \max_{i \leq n} \Ep[\eps_{i}^{\text{a}2}|x] +3 \En [\xi_{ijS}^4]^{1/2} \En[x_{ij}^4]^{1/2} .$$  
Next, $\(\En[\xi_{ijS}^4 ] \)^{1/2} \leq O(1)s_0^2 $ on  $\mathscr X_1 \cap \mathscr X_3$.  To see this, note $\xi_{jS} = {\mathscr Q}_{jS}x\theta_0 = \sum_{l=1}^{s_0} {\mathscr Q}_{jS} x_l\theta_{0,l} = \sum_{l=1}^{s_0} \eta_{l,(jS)} x_{ljS} = \tilde \eta_{jS} x_{S_0 \cup jS}$
for some new linear combination $\tilde \eta_{jS}$.  Note that $\| \tilde\eta_{jS} \|_1 \leq s_0 O(1)$.
Then  $ \(\En[\xi_{ijS}^4 ] \)^{1/4} \leq \| \tilde \eta_{jS} \|_1 \max_{k \leq p} \En[ x_{ik}^4]^{1/4} $ from which the bound follows.

Next consider $\En[ x_{ik}^2 \ddot \eps_{i}\xi_{ijS}]$.  
Consider two cases.  In Case 1, $$\En[x_{ik}^4\xi_{ijS}^2]^{1/2} \leq  \frac{1}{2 \tilde c_{\tau}}  \frac{(\tilde c_{\tau} - 1)}{2} \frac{\En [ x_{ik}^2 \ddot \varepsilon_i^2 ]}{ \En[\ddot \eps_{i}^2]^{1/2}}.$$  In this case, 
$2 \tilde c_{\tau} \En[x_{ik}^2 \ddot \eps_i \xi_{ijS} ] \leq \En[x_{ik}^4\xi_{ijS}^2]^{1/2} \En[\ddot \eps_{i}^2]^{1/2} \leq \frac{\tilde c_{\tau} - 1}{2} $, and the requirement of $\mathscr V''$ for $k,j,S$ holds.  

For Case 2, suppose the alternative that $\En[x_{ik}^4\x_{ijS}^2] >  \frac{1}{2 \tilde c_{\tau}}  \frac{(\tilde c_{\tau} - 1)}{2} \frac{\En [ x_{ik}^2 \ddot \varepsilon_i^2 ]}{ \En[\ddot \eps_{i}^2]^{1/2}}$ holds.  Then $\Ep[ \En[x_{ik}^4\xi_{ijS}^2 \ddot \eps_i^2] | x]$ is bounded away from zero by conditions on $\Ep[\eps_i^2|x]$ and $\max_{i} |\eps_i^{\text{a}}|$.  In addition,  $\Ep[ \En[|x_{ik}|^6|\xi_{ijS}|^3 |\ddot \eps_i|^3] |x ] \leq \max_{i} \Ep[| \ddot{\eps_i}|^3|x]  \En[ |x_{ik}|^6 |\xi_{ijS}|^3]  \leq O(1) \En[ |x_{ik}|^6 |\xi_{ijS}|^3]  $.  This term is further bounded by 
$ O(1)\En[ x_{ik}^{12}]^{1/2} \En[|\xi_{ijS}|^6]^{1/2}.$
Using the same reasoning as bounding $\En[ \xi_{ijS}^4 ]$ earlier, it follows that $\En[|\xi_{ijS}|^6]^{1/2} = O(1)s_0^3$.  In addition, $\En[ x_{ik}^{12}] = O(1)$.  
As a result, for those $k,j,S$ which fall in Case 2, the self-normalized sum
 $$   = \max_{j,k,S \in \text{Case 2}}  \frac{\sqrt{n} |\En[ x_{ik}^2 \xi_{ijS} \ddot \eps_i] | }{ \sqrt{\En[ x_{ik}^4 \xi_{ijS}^2 \ddot \eps_i^2] } } $$
is $O(\log (p^{K_n}))$ with probability $1 - o(1)$ provided $\sqrt{\log(p^{K_n})} = o(n^{1/6}/(s_0^3)^{1/3})$.  This holds under the assumed rate conditions.  
Then $\max_{j,k,S} |\En[ x_{ik}^2 \xi_{ijS} \ddot \eps_i] |$ is bounded by $\frac{1}{\sqrt{n}}O(\log(p^{K_n}) \max_{j,k,S}  \sqrt{\En[ x_{ik}^4 \xi_{ijS}^2 \ddot \eps_i^2] } $.  Furthermore,     $\En[ x_{ik}^4 \xi_{ijS}^2 \ddot \eps_i^2]  \leq \En[x_{ik}^8 \xi_{ijS}^4]^{1/2} \En[\ddot \eps_{i}^4]^{1/2} \leq (\En[x_{ik}^{12}]^{2/3} \En[\xi_{ijS}^{12}]^{1/3})^{1/2} \En[\ddot \eps_{i}^4]^{1/2} \leq O(1) s_0^2 \En[ \ddot \eps_i^4]^{1/2}$.   Note that $\En[ \ddot \eps_i^4]^{1/2} \leq O(1)$ with $\Pr^{\mathscr X} \rightarrow 1$.  Together, these give that $\max_{j,k,S}\En[x_{ik}^2 \ddot \eps_i \xi_{ijS}] =o(1)$ with $\Pr^{\mathscr X} \rightarrow 1$.  Finally, $\Pr^{\mathscr X}(\mathscr V) \rightarrow 1$.

%Therefore, 

%$$\( \En [\xi_{ijS}^4 ]\)^{1/2}  \leq \sqrt{(s_0 + K_n)^2 s_0^4 O(1) (s_0 + K_n)^2  \max_{k \leq p} \En[ x_{ik}^4 ]} \ \leq O(1)K_n^4$$

%Since $K_n^4/n \rightarrow 0$, this gives $\Pr^{\mathscr X}(\mathscr V'') \rightarrow 1$.  

\subsubsection{Proof of Lemma 14} 

%This step shows that $\Pr^{\mathscr X}(\mathscr B) \rightarrow 1$. 

First, $A_{jS}$ depend only on $x$ and are bounded above and below by constants which do not depend on $n$ on $\mathscr X$ from the assumption on the sparse eigenvalues of $G$.  For bounding $\hat \tau_{jS}$ above and away from zero, since $1 \leq \| \eta_{jS} \|_1, \|\eta_{jS} \|_2 \leq O(1)$ on $\mathscr X$, it is sufficient to show that the eigenvalues of $\Psi_{jS}^{\hat \eps} = \En[x_{ijS}x_{ijS}' \hat \eps_{ijS}^2]$ remain bounded above and away from zero and that the diagonal terms of $\Psi_{jS}^{\hat \eps}$ remain bounded above and away from zero.  Note that by arguments in last section, it was shown that  $\En[x_{ik}x_{il}(\hat \eps_{ijS} - \eps_{ijS})] = O(\sqrt{\log p /n})$ with $\Pr^{\mathscr X} \rightarrow 1$.  Therefore,  
$ \| \En[x_{ijS}x_{ijS}' \hat \eps_{ijS}^2] - \En[x_{ijS}x_{ijS}' \eps_{ijS}^2] \|_{\mathscr F} = O(K_n \sqrt{\log p /n} \ ) $ with $\Pr^{\mathscr X} \rightarrow 1$.
Here, $\mathscr F$ is the Frobenius norm.  By the assumed rate condition, the above quantity therefore vanishes with $\Pr^{\mathscr X} \rightarrow 1$. 

Next, $$\En[ x_{ijS}x_{ijS}'  \eps_{ijS}^2] = \En[ x_{ijS}x_{ijS}'  \eps_{i}^2] + 2 \En[ x_{ijS}x_{ijS}'  \eps_{i}( \xi_{ijS}+ \Ep[\eps_i^{\text{a}}|x] ) ] $$ $$+  \En[ x_{ijS}x_{ijS}'  (\xi_{ijS} + \Ep[\eps_i^{\text{a}}|x])^2]$$

The first term above, $ \En[ x_{ijS}x_{ijS}'  \eps_{i}^2] $, has eigenvalues bounded away from zero for all $j,S$ with $\Pr^{\mathscr X} \rightarrow 1$.  The third term above, $\En[ x_{ijS}x_{ijS}'  (\xi_{ijS} + \Ep[\eps_i^{\text{a}}|x])^2]$ is positive semidefinite by construction.  The second term above has Frobenius norm tending to zero for all $j,S$ with $\Pr^{\mathscr X} \rightarrow 1$.  This, in conjunction with the fact that the eigenvalues of $\En[ x_{ijS}x_{ijS}' \hat \eps_{ijS}]$ are bounded above and away from zero with $\Pr^{\mathscr X} \rightarrow 1$ shows that the eigenvalues of $\Psi_{jS}^{\hat \eps} = \En[x_{ijS}x_{ijS}' \hat \eps_{ijS}^2]$ are bounded above and away from zero with $\Pr^{\mathscr X} \rightarrow 1$.  Finally, for bounding $\hat V_{jS}$, it is sufficient to show that $\max_{k \leq p}\En[ \eps_i^2 (\eta_{jS}'x_{ijS})^2]$ be bounded above.  This follows immediately from $\Ep[ \eps_i^4 | x]$ being uniformly bounded and $\max_{j,S} \| \eta_{jS }\|_1 = O(1)$ and $\max_{k \leq p}\En[x_{ik}^4] = O(1)$.  These imply that $\Pr^{\mathscr X}(\mathscr B) \rightarrow 1$.

\section*{References}
A. Belloni, D. Chen, V. Chernozhukov, and C. Hansen. Sparse models and methods for optimal 

instruments with an application to eminent domain. \textit{Econometrica}, 80:2369--2429, 2012.

\noindent
B. Y. Jing, Q. M. Shao, and Q. Wang. Self-normalized Cram\`er-type large deviations for inde-

pendent random variables. \textit{Annals of Probability}, 31(4):2167--2215, 2003.

\end{document}